\documentclass[]{arXivclass} %

\DoubleSpacedXI %

\usepackage{endnotes}
\let\footnote=\endnote
\let\enotesize=\normalsize
\def\notesname{Endnotes}%
\def\makeenmark{\hbox to1.275em{\theenmark.\enskip\hss}}
\def\enoteformat{\rightskip0pt\leftskip0pt\parindent=1.275em
  \leavevmode\llap{\makeenmark}}

\usepackage{natbib}
 \bibpunct[, ]{(}{)}{,}{a}{}{,}%
 \def\bibfont{\small}%

\usepackage{multibib}
\newcites{EC}{References}

\TheoremsNumberedThrough     %
\ECRepeatTheorems

\EquationsNumberedThrough    %

\usepackage[hypertexnames=false]{hyperref}
\usepackage{mathtools, nccmath}
\usepackage{braket} %
\usepackage{siunitx}
\usepackage{enumitem}
\setlist[enumerate]{align=left}
\usepackage{amsmath,amssymb}
\usepackage{mathabx}
\usepackage{cleveref}
\usepackage{xcolor}
\usepackage[short]{optidef}
\DeclareDocumentEnvironment{miniinf}{D||{\defaultProblemFormat} O{\defaultConstraintFormat} D<>{} m m m m}{%
\ifthenelse{\equal{#3}{b}}{%
    \ifthenelse{\equal{#1}{s}}%
    {\setFormatShort{inf}{#4}\BaseMiniStar{#2}{#4}{#5}{#7}{inf}{#3}}%
    {\setFormatLong{infimum}{#4}\BaseMiniStar{#2}{#4}{#5}{#7}{infimum}{#3}}%
}{%
    \ifthenelse{\equal{#1}{s}}%
    {\setFormatShort{inf}{#4}\BaseMini{#2}{#4}{#5}{#6}{#7}{inf}}%
    {\setFormatLong{infimum}{#4}\BaseMini{#2}{#4}{#5}{#6}{#7}{infimum}}%
}%
}%
{\endBaseMini\toggletrue{bodyCon}}
                        
\usepackage{booktabs}
\usepackage{graphicx}
\usepackage{subfig} %
\usepackage{multirow}
\usepackage[tworuled,linesnumbered]{algorithm2e}

\newcommand{\RR}{\mathbb R}
\newcommand{\EE}{\mathbb E}

\newcommand{\INT}[1]{\mathrm{int}( #1 )}
\newcommand{\CL}[1]{\mathrm{cl}( #1 )}
\newcommand{\RI}[1]{\mathrm{rint}( #1 )}

\newcommand{\norm}[1]{\left\lVert#1\right\rVert}
\newcommand{\abs}[1]{\left\lvert#1\right\rvert}

\newcommand{\qedsymbol}{\hfill\ensuremath{\Box}}

\allowdisplaybreaks

\newcommand{\cost}[3]{\varphi({#1}, {#3})}

\newcommand{\HK}{{\mathcal{H}}_{\mu}^K}
\newcommand{\barHK}{{\bar{\mathcal{H}}}_{\mu}^K}

\newcommand{\cC}{\mathcal C}

\newcommand{\cE}{\mathcal E}
\newcommand{\cF}{\mathcal F}
\newcommand{\cG}{\mathcal G}
\newcommand{\cH}{\mathcal H}
\newcommand{\cI}{\mathcal I}
\newcommand{\bI}{\mathbf{I}}
\newcommand{\wtbI}{\widetilde{\mathbf I}}
\newcommand{\wtI}{\widetilde{I}}
\newcommand{\wtxi}{\widetilde{\xi}}

\newcommand{\cK}{\mathcal K}

\newcommand{\cN}{\mathcal N}

\newcommand{\cO}{\mathcal O}
\newcommand{\cP}{\mathcal P}
\newcommand{\cR}{\mathcal R}
\newcommand{\cS}{\mathcal S}
\newcommand{\cT}{\mathcal T}
\newcommand{\cU}{\mathcal U}

\newcommand{\cY}{\mathcal Y}
\newcommand{\cX}{\mathcal X}
\newcommand{\cZ}{\mathcal Z}

\newcommand{\Xcube}{\bar{\mathcal X}}
\newcommand{\f}{\bar{f}_0}

\newcommand{\RS}{\widehat{R}_{\Xi}}

\newcommand{\wtheta}{\bar{\theta}}
\newcommand{\IJen}{\wtbI^{\varepsilon}_{\nu}}

\newcommand{\bthetan}{{\theta}^{\nu + \frac{1}{2}}}
\newcommand{\distSd}[1]{\mathrm{dist}(#1, \cS^d)}
\newcommand{\dist}[2]{\mathrm{dist}(#1, #2)}

\newcommand{\dir}{d \in \cT(\theta;\Theta),\norm{d}\leq 1}

\newcommand{\HatFmin}{\widehat{F}^\varepsilon_{\min}}
\newcommand{\RKHSbk}{$\text{RKHS-DR}_\text{BK}$ }
\newcommand{\RKHSnp}{$\text{RKHS-DR}_\text{NP}$ }
\newcommand{\redlog}[1]{\textcolor[rgb]{ .753,  0,  0}{\textbf{\underline{#1}}}}
\newcommand{\bluelog}[1]{\textcolor[rgb]{ 0,  .439,  .753}{\textbf{#1}}}
\newcommand{\nonlinearity}{\omega}

\begin{document}

\RUNAUTHOR{Zhang, Liu, and Zhao}

\RUNTITLE{Data-driven PADR for SP with Covariate Information} 

\newcommand{\myfulltitle}{Data-driven Piecewise Affine Decision Rules for Stochastic Programming with Covariate Information}%
\TITLE{\myfulltitle}

\ARTICLEAUTHORS{%
\AUTHOR{Yiyang Zhang}
\AFF{Department of Industrial Engineering, Tsinghua University, Beijing 100084, China, \EMAIL{zyy21@mails.tsinghua.edu.cn}} %
\AUTHOR{Junyi Liu}
\AFF{Department of Industrial Engineering, Tsinghua University, Beijing 100084, China, \EMAIL{junyiliu@tsinghua.edu.cn}}
\AUTHOR{Xiaobo Zhao}
\AFF{Department of Industrial Engineering, Tsinghua University, Beijing 100084, China, \EMAIL{xbzhao@tsinghua.edu.cn}}
} %

\ABSTRACT{%
Focusing on stochastic programming (SP) with covariate information, this paper proposes an empirical risk minimization (ERM) method embedded within a nonconvex piecewise affine decision rule (PADR), which aims to learn the direct mapping from features to optimal decisions. We establish the nonasymptotic consistency result of our PADR-based ERM model for unconstrained problems, which illustrates the role of piece number in balancing the trade-off between the approximation and estimation errors. To solve the nonconvex and nondifferentiable ERM problem, we develop an enhanced stochastic majorization-minimization (ESMM) algorithm and establish the convergence to (composite strong) directional stationarity along with convergence rate analysis. We show that the proposed PADR-based ERM method applies to a broad class of nonconvex SP problems with theoretical consistency guarantees and computational tractability. Our numerical study demonstrates the effectiveness of PADR-based ERM methods with comparisons to state-of-the-art approaches under various settings, with significantly lower costs, less computation time, and robustness to feature dimensions and nonlinearity of the underlying dependency. 
}%

\KEYWORDS{stochastic programming, covariate information, decision rule, piecewise affine function, nonconvex and nondifferentiable optimization} 
\SUBJECTCLASS{programming: nondifferentiable; nonlinear; stochastic}

\maketitle
\vspace{-1\baselineskip}
\section{Introduction}\label{sec:Introduction}

With the fast-growing real-world applications and availability of large datasets, 
the field of stochastic programming (SP) has flourished with many theoretical and algorithmic advances \citep{birge2011introduction, shapiro2021lectures}, in which the empirical probability distribution based on historical scenarios plays a central role in approximating the true distribution. However, it is common in practice that some covariate information (also referred to as contextual information or features) that may affect the random variable could be observed before the decision is made. In such situations, the empirical probability distribution and classical SP methodologies ignoring the covariate information  lead to the poor performance. Given the covariate information $X = x\in \RR^p$, an unconstrained SP is formulated as follows,
\begin{equation}\label{eq:intro:contextalSP}
    \inf_{z\in \RR^d} \quad \EE_Y\left[\cost{z}{X}{Y} \mid X = x\right],
\end{equation}
where $z \in \RR^d$ is a decision vector, $\varphi: \RR^d \times \cY \to \RR$ is a cost function, and $Y$ is a random variable supported by $\cY$ whose distribution is related to a feature vector $X$ supported by $\cX \subseteq \RR^p$.
A simple instance is an inventory planning problem, in which $z$ represents order decisions for $d$ products, $\varphi$ represents the inventory cost, $Y$ represents uncertain demands, and $X$ represents covariate information such as seasonality, promotion, weather, past orders, etc. With historical scenarios of $(X, Y)$, a classical sample average approximation (SAA) model disregards the dependency of $Y$ on $X$ and thus may lead to inconsistent decisions; such deficiency is also observed in numerical results of \citet{ban2019big} and \citet{bertsimas2022data}. 
Motivated by the latent dependency of the optimal solution to~\eqref{eq:intro:contextalSP} on the feature vector, this paper aims to construct a data-driven decision rule that directly maps features to optimal decisions so as to approximately solve the SP problem~\eqref{eq:intro:contextalSP} for any possible observation of the feature.

\subsection{Problem setting} 
Solving problem~\eqref{eq:intro:contextalSP} for any realization of $X$ can be generalized as seeking an optimal decision rule (DR) to the following problem
\begin{equation}
    \label{eq:intro:contextalSP:Full-E}
    \inf_{f \in \mathcal F} \quad R(f) \coloneqq \EE_{X, Y}\left[\cost{f(X)}{X}{Y}\right],
\end{equation}
where $\mathcal F$ contains all measurable functions 
mapping from the feature space $\mathcal X$ to the decision space~$\RR^d$.
By the interchangeability principle~\citep[Theorem 7.80]{shapiro2021lectures}, if a mapping $f^* \in \cF$ is optimal to problem~\eqref{eq:intro:contextalSP:Full-E}, then the decision $f^*(x)$ is optimal to~\eqref{eq:intro:contextalSP} for almost every $x \in \cX$. 
In general, problem~\eqref{eq:intro:contextalSP:Full-E} is intractable due to the infinite-dimensional function space and the unknown joint distribution of $X$ and $Y$. 
It is common to consider an approximation method which constructs an empirical risk minimization (ERM) problem based on a dataset $\Xi \coloneqq \{(x^s, y^s)\}_{s=1}^n$ with the decision rule $f$ restricted to a parametric DR class $\mathcal H$ (also known as {hypothesis class}):
\begin{equation}
    \label{eq:intro:ERM-H}
    \inf_{f \in \cH} \quad \widehat{R}_{\Xi}(f) \coloneqq 
    \frac{1}{n}\sum_{s=1}^n \cost{f(x^s)}{x^s}{y^s}.
\end{equation}

In this paper, we adopt the above DR-based ERM framework, recognizing that the efficacy of such a framework is influenced by the trade-off between approximation accuracy of $\mathcal H$ to the optimal DR and computational complexity of solving~\eqref{eq:intro:ERM-H}.
For example, problem~\eqref{eq:intro:ERM-H} with a convex cost function and linear decision rules (LDR) can be solved efficiently \citep{beutel2012safety, ban2019big}, while the LDR has difficulty describing the potential nonlinear relation between decisions and features.
On the other side, a DR class consisting of neural networks could have a strong approximation capability~\citep{qi2022practical}; nevertheless, the resulting nonconvex problem~\eqref{eq:intro:ERM-H} with a multi-layer composite structure is difficult to solve with theoretical optimality guarantees, and the model may face overfitting issues when the dataset is not large enough.

Therefore, this paper aims to address a critical question: 
\emph{what could be an amenable hypothesis class for the DR-based ERM framework so as to maintain both high reliability and computational efficiency for solving SP with covariate information?} 
To this end, we employ a piecewise affine decision rule (PADR) in the ERM problem~\eqref{eq:intro:ERM-H} that is broader than LDR but simpler than many general nonlinear DRs.
According to \citet{scholtes2012introduction}, any continuous piecewise affine (PA) function $f:\cX \rightarrow \RR$ can be expressed as the difference of max-affine functions:
\begin{equation}\label{eq:intro:continuous PA function}
    f(x,\theta) = \max_{k \in [K]} \left\{(\alpha^k)^\top x + a^k \right\}
                    - \max_{k\in [K]} \left\{ (\beta^k)^\top x + b^k \right\},
\end{equation}
where $\theta \coloneqq \left(\alpha^k, a^k, \beta^k, b^k\right)_{k=1}^K \in \RR^{2 K(p+1)}$ and $[K]\coloneqq \{1,\dots, K\}$. The PADR of such a form is beneficial because it requires minimal parameter tuning in multi-dimensional cases without the need to specify breakpoints for piecewise
linear parts in each dimension, and it can asymptotically approximate any continuous function according to \cite{royset2020approximations}.  
Throughout the theoretical analysis in this paper, we assume that the numbers of pieces in the two max-affine components of~\eqref{eq:intro:continuous PA function} are equal to simplify the exposition, which could be easily extended to the case with different numbers of pieces.

To address the aforementioned question, we intend to tackle two critical issues of the PADR-based ERM model for solving the unconstrained decision-making problem~\eqref{eq:intro:contextalSP} in this paper: the consistency guarantees of the PADR-based ERM problem  
and the computational method for efficiently solving a large-scale, nonconvex, and nondifferentiable composite optimization problem. It is more challenging to extend the data-driven decision rule methods for constrained decision-making problems. Motivated by the scenario approximation method in \cite{calafiore2005uncertain}, a straightforward approach is to impose the strengthened feasibility requirement over the observed scenarios in learning the data-driven decision rule. In 
Appendix \ref{section:Extension:constrained}, we illustrate that the obtained decision rule at a new scenario is feasible with high probability and is asymptotically consistent with the original SP problem, and we also discuss several special cases where the feasibility of the decision rule can be guaranteed.

In the rest of this section, we give a brief literature review and present the distinguished features of the proposed PADR method over other related work. We will also summarize the major contributions of the present paper.

\subsection{Literature review}

There are two main schemes of data-driven approaches in the literature of SP dealing with the latent dependency on covariate information: the sequential and integrated schemes. 
In the sequential scheme, also referred to as the predict-then-optimize (PO) scheme \citep{sen2018learning, kannan2022data}, one first predicts uncertain parameters and then solves the optimization problem composite with the prediction model. 
Several works such as \citet{kao2009directed, donti2017task}, \citet{elmachtoub2022smart}, and \citet{qi2021integrated}
propose modifications of the loss function by considering the cost error of the downstream optimization problem. %
Since those modified loss functions involve optimal solutions to the decision-making problems that are parameterized by prediction models, such frameworks usually lead to bi-level nonconvex optimization problems and thus encounter substantial computational challenges in dealing with nonlinear constrained programs.

Approaches in the integrated scheme yield decisions via one-step data-driven optimization.
Prescriptive methods by \citet{bertsimas2020predictive} and \citet{kallus2022stochastic} approximate the conditional expected cost in~\eqref{eq:intro:contextalSP} through reweighing empirical costs. 
\citet{elmachtoub2023estimate} point out that, when the model class is well-specified and there is sufficient data, methods that directly estimate the conditional distribution outperform those that consider the downstream optimization problem; in a misspecified scenario, the reverse holds.

However, as claimed by \citet{bertsimas2022data}, the prescriptive methods are highly affected by the data with features close to the current observation and thus suffer from the curse of dimensionality. 
The DR-based ERM framework presented in~\eqref{eq:intro:ERM-H}, which is another one-step approach in the integrated scheme, is relieved from the aforementioned issue due to its functional approximation in nature.
For constrained SP with general nonlinear constraints, the current study of decision rule methods could only guarantee feasibility with high probability, and the predict-then-optimize framework may be more suitable for solving the contextual problems when feasibility is essential to decision-makers. 
{We refer interested readers to \citet{qi2022integrating} and \citet{sadana2025survey} for comprehensive reviews on recent advances in both two-stage and integrated approaches in contextual SP.}

Several works on DR-based ERM methods are most relevant to this paper.
\citet{ban2019big} employ LDR for the Newsvendor problem with nonasymptotic optimality guarantees under the assumption of the linear demand model. 
Besides the classical LDR, several generalized LDRs are proposed for standard SP problems which lift random vectors into a higher-dimensional probability space to describe nonlinear relations \citep{chen2008linear, bampou2011scenario, georghiou2015generalized}.
However, they tend to have unsatisfactory performance without cautious prior designs especially when the random vector is high-dimensional; moreover, they usually lack theoretical guarantees on approximation accuracy except for the posterior estimation in \citet{kuhn2011primal}. 

\citet{notz2021prescriptive} and \citet{bertsimas2022data} study decision rules in reproducing kernel Hilbert spaces (RKHS), which we refer to as RKHS-DR, and provide asymptotic optimality guarantees for convex SP problems.
{
With the rise of deep neural networks (NN), several studies apply NN-based decision rules (NN-DR) in practical applications, such as inventory management~\citep{zhang2020deep, oroojlooyjadid2020applying, qi2023practical} and portfolio optimization~\citep{zhang2020deep}. NN-DRs are powerful in approximating the underlying relationship between decisions and covariate information; however, they may face numerical challenges for complex decision-making problems and also lack theoretical consistency guarantees.} %
Compared to the above existing DR methods, given the asymptotic approximation of PA functions to any continuous function, we show that the PADR method is applicable to a class of more general decision-making problems with nonconvex cost function with both numerical efficiency and theoretical consistency.

Regarding solving the PADR-based ERM problem with potentially large sample size $n$,  it is prohibitively challenging to obtain a global minimizer to such a coupled nonconvex and nondifferentiable optimization problem.  As provided by \citet[Section 2.4]{bertsimas2015design}, the PADR-based ERM problem can be formulated as a mixed-integer nonlinear program using Big-M method with additional $2 n d K $ binary variables to describe the PADR of interest. For a PADR-based ERM problem of the modest size with data size $n=1000$, the obtained mixed-integer nonlinear program could have $10^4 \sim 10^5$ binary variables and state-of-the-art solvers find it hard to solve within a few hours to achieve an approximately optimal solution. 

Instead of seeking global optimality, majorization-minimization (MM) methods provide a unified framework that is practically efficient and theoretically convergent for solving the PADR-based ERM problem with the finite-sum composite objective.  MM methods iteratively construct convex
subproblems with an upper surrogate objective function for obtaining ``stationary" solutions.  According to variational analysis by \cite{rockafellar2009variational}, there are many kinds
of stationary concepts which differ in the ability to characterize necessary optimality conditions
for local/global optimality. As pointed out in \citet[Chapter 7]{cui2021modern}, MM algorithms may converge to different types of stationary solutions depending on the properties of surrogate and objective functions. Algorithms with convergence to a stronger stationarity are favorable not only
theoretically but also computationally because the obtained solution from numerical implementation with random initialization has the greater potential to be a global optimum. 

Difference-of-convex (DC) algorithms, which is a special MM method, have been extensively studied over the last several decades \citep{le2018dc} and can compute a d(irectional)-stationary point for DC programs with a differentiable concave component. When the cost function $\varphi(\cdot, y)$ is convex or piecewise affine, the PADR-based ERM problem can be formulated as a nondifferentiable DC program for which DC algorithms including \citet{tao1997convex} and \citet{lipp2016variations} are applicable.  %
However, due to the nondifferentiable concave component, DC algorithms for the ERM problem converge to a critical point whose defining condition is dependent on the DC decomposition form and is weaker than d(irectional)-stationarity. For nondifferentiable DC or composite DC problems, \citet{pang2017computing} and \citet{cui2018composite} propose an enhanced surrogation technique that iteratively solves convex subproblems that are constructed based on $\varepsilon$-active pieces of the concave component. However, the enhanced MM algorithm for large-scale ERM problems requires iteratively solving finite-sum convex subproblems over all the data, which raises the fundamental issue of computational efficiency. Sampling-based algorithms for stochastic DC programs have been recently studied in \citet{nitanda2017stochastic, xu2019stochastic,le2022stochastic} addressing smooth or partially smooth objectives, and in \citet{liu2022solving} addressing nonconvex and nonsmooth compound SP to compute critical solutions. There is still a lack of a comprehensive study on whether the $\varepsilon$-active set technique is amenable for sampling-based methods to enhance the convergence properties for nondifferentiable stochastic DC programs and PADR-based ERM problems.  %
 
\subsection{Contributions}

{
Our work contributes to the literature on contextual SP along three major dimensions:

1. \emph{Modeling aspect.}
We propose a PADR-based ERM method for a general class of SP problems with covariate information, where the cost function could be nonconvex and nondifferentiable.
PADR provides a compact structure to approximate non-parametric decision rules.
Though PA hypothesis class is not new in classical SP, we provide the first quantification of the approximation error and the excess risk bound of the data-driven PADR model for general nonconvex contextual SP problems.
The excess risk bound illustrates the role of the piece numbers $K_1$ and $K_2$ in balancing the approximation and estimation errors, which is further evidenced by numerical experiments. 
In contrast, previous DR methods for contextual SP provide asymptotic guarantees only limited to convex or linear settings. %
The PADR method thereby provides a rigorous framework to solve nonconvex contextual SP problems with theoretical guarantees.

2. \emph{Computational aspect.} By exploiting the compact max-affine structure of PADR, we develop the ESMM algorithm to solve the highly nonconvex and nonsmooth PADR-based ERM problem, with the natural combination of sequential sampling and the enhanced technique over the $\varepsilon$-active set within the majorization-minimization framework. 
We analyze the asymptotic convergence of ESMM in terms of composite strong d-stationarity which is strictly stronger than d-stationarity, the sharpest type of first-order stationary solutions for general nonconvex optimization.
Furthermore, with the design of a suitable sampling scheme, the ESMM method significantly improves computational efficiency compared to existing MM approaches. 
Numerical results show that the ESMM algorithm for PADR scales well to large problem sizes (e.g., up to 80 pieces in PADR) and could be up to 4 times faster than the basic MM method.

3. \emph{Empirical aspect.} Numerical experiments with comparisons to a variety of prominent methods including NN-DR and RKHS-DR (see Table~\ref{tab:exp:methods}) show that PADR requires fewer data points to achieve satisfactory performance, and thus it is suitable for cases with insufficient data sets. 
Additionally, it yields superior performance in terms of the test costs with both synthetic and real data sets in various instances including the Newsvendor problem and the product placement problem with feature dimension less than 20. 
Finally, PADR is more robust to the level of nonlinearity in demand models and consistently obtains satisfactory performance with highly nonlinear demand models, compared to all other benchmarks.
}

The rest of this paper is organized as follows.
In Section~\ref{sec:Model}, we provide the theoretical consistency result of the PADR-based ERM model. 
We develop the ESMM algorithm for learning PADR and present the convergence analysis in Section~\ref{section:algorithm}. %
In Section~\ref{section:experiments} we report numerical results. We end with conclusions in Section~\ref{sec:Conclusion}.

\section{Consistency Guarantees for PADR-based ERM}\label{sec:Model}

In this section, we aim to derive the consistency guarantees for the PADR-based ERM problem. Considering problem~\eqref{eq:intro:ERM-H} equipped with a class of continuous PA functions with each max-affine component having $K$ pieces and parameters $\{\theta_\iota\}_{\iota \in [d]}$ satisfying $\norm{\theta_\iota}_\infty \leq \mu$ for $\iota \in [d]$, we obtain the following ERM problem:
\begin{equation}
    \label{eq:intro:ERM-PADR}
    \min_{\norm{\theta_{\iota}}_\infty \leq \mu, \forall \iota \in [d]} \quad  
    \frac{1}{n}\sum_{s=1}^n \varphi\left( \Big( \max_{k\in [K]} \left\{ (\alpha^k_{\iota})^\top x^s + a^k_{\iota} \right\} 
          - \max_{k\in [K]} \left\{ (\beta^k_{\iota})^\top x^s + b^k_{\iota} \right\} \Big)_{\iota=1}^d, y^s \right).
\end{equation}

Excess risk is commonly used to quantify the performance of the ERM model in machine learning, which can be decomposed as the approximation error for the hypothesis class and the estimation error for the empirical risk. Intuitively, when the number of pieces in the PADR class increases, the decision rule class becomes more complex, which is expected to reduce the approximation error while increasing the estimation error from the empirical risk.
In this section, assuming the existence of an optimal decision rule that is Lipschitz continuous, we provide consistency guarantees of the PADR-based ERM model~\eqref{eq:intro:ERM-PADR} by quantifying the approximation and estimation errors, which sheds light upon the role of the number of pieces and sample size in balancing the two errors in the excess risk. To the best of our knowledge, the approximation error bound of the PA function class established in the succeeding section is the first result for Lipschitz continuous functions in the literature.  At the end of this section, we discuss the excess risk in practical implementation and compare our theoretical results with other prominent ERM models employing the LDR and RKHS-DR.

\subsection{Approximation error of piecewise affine function class}\label{subsec:the piecewise affine decision rule}

The approximation errors of the PA function class to convex and DC functions are quantified in \citet{balazs2015near} and \citet{siahkamari2020piecewise} respectively. PA functions that are of ``nonparametric" nature are known to asymptotically approximate Lipschitz continuous functions \citep{royset2020approximations}; nevertheless, the approximation error of PA functions with finite pieces remains unsolved. For the class of Lipschitz continuous functions that allow a more general setting for SP problems of interest, we provide the approximation error bound of the PA function class, which to the best of our knowledge is the first quantification result for such a general case. 

With $L > 0$ and $M > 0$, let
\[
    \mathcal F_{L, M} \coloneqq 
    \left\{
        f:\cX \rightarrow \RR^d 
        \mid
        \norm{f}_\infty \leq M,\,
        \text{$f$ is $L$-Lipschitz continuous}
    \right\},
\]
where $\norm{f}_\infty \coloneqq  \sup_{x\in \mathcal X} \max_{\iota \in [d]}|f_\iota(x)|$. 
We define the PADR class $\HK$ as 
\begin{equation}
\label{eq:intro:PA hypothesis class}
  \mathcal H^K_{\mu} \coloneqq 
  \left\{
      \left(f_\iota(\,\cdot\,, \theta_\iota):\mathcal X \rightarrow \RR\right)_{\iota=1}^d
      \ \middle\vert \ 
      \begin{aligned}
          &f_\iota(x, \theta_\iota) 
          = \max_{k\in [K]} \left\{ (\alpha^k_{\iota})^\top x + a^k_{\iota} \right\} 
          - \max_{k\in [K]} \left\{ (\beta^k_{\iota})^\top x + b^k_{\iota} \right\},\\
          &\theta_{\iota} \coloneqq \left(\alpha_{\iota}^k, a_{\iota}^k, \beta_{\iota}^k, b_{\iota}^k\right)_{k=1}^K \in \RR^{2 K(p+1)}, \norm{\theta_\iota}_{\infty} \leq \mu,\, \forall\,\iota\in [d]
      \end{aligned} \        
  \right\}.
\end{equation}
When $K=1$, $\cH^1_{\mu}$ reduces to the LDR class.  
The following proposition provides the bound for the approximation error of the PA function class to Lipschitz continuous functions.

\begin{proposition}[Approximation error bound of $\HK$]
    \label{prop:model:universal approximation of PADR}
    Suppose that $\mathcal X \subseteq \RR^p$ is a compact set with $\max_{x \in \cX}\norm{x}_\infty = \bar{X}$.
    Then given $L_0 > 0$ and $M_0 > 0$, for any $f_0 \in \mathcal F_{L_0, M_0}$ and $\HK$ defined in~\eqref{eq:intro:PA hypothesis class} with $K \in \mathbb N$ 
    and $\mu \geq \max\left\{\frac{1}{2}L_0 K^{1/p}, \frac{1}{4} p L_0 \bar{X} K^{1/p} + \frac{1}{2}M_0\right\}$, we have
    \[
        \min_{f\in \mathcal H_{\mu}^K} \norm{f - f_0}_{\infty} 
        \leq 2(p^{1/2} + 3)p^{1/2} L_0 \bar{X} K^{-1/p}.
    \]
\end{proposition}

The idea is that for any target function $f_0 \in \mathcal F_{L_0, M_0} $ we construct an interpolation PA function motivated by \citet{siahkamari2020piecewise}. The difference is that with the absence of convexity, the parameters that define the interpolation functions are carefully chosen to ensure that the Lipschitz continuity modulus of the interpolation function is independent of the approximation accuracy. The constructive proof in detail is given in Online Appendix~\ref{proof:prop:universal approximation}. From Proposition~\ref{prop:model:universal approximation of PADR}, the approximation error for Lipschitz continuous functions is bounded by $\cO(K^{-1/p})$,  which is weaker than the approximation error $\cO(K^{-2/p})$  in approximating convex and DC functions by \citet{balazs2015near} and \citet{siahkamari2020piecewise}. However, it applies to general PA-based statistical learning problems in which the underlying true model is not necessarily convex or DC. Moreover, it renders the succeeding nonasymptotic consistency result of the PADR-based ERM model without any parametric or convex assumptions on the optimal decision rule or latent dependency.

\subsection{Nonasymptotic consistency of PADR-based ERM}\label{subsec:nonasymptotic consistency of ERM with PADR}

We impose the following assumptions for the nonasymptotic consistency analysis:
\begin{enumerate}[label=(A\arabic*), series=A]
    \item\label{assum:A1:iid dataset S}  Samples in $\Xi  = \{(x^s, y^s)\}_{s=1}^n$ of $(X,Y)$ are independently and identically distributed.
    \item\label{assum:A2:compactness of X} \label{assum:A2:Lipschitz of f} 
    The support set $\mathcal X \subseteq \RR^p$ of $X$ is a compact set with $\bar{X} \coloneqq \max_{x\in\mathcal X} \norm{x}_\infty $.
    \item \label{assum:A3:Lipschitz of varphi} For some $L_\varphi>0$, $C_\varphi > 0$, $z_0 \in \RR^d$, and any $y \in \mathcal Y$, the cost function $\cost{\cdot}{x}{y}$ is Lipschitz continuous on $\RR^d$ with modulus $L_\varphi$ and satisfies $\abs{\varphi(z_0,y)} \leq C_\varphi$. %
    \item\label{assum:A4:opt lipschitz solution f*} There exists a Lipschitz continuous function $f^*: \cX \rightarrow \RR^d$ that is optimal to problem~\eqref{eq:intro:contextalSP:Full-E}.
\end{enumerate}
The boundedness of $\cX$ in Assumption~\ref{assum:A2:compactness of X} is common in practical applications.
{Assumption~\ref{assum:A4:opt lipschitz solution f*}  may appear to be strict, however, in Online Appendix~\ref{EC:sec:examples-ensuring-Lipschitz-assumption}, we show that it applies to a broad class of problems. We also provide examples of problems that violate this assumption to further clarify its scope of applicability.}
Under Assumption~\ref{assum:A2:compactness of X}, the optimal DR $f^*$ in Assumption~\ref{assum:A4:opt lipschitz solution f*} is bounded on $\cX$, and thus $f^* \in \cF_{L_z, M_z}$ for some $L_z > 0$ and $M_z > 0$.

Let $f^*_{\Xi, \HK}$ denote the optimal solution to the PADR-based ERM model~\eqref{eq:intro:ERM-PADR}.
The consistency of the PADR-based ERM model is evaluated by the excess risk 
$R(f^*_{\Xi, \HK}) - \underset{f \in \cF}{\inf} \, R(f)$,
which is the gap between the expected cost of the optimal empirical PADR and the minimal cost.
With $\min_{f\in \HK}R(f)$ as an intermediate problem, the excess risk can be decomposed into the approximation and estimation errors as follows,
\begin{equation}\label{eq:model-3-excess-risk:excess-risk-decomposition-1}
    R(f^*_{\Xi, \HK}) - \inf_{f \in \cF} R(f)
    = \underbrace{R(f^*_{\Xi, \HK}) - \min_{f\in \HK}R(f)}_{\text{estimation error}} 
      + \underbrace{\min_{f\in \HK}R(f) - \inf_{f\in \cF} R(f)}_{\text{approximation error}}.
\end{equation}
The approximation error, measuring the gap between the minimal expected costs over the PADR class and the class of all measurable functions, can thus be bounded based on the approximation error bound of PADR class provided in Proposition~\ref{prop:model:universal approximation of PADR}. %
We provide the uniform generalization bound based on the Rademacher complexity theory in Proposition~\ref{prop:model:Uniform generalization error bound} which is used to bound the estimation error.
Combining the approximation and estimation error bounds, we thus obtain the excess risk bound of PADR-based ERM in Theorem~\ref{thm:excess-wo-constr}.  The proof of the following two results are provided in Online Appendix~\ref{EC:subsec:proofs for model:nonasymptotic consistency}.
Note that in the present paper the notation $\norm{\cdot}$ without a subscript refers to the Euclidean norm.
 
\begin{proposition}[Uniform generalization bound of PADR-based ERM]\label{prop:model:Uniform generalization error bound}
Under Assumptions~\ref{assum:A1:iid dataset S}-\ref{assum:A3:Lipschitz of varphi}, with $M_{p,d}(\mu) \coloneqq C_\varphi + L_\varphi\left(2\sqrt{d}\mu(p\bar{X} + 1) + \norm{z_0}\right)$, the following holds with probability at least $1-\delta$ over the draw of $\Xi$,
\begin{align*}
    \max_{f\in \HK}\abs{R(f) - \widehat{R}_{\Xi}(f)} 
        \leq &\  \frac{8\mu L_\varphi\sqrt{2dK(p+1)(p\bar{X}^2 + 1)}}{n} + 2M_{p,d}(\mu)\sqrt{\frac{2dK(p+1)\log n}{n}} \\
        &\ + \sqrt{2} M_{p,d}(\mu) \sqrt{\frac{\log(2 / \delta)}{n}}.
\end{align*}
\end{proposition}
\begin{theorem}[Excess risk bound of PADR-based ERM]
    \label{thm:excess-wo-constr}
    Under Assumptions~\ref{assum:A1:iid dataset S}-\ref{assum:A4:opt lipschitz solution f*}, with $\mu \geq \max\left\{\frac{1}{2}L_z K^{1/p}, \frac{1}{4} p L_z \bar{X} K^{1/p} + \frac{1}{2}M_z\right\}$ and $M_{p,d}(\mu) \coloneqq C_\varphi + L_\varphi\left(2\sqrt{d}\mu(p\bar{X} + 1) + \norm{z_0}\right)$, with probability at least $1-\delta$ over the draw of $\Xi$, we have
    \begin{align*}
        R(f^*_{\Xi, \HK}) \, \leq \, &
        \inf_{f \in \cF} R(f)
        + 2d^{1/2}(p^{1/2} + 3)p^{1/2}L_\varphi L_z \bar{X} K^{-1/p}\\
        & +\frac{16\mu L_\varphi\sqrt{2dK(p+1)(p\bar{X}^2 + 1)}}{n}
        + 4M_{p,d}(\mu)\sqrt{\frac{2dK(p+1)\log n}{n}}
        + 2\sqrt{2}M_{p,d}(\mu) \sqrt{\frac{\log(2 / \delta)}{n}}.
    \end{align*}
\end{theorem}

The above theorem implies that with $\mu=\cO(K^{1/p})$, the complexity of the excess risk bound with respect to the sample size $n$ and parameter $K$ is
\begin{equation}\label{eq:model:excess risk complexity of N k}
    R(f^*_{\Xi, \HK}) - \inf_{f \in \cF} R(f) 
    = \underbrace{\widetilde{\cO}\left(K^{\frac{p+2}{2p}
    }n^{-\frac{1}{2}}\right)}_{\substack{\text{estimation}\\\text{error bound}}} 
    + \underbrace{\cO\left(K^{-\frac{1}{p}}\right)}_{\substack{\text{approximation}\\\text{error bound}}},
\end{equation} 
where $\widetilde{\cO}(\cdot)$ suppresses logarithmic factors.
It can be seen that the sample size $n$ only affects the estimation error, whereas $K$, the number of max-affine pieces for PA functions in $\HK$, plays the central role of balancing the two error terms. 
As $K$ grows, the PA function class becomes larger, so the approximation error decreases due to the stronger approximation ability of $\HK$, and the estimation error increases due to the more pronounced overfitting issue.
With an appropriate choice of $K$ that balances these two errors, we can derive the optimal excess risk bound. Furthermore, with the existence of a convex or DC optimal DR, this bound can be improved based on the approximation bounds $\cO(K^{-2/p})$ of $\HK$ from \citet{balazs2015near} and \citet{siahkamari2020piecewise}. The extensions are summarized in the following corollary.
\begin{corollary}\label{coro:optimal excessriskbound}
    Under the settings in Theorem~\ref{thm:excess-wo-constr}, with $K = \cO(n^{p/(p+4)})$, the optimal excess risk bound is $\widetilde{\cO}(n^{-{1}/(4+p)})$. 
    If there exists a convex or DC optimal decision rule to problem~\eqref{eq:intro:contextalSP:Full-E}, the optimal excess risk bound is $\widetilde{\cO}(n^{-1/(3 + p/2)})$ with $K = \cO(n^{p/(p+6)})$.
\end{corollary}

We make a remark that Proposition~\ref{prop:model:universal approximation of PADR} provides a bound for the worst-case approximation error and thus the excess risk bound~\eqref{eq:model:excess risk complexity of N k} could be conservative in practice.  In the practical implementation of PADR in Section~\ref{section:experiments}, we observe that a much more moderate number of $K$ is able to produce satisfactory results including the two-stage product placement problem and Newsvendor problems with real data. This is because the (near) optimal decision rule may possess some underlying structural property besides Lipschitz continuity so that it could be approximated well by the PADR with a moderate number of pieces, for which a smaller estimation error and excess risk could be achieved. 

The excess risk bound~\eqref{eq:model:excess risk complexity of N k} illustrates the role of $K$ in balancing the two errors, which is further evidenced by numerical experiments summarized in the table below. Table~\ref{tab:PADRapproximationII} reports the validation and training costs provided by the PADR-based ERM method in solving a Newsvendor problem with the sine plus max-affine demand model~\eqref{eq:exp:sin+ma model} (details are given in Section~\ref{subsec:exp:Newsvendor}), under multiple choices of $K_1$ and $K_2$ which represent the number of pieces for the two max-affine components in the PADR respectively. The validation cost, which is the empirical cost of the obtained PADR solution over the validation data set of size 1000, can be regarded as the sample average approximation of $R(f^*_{\Xi, \HK})$. Noticing that  $\inf_{f \in \cF} R(f) $ is irrelevant to the PA hypothesis class,  the performance of excess risk could be revealed through the pattern of validation costs in the table below. 
{The training cost decreases consistently as $K_1$ and $K_2$ increase, indicating the stronger fitting capacity of PADR on the training data. 
For small values of $K_1$ and $K_2$, the validation cost also decreases significantly, which can be attributed to the reduction in approximation error. 
However, this decreasing trend does not hold as $K_1$ and $K_2$ become large;
instead, the validation cost begins to plateau or even slightly increase.
For instance, when $K_2$ is fixed at 40, the validation cost first decreases to 4.13 (at $K_1=5$) and then increases to 4.24 (at $K_1=20$).
This suggests that as the PADR model becomes more complex, the increasing estimation error gradually dominates the decreasing approximation error. 
Hence, the best validation costs shall be achieved with moderate values of $K_1$ and $K_2$, which are usually considerably smaller than the ones implied by the theoretical results.}

\begin{table}[hbtp]
  \TABLE
  {Validation and training costs for  basic Newsvendor with demand model~\eqref{eq:exp:sin+ma model} \label{tab:PADRapproximationII}}
  {%
      \resizebox{0.99\columnwidth}{!}{\begin{tabular}{ccccccccc}
      \toprule
      K1 $\backslash$ K2 & 1     & 2     & 3     & 4     & 5     & 10    & 20    & 40 \\
      \midrule
      1     & 24.06 (24.09) & 23.62 (23.57) & 23.61 (23.61) & 23.64 (23.55) & 23.56 (23.62) & 23.55 (23.64) & 23.54 (23.62) & 23.56 (23.64) \\
      2     & 8.68 (8.73) & 7.01 (7.09) & 7.07 (7.00) & 7.02 (7.08) & 7.08 (6.99) & 7.05 (7.02) & 7.05 (7.00) & 7.05 (7.01) \\
      3     & 7.95 (7.79) & 5.60 (5.40) & 5.47 (5.39) & 5.40 (5.35) & 5.42 (5.34) & 5.36 (5.33) & 5.41 (5.37) &  5.39 (5.35)\\
      4     & 7.29 (7.47) & 4.26 (4.27) & 4.25 (4.14) & 4.28 (4.29) & 4.14 (4.17) & 4.13 (4.25) & 4.15 (4.15) & 4.14 (4.21) \\
      5     & 7.41 (7.23) & 4.21 (4.17) & 4.34 (4.14) & 4.32 (4.15) & 4.17 (4.13) & 4.34 (4.16) & 4.14 (4.15) & 4.11 (4.13) \\
      10    & 7.41 (7.25) & 4.15 (4.18) & 4.15 (\bluelog{4.11}) & 4.10 (4.13) & 4.14 (4.14) & 4.06 (4.24) & 4.04 (4.25) & 4.07 (4.14) \\
      20    & 7.41 (7.24) & 4.16 (4.21) & 4.15 (4.14) & 4.10 (4.14) & 4.32 (\bluelog{4.11}) & 4.04 (4.14) & 4.00 (4.15) & 4.03 (4.24) \\
      40    & 7.40 (7.23) & 4.15 (4.14) & 4.12 (4.15) & 4.08 (4.15) & 4.11 (4.13) & 4.01 (4.24) & \bluelog{3.99} (4.26) & 4.02 (4.15) \\
      \bottomrule
      \end{tabular}}
  }
    {The PADR is obtained by Algorithm~\ref{ALGO:unconstrained case} for solving the PADR-based ERM model.  Numbers listed in brackets represent training costs, and the lowest costs are highlighted in bold and blue text.}
\end{table}%

\subsection{Consistency comparison with existing DR-based methods}\label{subsec:model:comparison of theoretical results with benchmarks}
We compare the consistency result of the PADR-based ERM model with two prominent types of DR-based ERM methods in Table~\ref{tab:comparison of risk bounds}.
NV-ERM1 and NV-ERM2 represent the two LDR methods from~\citet{ban2019big}, the latter of which additionally applies regularization in the objective.
 \RKHSnp and \RKHSbk represent RKHS-DR methods from \citet{notz2021prescriptive} and \citet{bertsimas2022data} respectively, with the former imposing regularization via constraints and the latter applying regularization in the objective.
\begin{table}[htbp]
    \TABLE{Comparison of error bounds for multiple DR-based ERM models.\label{tab:comparison of risk bounds}}
    {\begin{tabular}{ccccc}
      \toprule
      {Problem type} & {Method} & {Approximation error} & {Estimation error} & {Excess risk} \\
      \midrule
      \multirow{4}{*}{Convex} & NV-ERM1                  & 0 (Linear asm.) & $\widetilde{\cO}(n^{-\frac{1}{2+p/2}})$ & $\widetilde{\cO}(n^{-\frac{1}{2+p/2}})$ \\
                              & NV-ERM2                  & 0 (Linear asm.) & $\geq \cO(n^{-\frac{1}{2+p/4}})$                   & /                                       \\
                              & {\RKHSnp}  & /               & $\cO(n^{-\frac{1}{2}})$     &    /  \\
                              & {\RKHSbk}  & 0               & $\geq \cO(n^{-\frac{1}{2}})$                       &   /   \\
      \midrule
      Nonconvex & PADR & $\cO(K^{-\frac{1}{p}})$  & $\widetilde{\cO}(K^{\frac{p+2}{2p}}n^{-\frac{1}{2}})$ & $\widetilde{\cO}(n^{-\frac{1}{4+p}})$\\
      \bottomrule
    \end{tabular}}
    {The complexity with ``$\geq$'' means that the error is only bounded partly. The notation ``/'' means that the error bound is not provided.}
\end{table}%

Out of the four benchmarks considered in the literature, only \citet{ban2019big} provide an excess risk bound for NV-ERM1,
 leveraging the special structure of the Newsvendor problem with a linear true demand model.
The estimation error bounds for NV-ERM2 and \RKHSbk are not provided due to the unknown regularization bias.
Since an RKHS with a universal kernel is dense in the space of continuous functions, there is no approximation error for \RKHSbk, while such an error of \RKHSnp with a norm constraint is not analyzed in \citet{notz2021prescriptive}. It should be noted that the theoretical consistency results for all benchmarks in the aforementioned references, except \RKHSnp, are restricted to convex SP problems.

In comparison, 
our PADR-based ERM model is nonasymptotically consistent in more general cases, including convex SP problems and nonconvex ones.
Moreover, by Corollary~\ref{coro:optimal excessriskbound}, if there exists a convex or DC optimal DR, our excess risk bound can be improved to $\widetilde{\cO}(n^{-1/(3 + p/2)})$, which is close to the bound of NV-ERM1 (a special case of our model) when $p$ is large.
This gap will vanish, that is, the bound can be further refined to $\widetilde{\cO}(n^{-1/(2 + p/2)})$, for problem~\eqref{eq:intro:contextalSP} with a box-constrained set $[-M, M]^d$ %
by employing a truncated PADR class $\barHK \coloneqq \{\bar{f}\mid \bar{f} = \max\{\min\{f, M\}, -M\},\,f\in \HK\}$ in~\eqref{eq:intro:ERM-H}, under which $M_{p,d}(\mu)$ in Theorem~\ref{thm:excess-wo-constr} could be independent of $\mu$.

Given the statistical consistency merits of the PADR method demonstrated by the analysis above, it is still prohibitively challenging to obtain the globally optimal solution to the nonconvex PADR-based ERM problem~\eqref{eq:intro:ERM-PADR}. From the continuous optimization perspective, the PADR-based ERM problem is a composite nonconvex and nondifferentiable optimization problem embedded with a particular piecewise affine structure. This leads to the development of a sampling-based surrogate-based algorithm in the next section to efficiently solve problem~\eqref{eq:intro:ERM-PADR} of large sizes. %

\section{Efficient Sampling-based Algorithm for Learning PADR}\label{section:algorithm}

In this section, we focus on learning PADR by solving the following parameter estimation problem
\begin{equation}
    \label{eq:alg:ERM-PADR}
    \min_{\norm{\theta_{\iota}}_\infty \leq \mu, \forall \iota \in [d]} \quad  
    \frac{1}{n}\sum_{s=1}^n \varphi\left( \Big( \max_{k\in [K]} \left\{ (\alpha^k_{\iota})^\top x^s + a^k_{\iota} \right\} 
          - \max_{k\in [K]} \left\{ (\beta^k_{\iota})^\top x^s + b^k_{\iota} \right\} \Big)_{\iota=1}^d, y^s \right),
\end{equation}
where the objective is the finite sum of composite-PA functions. Instead of attempting to achieve global optimality, it is more practical to apply majorization-minimization (MM)-type methods to the contextual optimization problems with the practical scale in seeking ``stationary'' solutions. 
Among various kinds of stationary concepts, the d-stationary solution is the sharpest type of first-order stationary solutions for general nonconvex optimization.
For the PADR learning problem~\eqref{eq:alg:ERM-PADR} with a general nonconvex and nondifferentiable composite structure, our algorithmic development is mainly motivated by the following question: 
  
\emph{Given the potentially large sample size, how to efficiently compute a sharper ``stationary'' solution to the PADR-based ERM problem~\eqref{eq:alg:ERM-PADR} so that a global minimizer could be more easily obtained through replications with random initialization?}

Building on this motivation, we develop an MM-type algorithm specializing to the structured PADR-based ERM problem with the theoretical convergence analysis. It is further observed with numerous experiments in Section~\ref{section:experiments} that satisfactory (nearly optimal) solutions can be obtained with a few number of repetitions, which is up to 4 times faster than the basic MM method.

To make our discussion self-contained, we first present the standard MM algorithms in Section~\ref{sec:standardMM}, which serves as a reference for the motivation behind our proposed sampling-based algorithm in Section~\ref{sub:esmm}. We make some remarks on notations to facilitate algorithmic development and convergence analysis. To emphasize the role of PADR as a function of $\theta$, with $\xi=(x, y) \in \Xi $, we denote $g(\theta, \xi) \coloneqq \max_{i \in [K]} g_{i}(\theta, x)$ and
$h(\theta, \xi) \coloneqq \max_{i \in [K]} h_{i}(\theta, x)$ where $\{g_{i}(\theta, x)\}$ and $\{h_{i}(\theta, x)\}$ are linear functions with respect to $\theta$, and we denote 
$f(\theta, \xi) \coloneqq g(\theta, \xi) - h(\theta, \xi)$ as the PA function  
 with a slight abuse of notation.  

With $\widetilde{\xi}_n$ representing a random variable uniformly distributed on $\Xi$,
the PADR-based ERM problem~\eqref{eq:intro:ERM-PADR} can be generalized as the following stochastic program with finite scenarios:
\begin{equation}
    \label{eq:algo:formulation:finite sum primal problem}
    \min_{\theta \in \Theta}\quad  
    F(\theta) \coloneqq \mathbb{E} \,\left[F(\theta,\widetilde{\xi}_n)\right],
\end{equation}
where $F(\theta, \xi) \coloneqq \varphi(f(\theta,  {\xi}),  {\xi}) $ and $\Theta \subseteq \RR^q$ is a compact convex set. 
In particular, it reduces to problem~\eqref{eq:intro:ERM-PADR} when $\Theta = [-\mu, \mu]^q$ with $q = 2dK(p+1)$, the outer function $\varphi(\cdot, \xi)$ with $\xi = (x, y)$ depends solely on $y$, and the inner PA function $f(\cdot, \xi)$ depends solely on $x$. We say a function $F: \Theta \subseteq \mathbb R^q \to \mathbb R$ is directionally differentiable if the directional derivative $F'(\theta; d) \coloneqq \lim_{\tau \to 0} \frac{F(\theta + \tau d) - F(\theta)}{\tau} $ exists for any $\theta \in \Theta$ and $d \in \mathbb R^{q}$. Supposing that $\Theta$ is a closed convex set and $F:\Theta\rightarrow\RR$ is directionally differentiable,  we say that $\theta^d \in \Theta$ is a $\text{d(irectional)-stationary}$ point of problem~\eqref{eq:algo:formulation:finite sum primal problem} if $F'(\theta^d; \theta - \theta^d) \geq 0$ for any $\theta \in \Theta$. %

 \subsection{Preliminary analysis: standard MM algorithms}
 \label{sec:standardMM}

The core of MM-type methods for solving nonconvex minimization problem 
is to iteratively solve subproblems with a surrogate function $\widehat{F}(\theta; \theta^\nu)$ majorizing $F(\theta)$ at an iterate point $\theta^\nu \in \Theta$. 
The iterate point is updated with $\eta \geq 0$ as follows,
\begin{equation}\label{eq:iterative process}
    \theta^{\nu + 1}  \leftarrow   
\underset{\theta \in \Theta}{\mathop{\arg\min}}
\left\{\widehat{F}(\theta; \theta^\nu) + \frac{\eta}{2}\norm{\theta - \theta^\nu}^2\right\}.
\end{equation}
    
A surrogate function $\widehat{F}(\theta; \bar{\theta})$ at a reference point $\bar{\theta} \in \Theta$ typically satisfies the following conditions:
\begin{enumerate}[label = (P\arabic*), series=P]
    \item\label{enum:algo:P1:touching} Touching condition:  $\widehat{F}(\bar{\theta}; \bar{\theta}) = F(\bar{\theta})$.
    \item\label{enum:algo:P2:majorization} Majorizing condition: $\widehat F({\theta};{\wtheta}) \geq F(\theta)$ for any $\theta\in \Theta$.
    \item\label{enum:algo:P3:convexity} Convexity: $\widehat F({\cdot};{\wtheta})$ is convex on $\Theta$.
\end{enumerate}
A function satisfying~\ref{enum:algo:P1:touching} and~\ref{enum:algo:P2:majorization} is referred to as an upper surrogate function according to \citet[Definition 7.1.1]{cui2021modern}.  Under~\ref{enum:algo:P1:touching} and~\ref{enum:algo:P2:majorization}, it yields the sequence of decreasing objective values as follows,
\[
F(\theta^{\nu+1} ) \leq \widehat{F}(\theta^{\nu+1}; \theta^\nu)  \leq \widehat{F}(\theta^\nu; \theta^\nu) -\frac{\eta}{2}\norm{\theta^{\nu+1} - \theta^\nu}^2  = F(\theta^\nu) -\frac{\eta}{2}\norm{\theta^{\nu+1} - \theta^\nu}^2 \leq F(\theta^\nu).
\]
Condition~\ref{enum:algo:P3:convexity} ensures that surrogation-based subproblems in the algorithm can be solved by state-of-the-art convex-programming solvers. 
By the basic analysis (see Chapter 7 in \cite{cui2021modern}), the limit point of 
$\{\theta^{\nu}\}$ produced by the process~\eqref{eq:iterative process} could be easily shown to be a fixed point of the surrogation map satisfying $\theta^{*}  = \underset{\theta \in \Theta}{\mathop{\arg\min}}
\left\{\widehat{F}(\theta; \theta^*) + \frac{\eta}{2}\norm{\theta - \theta^*}^2\right\}$.
Such a fixed point property could be related to various types of stationarity of the original problem depending on the surrogate functions as illustrated in \citet{cui2018composite} and \citet{liu2022solving}. %

\begin{figure}[htbp]
    \FIGURE
    {
    {
        \subfloat[DC function]{%
          \includegraphics[width=0.326\textwidth]{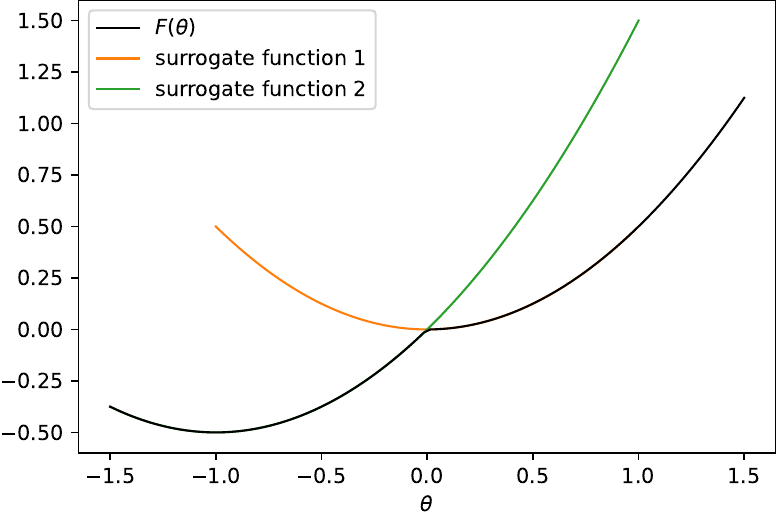}%
          \label{subfig:algorithm:dcfunction}
        }
        \subfloat[Surrogation with $\varepsilon=0$]{%
          \includegraphics[width=0.326\textwidth]{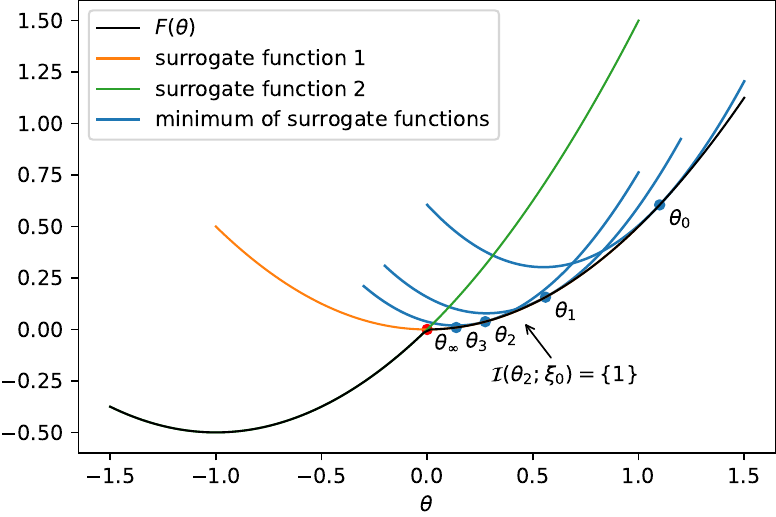}%
          \label{subfig:algorithm:mm}
        }
        \subfloat[Surrogation with $\varepsilon = 0.3$]{%
          \includegraphics[width=0.326\textwidth]{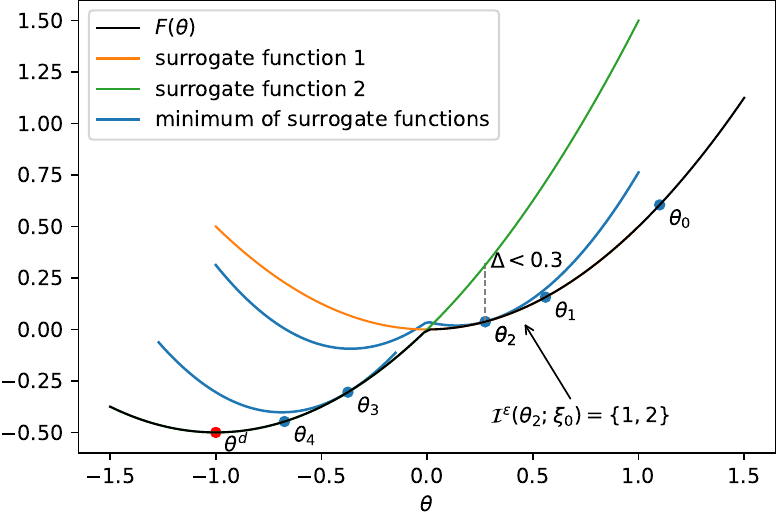}%
          \label{subfig:algorithm:emm}
        }
      }
    }
    {Iterative process of the standard MM algorithm and enhanced MM algorithm \label{fig:algorithm:iterating}\label{fig:MM-procedure}}
    {}
  \end{figure}

Regarding the PADR-based ERM problem~\eqref{eq:algo:formulation:finite sum primal problem}, due to the nonconvex and nondifferentiable piecewise affine functions, a major issue of the standard MM algorithms including \citet{lipp2016variations} and \citet{le2018dc} is that the fixed point corresponding to the surrogation map is a (compound) critical solution (see definitions in \cite{pang2017computing,liu2022solving}) which is much weaker than directional stationarity. Specifically, based on the composite structure in~\eqref{eq:algo:formulation:finite sum primal problem}, approaches to construct surrogate functions of $F(\theta, \xi)$ in \cite{cui2018composite} %
typically involve linearization of the two max-affine components of the inner piecewise affine function as follows,
\begin{align}
\label{eq:upper_PA}
    \widehat{f}(\theta, \xi; I_2)& \coloneqq g(\theta, \xi) - h_{I_2}(\theta, \xi) \mbox{ with } I_{2} \in \cI_h(\wtheta, \xi) \coloneqq \mbox{argmax}_{i \in [K]} h_i(\bar{\theta}, \xi), \\
    \label{eq:lower_PA}
    \widecheck{f}(\theta, \xi;I_1) & \coloneqq  g_{I_1}(\theta, \xi) - h(\theta, \xi) \mbox{ with } I_{1} \in \cI_g(\wtheta, \xi) \coloneqq  \mbox{argmax}_{i \in [K]} g_i(\bar{\theta}, \xi). 
\end{align}
{
Utilizing the above upper and lower surrogate functions of $f(\cdot, \xi)$ at $\bar{\theta}\in \Theta$, the following example illustrates the construction of the surrogate function when outer function $\varphi(\cdot, \xi)$ is linear. 
\begin{example}\label{eg:surrogation under linear varphi}
    For $\varphi(z,\xi) = c_\xi^\top z$, where $c_\xi \in \RR^d$ is a parameter dependent on $\xi$, the composite objective function is $F(\theta, \xi) = c_\xi^\top f(\theta, \xi)$, and its surrogate function is constructed as
    \[
    \widehat{F}(\theta, \xi; \bar{\theta}, I) = \max\{c_\xi, 0\}^\top\widehat{f}(\theta, \xi; I_2) + \min\{c_\xi, 0\}^\top\widecheck{f}(\theta, \xi; I_1),\]
where $I \coloneqq \big(I_{1},I_{2}\big) \in \cI_g(\bar{\theta},\xi) \times \cI_h(\bar{\theta},\xi)$. For $d > 1$, the maximum is taken componentwisely.
\end{example}
Surrogation in Example~\ref{eg:surrogation under linear varphi} obviously satisfies Conditions~\ref{enum:algo:P1:touching} to~\ref{enum:algo:P3:convexity}, and the standard MM algorithm for such a case is exactly a DC algorithm. 
}
Clearly $f'(\cdot, \xi)(\theta; d)= g'(\cdot, \xi)(\theta;d) - \max\{ \, \nabla h_i(\wtheta, \xi)^\top d: \, i \in \mbox{argmax}_{i \in [K]} h_i(\bar{\theta}, \xi)\}$, whereas the surrogate function $\widehat{F}(\theta, \xi; \bar{\theta}, I)$ employs an arbitrary active index $I \in \cI_g(\wtheta, \xi) \times \cI_h(\wtheta, \xi) $. Hence, the fixed point of the surrogation map is a critical solution. Solutions of criticality are not desirable because in practice solutions may be at the ``cliff'' of the cost landscape, which is illustrated in Example \ref{eg:plot of iterative process} and Figure \ref{subfig:algorithm:mm}.

To obtain a directional stationary solution, 
\cite{pang2017computing} considered a class of nondifferentiable DC programs which have a max-smooth structure in the second convex part, and introduced $\varepsilon$-active sets to construct the enhanced surrogate functions for DC functions. Specifically, for the piecewise affine function $f(\cdot, \xi)$ given $\xi \in \Xi$, we denote the $\varepsilon$-active set with $\varepsilon > 0$ by 
\begin{equation}
\label{eq:epsilon-active-set}
    \cI_g^\varepsilon(\bar{\theta}, \xi) \coloneqq \Big\{
        i \in [K]
        \mid
        \,g_{i}(\bar{\theta}, \xi) \geq g(\bar{\theta}, \xi) - \varepsilon
    \Big\},
\end{equation}
and the set $\cI_h^\varepsilon(\bar{\theta}, \xi)$ is similarly defined for $h(\cdot, \xi)$. Let $\cI^\varepsilon(\bar{\theta},\xi) \coloneqq  \cI_g^\varepsilon(\bar{\theta}, \xi)  \times  \cI_h^\varepsilon(\bar{\theta}, \xi) $. \citet{pang2017computing} proposes the enhanced MM algorithm in which the subproblem takes the minimum over the family of surrogate functions with indices in $\cI^\varepsilon(\bar{\theta},\xi)$ and provides the theoretical convergence analysis to a d-stationarity point. The following example provides an intuitive illustration that the iterate point can escape from a critical point to better local minima and make valid descents via the enhanced technique.

\begin{example}\label{eg:plot of iterative process}
Consider problem~\eqref{eq:algo:formulation:finite sum primal problem} with $d=2$, $q=1$, $\Xi = \{\xi_0\}$, $
\varphi(z, \xi_0) = z_1^2 + z_2, f_1(\theta, \xi_0) = \theta$, and $ f_2(\theta, \xi_0) = -\max\{-\theta, 0\}$. {In this case, the composite objective function is $F(\theta) = \theta^2 + \max\{-\theta, 0\}$.} Figure~\ref{subfig:algorithm:dcfunction} illustrates the two candidates of surrogate functions depending on which piece of the inner PA function $f_2(\theta, \xi_0) $ is active at the reference point of interest.  
Figure~\ref{subfig:algorithm:mm} and Figure~\ref{subfig:algorithm:emm} illustrate the iterative procedures of the standard MM algorithm and the enhanced MM algorithm respectively, with the initial point $\theta_0=1.1$. The standard MM algorithm converges to a critical point $\theta_\infty=0$. For the enhanced MM algorithm with $\varepsilon=0.3$, at the second iterate point $\theta_2$, we have $\abs{\mathcal I^\varepsilon(\theta_2, \xi_0)} = 2$ at $\theta_2$, which thus leads to the next iterate point passing over 0 and continues to move towards a d-stationary point $\theta^d$ that is also globally optimal to the problem.
\end{example}

Despite possessing stronger convergence properties, the enhanced MM algorithm applied for large-scale ERM problems could be exceptionally inefficient as it requires  
iteratively solving a potentially very large number of finite-sum convex subproblems over all the data.  Sampling-based MM-type algorithms have been recently studied for (partially) smooth stochastic DC programs in \citet{nitanda2017stochastic, xu2019stochastic}, and \citet{le2022stochastic} and nondifferentiable stochastic composite DC problems in \citet{liu2022solving} to compute compound critical solutions. Based on a natural idea of combining a sequential sampling scheme with the enhanced technique,  we propose an enhanced sampling-based MM algorithm with rigorous convergence analysis in the next subsection to solve the stochastic composite DC programs~\eqref{eq:algo:formulation:finite sum primal problem} with both computational efficiency and theoretical convergence guarantee.

\subsection{Enhanced sampling-based MM algorithm for learning PADR}\label{sub:esmm}

In this subsection, we first introduce the family of enhanced surrogate functions of the composite objective $F$ and provide several function classes of the outer function $\varphi$ as particular instances.
This provides the basis for the characterization of a so-called \emph{composite $(\varepsilon, \rho)$-strong d-stationarity}. %
We then propose the enhanced sampling-based MM algorithm (ESMM) which combines the enhanced $\varepsilon$-active technique with sequential sampling and analyze the asymptotic and non-asymptotic convergence properties in terms of (composite strong) d-stationarity under suitable incremental sampling schemes.  Such a combination enables the algorithm to solve subproblems of small sizes in the early iterations that effectively reduces computational costs. The computational improvement of such a technique is further supported by numerical experiments in Section~\ref{section:experiments}.

\subsubsection{Family of enhanced surrogate functions.}\label{subsec:surrogate}

Based on surrogation of $f(\theta, \xi)$ in~\eqref{eq:upper_PA} and~\eqref{eq:lower_PA} and corresponding $\varepsilon$-active index set $\mathcal I^{\varepsilon}(\bar \theta,\xi) \coloneqq  \cI_g^\varepsilon(\bar{\theta}, \xi)  \times  \cI_h^\varepsilon(\bar{\theta}, \xi) $ given $\bar \theta \in \Theta$, we introduce a family of enhanced surrogate functions  $\{\widehat{F}(\cdot, \xi; \bar{\theta}, I)\}_{I \in \mathcal I^\varepsilon(\bar{\theta},\, \xi)}$ for any $\xi \in \Xi$ that satisfy the following conditions:
\begin{enumerate}[label=(B\arabic*), series=B]
    \item \label{assum:B1:Touching}$\widehat{F}(\cdot, \xi; \bar{\theta}, I)$ satisfies~\ref{enum:algo:P1:touching}  for every $I \in \cI^0(\bar{\theta}, \, \xi)$;
    \item \label{assum:B2:Majorizing and Convex} $\widehat{F}(\cdot, \xi; \bar{\theta}, I)$ satisfies~\ref{enum:algo:P2:majorization} and~\ref{enum:algo:P3:convexity} for every $I \in \cI^\varepsilon(\bar{\theta}, \, \xi)$ with $\varepsilon \geq 0$;
    \item \label{assum:B3:Continuous given fixed I} $\widehat{F}(\cdot, \xi;\, \cdot, I)$ is continuous on $\Theta\times\Theta$ given fixed indices $I = (I_1, I_2) \in [K]^2 $;
    \item \label{assum:B4:minSurr dd-consistency} the minimum surrogate function $\widehat{F}^\varepsilon_{\min}(\theta, \xi; \bar{\theta})\coloneqq\min_{I \in \cI^\varepsilon(\bar{\theta}, \, \xi)} \widehat{F}(\theta, \xi; \bar{\theta}, I)$ satisfies the directional derivative consistency for any $\varepsilon \geq 0$; namely, 
    $\widehat F_{\min}^{\varepsilon}(\cdot\,; \bar{\theta})$ is directionally differentiable on $\Theta$ and satisfies
    $\widehat F_{\min}^{\varepsilon}(\cdot\,; \bar{\theta})'(\bar{\theta}; d) = F'(\bar{\theta}; d)$ for any $d \in \RR^q$.
\end{enumerate}
The directional derivative consistency condition is introduced in {\citet{cui2021modern}  to relate the directional derivatives of surrogation and original functions. 
We consider the following three cases of the outer function $\varphi$ covering a broad class of cost functions. 
\begin{enumerate}[label = {\it Case \arabic*.}, ref = {\it Case~\arabic*}]
    \item\label{enum:algo:phi1} $\varphi(\cdot, \xi)$ is piecewise affine for any $\xi \in \Xi$.
    \item\label{enum:algo:phi2} $\varphi(\cdot, \xi)$ can be decomposed into monotonic convex functions for any $\xi \in \Xi$.
    \item\label{enum:algo:phi3} $\varphi(\cdot, \xi)$ is smooth with Lipschitz gradient for any $\xi \in \Xi$.
\end{enumerate}
For the above cases, we can construct the surrogate function of $F(\cdot,\xi)$ satisfying Assumptions~\ref{assum:B1:Touching} to~\ref{assum:B4:minSurr dd-consistency}, by utilizing surrogate functions~\eqref{eq:upper_PA} and~\eqref{eq:lower_PA} of the inner PA function $f(\cdot,\xi)$ with the selected index $I = (I_1, I_2)$. The specific construction and verification are provided in Online Appendix~\ref{EC:sec:surrogation for sec3}. In particular, \ref{enum:algo:phi2} is a multivariate extension of the one considered in \citet{cui2018composite}.  
\ref{enum:algo:phi3} is suitable for nonconvex smooth cost functions that are applicable for many practical problems.  
\begin{proposition}\label{prop:verify Cases1-3}
    Let $\Theta$ be a compact set. 
    For~\ref{enum:algo:phi1} to~\ref{enum:algo:phi3}, given $\varepsilon \geq 0$ and $\bar{\theta} \in \Theta$, there exists a family of enhanced surrogate functions based on the index set $\mathcal I^{\varepsilon}(\bar \theta, \xi)$ of the inner PA function that satisfy Assumptions~\ref{assum:B1:Touching}-\ref{assum:B4:minSurr dd-consistency} for any $\xi \in \Xi$.
\end{proposition}

\subsubsection{Concepts on stationarity.}\label{subsec:stationarity}

For a general nonconvex problem on a convex set with the directional differentiable objective, d-stationarity is identified as the sharpest kind among the first-order stationarity concepts. Specialized to the composite DC program~\eqref{eq:algo:formulation:finite sum primal problem}, according to \citet[Proposition 2]{cui2018composite}, d-stationarity is equivalent to local optimality if $\varphi$ is in~\ref{enum:algo:phi1} or~\ref{enum:algo:phi2}. 
\citet{lu2019enhanced} proposes the strengthened d-stationarity for structured
DC programs which takes into account the $\epsilon$-active index set. 
We extend to the concept of composite $(\varepsilon, \rho)$-strong d-stationarity specialized for the PADR-based ERM problem~\eqref{eq:alg:ERM-PADR}.
The composite strongly d-stationary point can be interpreted as the fixed point corresponding to the proximal mapping of the expected minimum surrogate function $\mathbb E [\, \widehat{F}^\varepsilon_{\min}(\theta, \widetilde{\xi}^n; \bar{\theta})]$ based on the family of enhanced surrogate functions  $\{\widehat{F}(\cdot, \xi; \bar{\theta}, I)\}_{I \in \mathcal I^\varepsilon(\bar{\theta},\, \xi)}$ satisfying Assumptions~\ref{assum:B1:Touching}-\ref{assum:B4:minSurr dd-consistency}.

\begin{definition}\label{def:strong d-stat}
    Given $\varepsilon \geq 0$ and $\rho > 0$,
    we say that $\theta^d \in \Theta$ is a composite $(\varepsilon, \rho)$-strongly d-stationary point of problem~\eqref{eq:algo:formulation:finite sum primal problem} if 
    $\displaystyle{\{\theta^d\} = \mathop{\arg\min}_{\theta \in \Theta}\left\{ \mathbb E [\, \widehat{F}^\varepsilon_{\min}(\theta, \widetilde{\xi}^n; {\theta}^d) \,] + \frac{\rho}{2}\norm{\theta - \theta^d}^2
    \right\}}$. Correspondingly, let $\cS^d_{\varepsilon, \rho}$ denote the set of composite $(\varepsilon, \rho)$-strongly d-stationary points. %
\end{definition}

By Condition~\ref{assum:B4:minSurr dd-consistency}, we must have $\cS^d_{\varepsilon, \rho} \subseteq \cS_d$ for any $\varepsilon \geq 0$ where $\cS_d$ denotes the set of d-stationary points. 
Moreover, with $\varepsilon >0$, we have the inclusion property $\mathcal I^0(x, \xi) \subseteq \mathcal I^0(\hat x, \xi) \subseteq \mathcal I^\varepsilon(x, \xi) $ for $x$ sufficiently close to $\hat x$, which is crucial for the convergence analysis of the limit point in terms of (composite strong) d-stationarity.
With $\cS^*$ being the set of global minima of~\eqref{eq:algo:formulation:finite sum primal problem}, 
the following proposition gives the relation among $\cS^*, \cS^d_{\varepsilon, \rho}$ and $ \cS^d$. The proof is provided in Online Appendix~\ref{proof:EC:proof of S*Sd relationships}. 
\begin{proposition}\label{prop:convergence:relationship-S*-Sds}
    $\cS^* \subseteq \cS^d_{\varepsilon, \rho} \subseteq \cS^d$ for any $\varepsilon \geq 0$ and $\rho > 0$. In particular, $\cS^d_{0, \rho} = \cS^d$ for any $\rho > 0$.
\end{proposition}
Due to the requirement of a unique proximal mapping point in the defining condition, 
our concept is slightly stronger than $(\alpha, \eta)$-d-stationarity defined in \citet{lu2019enhanced} when the problem is reduced to composite DC program. 
Moreover, it is weaker than composite $\varepsilon$-strong d-stationarity introduced in \citet{qi2022asymptotic} (see also \citet[Definition 6.1.3]{cui2021modern}) due to the presence of an additional proximal term in the defining condition which turns out to be essential for the algorithm development and the follow-up convergence analysis. 
Furthermore, in Online Appendix~\ref{EC:subsec:proofs for Section4-1}, we show that $\cS^d_{\varepsilon, \rho}$ is equivalent to $\cS^*$ with sufficiently large $\varepsilon$ under a regularity assumption on $F$. It implies that the definition of $\cS^d_{\varepsilon, \rho}$ serves as a natural connection between d-stationarity and global optimality.

\subsubsection{ESMM algorithm and convergence analysis.}\label{subsec:the surrogation and the algorithm}

We propose the ESMM algorithm with a natural idea of combining the randomized $\varepsilon$-active surrogation with sequential sampling presented in Algorithm~\ref{ALGO:unconstrained case}. At iteration $\nu$, a subset of samples $\Xi_\nu = \{\wtxi_\nu^s\}_{s=1}^{N_{\nu}}$ of size $N_{\nu}$ is generated in step 2. 
Then in step 3 a sampling-based surrogate function $\widehat{F}_{\nu}(\theta; \theta^\nu, \wtbI_{\nu}^\varepsilon)$ is constructed which can be regarded as the upper surrogation function to the sample average approximation function $F_{\nu}(\cdot) \coloneqq \frac{1}{N_{\nu}} \sum_{s=1}^{N_{\nu}} F(\cdot, \wtxi_\nu^s)$ with the combination of indices $\wtbI_{\nu}^\varepsilon =\{\wtI^\varepsilon_{\nu,s}\}_{s=1}^{N_{\nu}}$ randomly selected from the $\varepsilon$-active index sets $\{\cI^\varepsilon(\theta^\nu, \wtxi_{\nu}^s)\}_{s=1}^{N_{\nu}}$, respectively.  
The randomized selection combined with the sampling scheme enables the algorithm to solve subproblems of small sizes in the early iterations. %
By the convexity of the surrogate function, in step~4 the proximal mapping point $\theta^{\nu+\frac{1}{2}}$ can be obtained via the state-of-the-art convex programming solvers.  However, since the indices $\wtbI_{\nu}^{\varepsilon}$ are randomly selected from the $\varepsilon$-active index sets, the touching condition~\ref{enum:algo:P1:touching} may not be satisfied when $\varepsilon >0$; namely, the surrogate function value $\widehat{F}_{\nu}(\theta^{\nu}; \theta^\nu, \wtbI_{\nu}^\varepsilon)$ might be strictly larger than $F_{\nu}(\theta^\nu)$. Thus the candidate solution $\theta^{\nu + \frac{1}{2}}$ may not lead to a valid update with the desired descent property.
Hence, in step~5 the candidate solution $\theta^{\nu + \frac{1}{2}}$ is accepted only when sufficient descent is satisfied. In the special case of $\varepsilon = 0$ for which the touching condition~\ref{enum:algo:P1:touching} is satisfied,  the proximal mapping point $\theta^{\nu + \frac{1}{2}}$ satisfies the sufficient descent inequality in step~5, and thus we could directly update $\theta^{\nu + 1} = \theta^{\nu + \frac{1}{2}}$. {We remark that under such an update scheme, since the size of $\varepsilon$-active index set grows
exponentially along with the number of pieces of PADR, the iterate point would be more easily rejected for PADR with larger pieces.}  

\begin{algorithm}[htbp]
    \caption{ESMM algorithm %
    }\label{ALGO:unconstrained case}
    \SetAlgoLined{}
    \KwIn{initial point $\theta^0$, iteration length $T$, parameters $\{\eta, \varepsilon\}$, and increasing sequence $\{N_{\nu}\}_\nu$;}
    \For{$\nu = 0$ {\rm \bf to} $T-1$}{      
        generate $\Xi_{\nu} \coloneqq \{\wtxi_\nu^s\}_{s=1}^{N_{\nu}}$ independently and uniformly from $\Xi$\;
        construct the surrogate function 
        $\widehat{F}_{\nu}(\theta; \theta^\nu, \wtbI_{\nu}^\varepsilon) \coloneqq \frac{1}{N_{\nu}} \sum_{s =1}^{N_{\nu}} \, \widehat F(\theta, \wtxi_\nu^s; \theta^\nu, \wtI^\varepsilon_{\nu,s})$ where $\wtbI_{\nu}^\varepsilon\coloneqq \{\wtI^\varepsilon_{\nu,s}\}_{s=1}^{N_{\nu}}$ and $\wtI^\varepsilon_{\nu,s}$ is uniformly generated from $ \cI^\varepsilon(\theta^\nu, \wtxi_{\nu}^s)$\;
        compute 
        $\theta^{\nu + \frac{1}{2}} \in 
        \underset{\theta \in \Theta}{\mathop{\arg\min}}
        \left\{\widehat{F}_{\nu}(\theta; \theta^\nu, \wtbI_{\nu}^\varepsilon) + \frac{\eta}{2}\norm{\theta - \theta^\nu}^2\right\}$\;
        update the iterate
            $\theta^{\nu + 1} = 
            \begin{cases}
                \theta^{\nu + \frac{1}{2}}, 
                & \text{if}\ \widehat{F}_{\nu}(\theta^{\nu + \frac{1}{2}}; \theta^\nu, \wtbI_{\nu}^\varepsilon) + \frac{\eta}{2}\lVert\theta^{\nu + \frac{1}{2}} - \theta^{\nu}\rVert^2 \leq F_{\nu}(\theta^\nu),\\
                \theta^{\nu}, 
                & \text{otherwise.}
            \end{cases}$
    }
    \KwOut{$\widetilde{\theta}^{\,T}$ uniformly sampled from the iterates $\{\theta^\nu\}_{\nu = 0}^{T-1}$.}
\end{algorithm}

Next, we establish the convergence analysis of Algorithm~\ref{ALGO:unconstrained case} in terms of (composite strong) d-stationarity for problem~\eqref{eq:algo:formulation:finite sum primal problem}. Due to the combination of sampling and enhanced technique, we introduce the averaged residual for convergence analysis, which is technically different from the analysis in \cite{pang2017computing,cui2018composite,lu2019enhanced}.  We present the main theoretical results with essential discussions while leaving technical details in Online Appendix~\ref{EC:sec:proofs for Section Convergence}. 
Throughout the analysis, we make several blanket assumptions: the outer function $\varphi$ is directionally differentiable and locally Lipschitz continuous, $\Theta$ is a compact convex set, and the family of enhanced surrogate functions  $\{\widehat{F}(\cdot, \xi; \bar{\theta}, I)\}_{I \in \mathcal I^\varepsilon(\bar{\theta},\, \xi)}$ satisfies Assumptions~\ref{assum:B1:Touching}-\ref{assum:B4:minSurr dd-consistency} for each $\xi \in \Xi$ and any $\bar \theta \in \Theta$. 

With $\Xi = \{\xi^s\}_{s\in[n]}=\{(x^s, y^s)\}_{s\in[n]}$,  let $\cI^\varepsilon(\bar{\theta})$ denote the set of all $\varepsilon$-active index combinations of $\{(g(\bar{\theta}, \xi),h(\bar{\theta}, \xi)\}_{\xi \in \Xi}$, i.e.,
\[
    \cI^\varepsilon(\bar{\theta}) \coloneqq \Big\{
        {\bI}^\varepsilon = \{(I_{1,s}^\varepsilon,I_{2,s}^\varepsilon)\}_{s=1}^n
        \mid I^\varepsilon_{1,s} \in \cI_g^\varepsilon(\bar{\theta}, \xi^s), I^\varepsilon_{2,s} \in \cI_h^\varepsilon(\bar{\theta}, \xi^s),
     \,\forall\,s=1, \ldots, n
    \Big\},
\] 
with $|\cI^\varepsilon(\bar{\theta})|$ %
representing the total number of combinations of $\varepsilon$-active index selection for the ERM objective function $\frac{1}{n} \sum_{s=1}^n F(\theta, \xi^s)$.  %
For the purpose of analysis, we introduce the surrogate function for such an ERM function, and adjusted proximal mapping point which is motivated by the iterative update in step 5 of Algorithm~\ref{ALGO:unconstrained case}.  At $\bar{\theta} \in \Theta$ with  $\bI^\varepsilon = \{I_s^\varepsilon %
\}_{s\in[n]} \in \cI^\varepsilon(\bar{\theta})$, let  
\begin{equation}\label{eq:mapping:barP}
    \begin{array}{ll}
    \widehat F(\cdot; \bar{\theta}, \bI^\varepsilon) & \coloneqq  \displaystyle{\frac{1}{n}} \sum_{s=1}^n \widehat F(\theta, \xi^s; \bar{\theta}, I^\varepsilon_s)\\[0.2in]
    P^{\rho}_{\bI^\varepsilon}(\bar{\theta}) & \coloneqq 
    \begin{cases}
        \mathrm{prox}_{\widehat{F}/\rho}(\bar{\theta},  \bI^\varepsilon) & \text{if}\ 
        \min_{\theta \in \Theta}\,\left\{\widehat{F}(\theta; \bar{\theta}, \bI^\varepsilon) + \frac{\rho}{2}\norm{\theta - \bar{\theta}}^2 \right\} \leq F(\bar{\theta}),\\
        \bar{\theta} & \text{otherwise}.
    \end{cases}
    \end{array}
\end{equation}
where the proximal map $
    \mathrm{prox}_{ \widehat{F}/\rho}(\bar{\theta}, \bI^\varepsilon) \coloneqq
    \mathop{\arg\min}_{\theta\in \Theta}
    \left\{\widehat{F}(\theta; \bar{\theta}, \bI^\varepsilon) + \frac{\rho}{2}\norm{\theta - \bar{\theta}}^2\right\}$. 
With the above prerequisite, we introduce the averaged residual 
\[
    \bar{r}^{\varepsilon, \rho}(\bar{\theta}) 
    \coloneqq \frac{1}{\abs{\cI^\varepsilon(\bar{\theta})}}\sum_{\bI^\varepsilon \in \cI^\varepsilon(\bar{\theta})}\norm{\bar{\theta} - P^{\rho}_{\bI^\varepsilon}(\bar{\theta})},
\]
which plays a central role in the follow-up convergence analysis. Under the sampling size sequence $\{N_{\nu}\}$, we present the convergence result of $\bar{r}^{\varepsilon,\rho}$ in expectation in the proposition below, with the proof provided in Online Appendix~\ref{proof:prop::algo:unconstrained:convergence rate of M mapping}.

\begin{proposition}\label{prop:convergence:convergence rate of bar-r}
    Let $\widetilde{\theta}^{\,T}$ be the output of Algorithm~\ref{ALGO:unconstrained case} after $T$ iterations with $\varepsilon \geq 0$, $\eta \geq 0$, and $\{N_{\nu}\}_{\nu=1}^T$. There exists $C_1 > 0$ such that for any $\rho > \eta$,
    \[
        \big( \EE\,[\bar{r}^{\varepsilon, \rho}(\widetilde{\theta}^{\,T})]\, \big)^2
        \leq \frac{
            F(\theta^0) - \min_{\theta\in\Theta} F(\theta) + C_1\cdot \sum_{\nu = 0}^{T-1}\sqrt{\log N_{\nu} / N_{\nu}}
            }{\frac{\rho - \eta}{2}\cdot T}.
    \]
    In particular, with $\varepsilon = 0$ and $\eta > 0$, there exists $C_2 > 0$ such that for any $\rho > \eta$,
    \[
        \big( \EE\,[\bar{r}^{0, \rho}(\widetilde{\theta}^{\,T})]\, \big)^2
        \leq \frac{
            F(\theta^0) - \min_{\theta\in\Theta} F(\theta) + C_2\cdot \sum_{\nu = 0}^{T-1}(1/N_{\nu})
            }{\frac{\rho - \eta}{2}\cdot T}.
    \]
\end{proposition}

By relating $\bar r^{\varepsilon, \rho}$ to the proximal mapping gap with respect to $\widehat F_{\min}^\varepsilon(\cdot; \bar{\theta})$, with $N_{\nu} = \cO(\nu^\alpha)$ for some $\alpha > 0$, we can establish the convergence of Algorithm~\ref{ALGO:unconstrained case} in terms of composite strong d-stationarity.

\begin{theorem}[Asymptotic convergence of ESMM]\label{thm:convergence:convergence results}
    For problem~\eqref{eq:algo:formulation:finite sum primal problem}, let $\{\widetilde{\theta}^{\,t}\}_t$ be a sequence of points produced by Algorithm~\ref{ALGO:unconstrained case} with $\varepsilon > 0$, $\eta \geq 0$, and $N_{\nu} = \cO(\nu^\alpha)$ for some $\alpha > 0$. Then any limit point of $\{\widetilde{\theta}^{\,t}\}_t$ which must exist is a composite $(\varepsilon', \rho)$-strongly d-stationary point of problem~\eqref{eq:algo:formulation:finite sum primal problem} with probability 1 for any $\varepsilon' \in [0, \varepsilon)$ and $\rho \in (\eta, +\infty)$
\end{theorem}

Theorem~\ref{thm:convergence:convergence results} shows the asymptotic convergence of the ESMM algorithm to composite strong d-stationarity that is stronger than d-stationarity by Proposition~\ref{prop:convergence:relationship-S*-Sds}.
{
Moreover, Proposition~\ref{prop:EC:relationship-S* equiv Sder} in Online Appendix~\ref{EC:subsec:proofs for Section4-1} shows that $\cS^d_{\varepsilon, \rho}$ coincides with $\cS^*$ under a sufficiently large $\varepsilon$ and a regularity assumption on $F$.
Consequently, any limit point of ${\widetilde \theta^{,t}}_t$ generated by the algorithm is a global minimizer of~\eqref{eq:algo:formulation:finite sum primal problem} with probability one.}
This indicates that the PADR-based ERM model together with the ESMM algorithm is asymptotically optimal to the original problem~\eqref{eq:intro:contextalSP:Full-E}. 
However, evidenced by the numerical implementations, when $\varepsilon$ is large, the candidate solutions produced in step~4 of Algorithm~\ref{ALGO:unconstrained case} may be rejected by the update rule in step~5 very often, which results in inefficient computation and slow convergence. Such a computational issue is resolved by a shrinking-$\varepsilon$ strategy which will be specified in  Section~\ref{section:experiments}.

Particularly, when $\varphi$ is piecewise affine and $\Theta$ is a polytope, we develop error-bound analysis and establish the nonasymptotic convergence to d-stationarity (also local optimality in this case) of Algorithm~\ref{ALGO:unconstrained case} by leveraging the special property that at all points within a polytope that are not d-stationary, 
the piecewise affine objective has descent rates bounded away from 0. 
Furthermore, the nonasymptotic convergence in this case is also assured with $\varepsilon=0$, in which the algorithm reduces to the standard SMM with randomized surrogation-index selection. %
Details of the error-bound analysis and convergence proofs can be found in Online Appendix~\ref{EC:subsec:analytical details on convergence}.
 
\begin{theorem}[Nonasymptotic convergence of ESMM]\label{thm:convergence:convergence results:PA} 
    For problem~\eqref{eq:algo:formulation:finite sum primal problem} with $\varphi$ being piecewise affine and $\Theta$ being a polytope, under the same settings of Algorithm~\ref{ALGO:unconstrained case} specified in Theorem~\ref{thm:convergence:convergence results} except that $\varepsilon \geq 0$,
    $\EE\,[\distSd{\widetilde \theta^{\, T}}] = \widetilde{\cO}\Big(T^{-\frac{1}{2}\min\left\{\frac{\alpha}{2},1\right\}}\Big)$.
    Moreover, with $\varepsilon = 0$ and $\eta > 0$, $\EE\,[\distSd{\widetilde \theta^{\,T}}] = \cO(T^{-\frac{1}{2}\min\left\{\alpha, 1\right\}})$, which is up to a logarithmic factor when $\alpha = 1$.
\end{theorem}
\begin{corollary}\label{coro:algo:opt-iteration-complexity}
    Under the settings specified in Theorem~\ref{thm:convergence:convergence results:PA},
    to achieve $\EE\,[\distSd{\widetilde \theta^{\,T}}] \leq \epsilon$, the optimal iteration complexity of Algorithm~\ref{ALGO:unconstrained case} is $\widetilde{\cO}(\epsilon^{-2})$, which holds with a sampling rate $N_{\nu} = \cO(\nu^\alpha)$ where $\alpha \geq 2$.
\end{corollary}
The nonasymptotic convergence to d-stationarity for general cases (including~\ref{enum:algo:phi2} and~\ref{enum:algo:phi3}) can be obtained under additional assumptions on the objective $F$ and its surrogate functions, and we refer interested readers to the relative results with technical details in Online Appendix~\ref{EC:sec:nonasymptotic convergence for general cases}. 
We recognize that it is possible, albeit challenging, to achieve nonasymptotic results in terms of composite strong d-stationarity with carefully designed residuals and error-bound analysis, which we will leave for future study. 

Finally, we make several remarks on the convergence results.

\begin{remark}
The enhanced $\varepsilon$-active surrogation with positive $\varepsilon$ is essential for general cases. 
As noted in \citet{liu2022solving}, due to the composite structure of objective functions, the SMM algorithm with $\varepsilon=0$ may yield convergence to compound criticality, which is 
weaker than d-stationarity. This is also supported by Example~\ref{eg:plot of iterative process} in Section~\ref{sec:standardMM}.
\end{remark}

\begin{remark}
By Corollary~\ref{coro:algo:opt-iteration-complexity}, the iteration complexity of Algorithm~\ref{ALGO:unconstrained case} is similar to that of the stochastic DC algorithm in \citet{le2022stochastic} with incremental sampling for differentiable stochastic DC programs, provided that convex subproblems are solved to the global optimum.
Notably, the optimal complexity of Algorithm~\ref{ALGO:unconstrained case} is achieved with the sampling rate that is faster than $\cO(\nu)$ in \citet{le2022stochastic} due to the nondifferentiable composite structure, and the convergence results of our algorithm hold for a much wider range of nonconvex and nondifferentiable composite problems.
\end{remark}

\begin{remark}

With minor modifications of the convergence analysis, we may obtain the same convergence results for an inexact version of Algorithm~\ref{ALGO:unconstrained case}, in which a $\delta_\nu$-suboptimal solution $\theta^{\nu + \frac{1}{2}}$ is computed in Step~4 and only accepted in Step~5 when it satisfies
$\widehat{F}_{\nu}(\theta^{\nu + \frac{1}{2}}; \theta^\nu, \wtbI_{\nu}^\varepsilon) + \frac{\eta}{2}\lVert\theta^{\nu + \frac{1}{2}} - \theta^{\nu}\rVert^2 \leq F_{\nu}(\theta^\nu) + \delta_\nu$ with $\delta_\nu = \cO(\nu^{-1})$.
In such a plausible extension, it is viable to design a two-loop algorithmic scheme addressing the computational complexity in solving convex subproblems. Together with our iteration complexity result, 
it may enable characterization of the actual work complexity of the entire algorithm similar to \citet{pasupathy2021adaptive} and complexity comparisons with other first-order algorithms for nonconvex SP.
\end{remark}

\begin{remark}
    Since the output $\widetilde \theta^{\,t}$ of Algorithm~\ref{ALGO:unconstrained case} is uniformly sampled from the iterate sequence $\{\theta^\nu\}_{\nu = 0}^{t-1}$,
   by the almost sure convergence of $\bar{r}^{\varepsilon, \rho}(\widetilde \theta^{\,t})$ to 0 in Proposition~\ref{prop:convergence:convergence rate of bar-r}, we can easily obtain that 
    there exists a subsequence $\{\theta^\nu\}_{\nu \in \kappa}$ with $\lim_{\nu ( \in \kappa) \to \infty} \bar{r}^{\varepsilon, \rho}(\theta^\nu) =0$ almost surely. With similar convergence analysis, we can obtain that the limit point of $\{\theta^\nu\}_{\nu \in \kappa}$ satisfies the same almost sure asymptotic convergence properties as the limit points of $\{\widetilde \theta^{\,t}\}_t$ presented in Theorem~\ref{thm:convergence:convergence results}.
\end{remark}

\section{Numerical Experiments}\label{section:experiments}
In this section, we compare the numerical performance of the PADR-based ERM method with several benchmarks including DR-based methods, PO methods, and multiple prescriptive methods which are listed in Table~\ref{tab:exp:methods}. The experiments are conducted on Newsvendor problems with both simulated and real data and a product placement problem on a practical scale.

For the simulated data, we set the true model of $Y$ as $Y = \overline{Y}(X) + \widetilde{\epsilon}$, where $\overline{Y}(\cdot)$ represents the mean model, and $\widetilde{\epsilon}$ denotes a random variable with a standard normal distribution independent of $X$.
Throughout all experiments with simulated data, we generate feature data uniformly from $[-1, 1]^p$, and ensure the mean model satisfies $\min_{x\in[-1, 1]^p}\overline{Y}(x)\geq 5$ which keeps all samples of $Y$ positive.
Unless otherwise noted, both training and test datasets contain 1000 samples.

\begin{table}[htbp]
    \TABLE
    {Methods implemented in numerical experiments\label{tab:exp:methods}}
      {\resizebox{0.99\columnwidth}{!}{\begin{tabular}{p{12.5em}p{40em}}
      \toprule
      Abbreviation & Explanation of method \\
      \midrule
      PADR  & Piecewise-affine decision rule method in this paper \\
      RKHS-DR  & Kernel optimizer prediction method from~\citet{bertsimas2022data} \\
      NN-DR & Neural network-based decision rule method with two or three fully connected hidden layers\\
      LDR   & Linear decision rule method from \citet{ban2019big}\\
      kNN, KERNEL, CART, RF & Prescriptive methods from~\citet{bertsimas2020predictive} based on $k$-nearest neighbiors, kernels, regression trees, and random forests respectively \\
      SOF   & Prescriptive method from~\citet{kallus2022stochastic} based on stochastic optimization forests \\
      PO-PA, PO-L & PO methods with piecewise-affine prediction model and linear prediction model respectively \\
      SIMOPT & Simulated optimal decision using true conditional distribution \\
      \bottomrule
      \end{tabular}}}
    {}
  \end{table}%

We refer to the PADR model with two max-affine components consisting of $K_1$ and $K_2$ pieces respectively as PADR$(K_1, K_2)$.
{The hyperparameters $K_1$ and $K_2$ are searched within the range [0, 40].}
Regarding the ESMM algorithm for learning the PADR, initial points are generated uniformly from $[-50, 50]^{d(K_1+K_2)(p+1)}$, and the sampling sequence $\{N_{\nu}\}_\nu$ follows $N_{\nu} = 40 \nu$. The parameter $\varepsilon$ in Algorithm~\ref{ALGO:unconstrained case} is set to be either constant or shrinking, the latter case of which is referred to as the shrinking-$\varepsilon$ strategy. Specifically, the shrinking strategy sets the parameter a constant value $\varepsilon_0$ for the first $T_0$ iterations and then a reduced value $\varepsilon_1(< \varepsilon_0)$ for the remaining ones.
The ESMM algorithm produces the iterate point obtaining the minimum ERM cost as the output after $T$ iterations. 
For each instance, the algorithm is executed for $r$ rounds with random initialization, and we report the best performance unless otherwise noted. We set $T=10$ and $r=10$ in Sections~\ref{subsec:exp:Newsvendor} and~\ref{subsec:exp:nonconvex Newsvendor}, and $T=20$ and $r=5$ in Section~\ref{subsec:exp:pp}. The reported running time of PADR is the total time of all rounds. {The specific construction of surrogate functions in ESMM for each problem is provided in Online Appendix~\ref{EC:sec:surrogation construction for exps}.}

We clarify settings for benchmark methods listed in Table~\ref{tab:exp:methods}. 
{For NN-DR, we use ReLU function ($\max\{\cdot, 0\}$) as activation, and the number of nodes in each hidden layer is set to half the size of the preceding layer. In each experiment, the model is trained for 100 epochs over 3 independent runs using the Adam optimizer~\citep{kingma2014adam} to minimize the empirical cost with $\ell_1$-regularization.}
We use a Gaussian kernel for kernel methods. 
For RKHS-DR, we adjust the diagonal of the kernel by values no larger than $10^{-3}$ so as to resolve numerical issues without undermining performance. 
For PO methods, we use $\ell_2$ loss in the prediction stage and solve the regression problem of PO-PA using the ESMM algorithm. We set the uniform scale limit $\mu=50$ in Sections~\ref{subsec:exp:Newsvendor} and~\ref{subsec:exp:nonconvex Newsvendor}, and $\mu = 100$ in Section~\ref{subsec:exp:pp}.
The hyperparameters of all benchmarks methods listed in Table~\ref{tab:exp:method hyperparameters} are selected via a standard validation process.
The performance of LDR, PO-L, and KERNEL is dominated by PADR, PO-PA, and other prescriptive methods in most cases, and thus their results are omitted from the report if necessary to avoid the distraction.

\begin{table}[htbp]
    \TABLE
    {Hyperparameters of each method\label{tab:exp:method hyperparameters}}
    {\resizebox{0.99\columnwidth}{!}{\begin{tabular}{p{12.5em}p{40em}}
        \toprule
        Methods & Hyperparameters \\
        \midrule
        PADR  & number of pieces (PADR); sampling rate parameters, shrinking-$\varepsilon$ strategy, and proximal parameter (ESMM) \\
        RKHS-DR  & kernel bandwidth, regularization parameter, penalty parameter \\
        NN-DR & first-hidden-layer size, regularization parameter, learning rate, weight decay, batch size\\
        LDR & N/A \\
        kNN, KERNEL, CART, RF & number of nearest neighbors (kNN); kernel bandwidth, regularization parameter (KERNEL); depth, minimum samples to split, minimum samples per leaf (CART and RF); \\
        SOF   & choice of apx-soln or apx-risk method in \cite{kallus2022stochastic}, depth, minimum samples to split, minimum samples per leaf \\
        PO-PA & number of pieces (prediction model); sampling rate parameters, shrinking-$\varepsilon$ strategy, and proximal parameter (ESMM solving $\ell_2$ regression) \\
        PO-L & N/A \\
        SIMOPT & N/A \\
        \bottomrule
        \end{tabular}}}
    {}
  \end{table}%

All methods are implemented in Python 3.9.7 on a MacBook Pro with a macOS Monterey system and an Apple M1 Pro chip.
{NN-DR is implemented in Pytorch version 2.2.2.}
All programs involved are built and solved by the optimizer Gurobi, version 9.5.0. We provide the code on
GitHub:  {\url{https://github.com/zhangyiy0223/PADR}. 

\subsection{Newsvendor problem: simulated data}\label{subsec:exp:Newsvendor}

We consider a Newsvendor problem with covariate information $X=x$ as follows,
\begin{equation}
    \label{eq:exp-unconstr-nv}
    \min_{z \geq 0} \quad \EE\left[
        c_b(Y - z)_+ + c_h(z - Y)_+ \mid X = x
    \right],
\end{equation}
where $c_b$ and $c_h$ are the unit back-order and holding costs respectively, and $(\cdot)_+$ denotes $\max\{\cdot, 0\}$. 
We set $(c_b, c_h) = (8, 2)$ and consider the following max-affine demand mean model with $\nonlinearity = 1$ unless otherwise specified:

\begin{equation}\label{eq:exp:max-affine model}
    \overline{Y}(X) = \nonlinearity\max\{5X_1 - 10X_2, -10X_1+ 5X_2, 15X_1\} + 10.
\end{equation}

\subsubsection{Basic results.}\label{subsubsec:NV:basic}
\begin{figure}[bp]
    \FIGURE
    {
      {
      \subfloat[Basic results $\phantom{xxxxxx}$]{
          \includegraphics[width=0.426\textwidth]{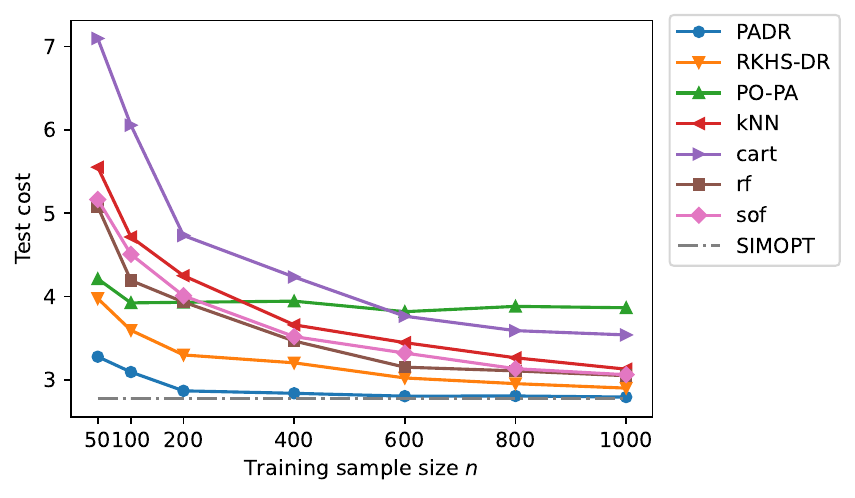}\label{subfig:unconstr-basic}
      }
      \subfloat[Running time $\phantom{xxxxx}$]{
          \includegraphics[width=0.46\textwidth]{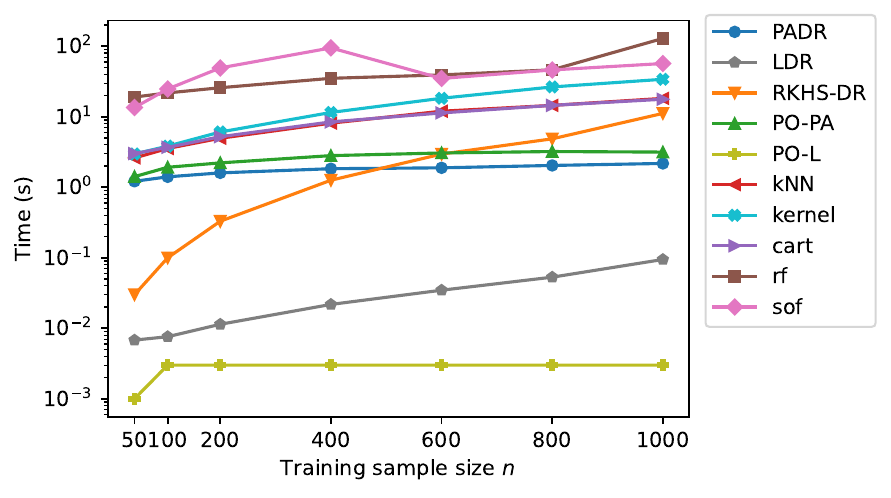}\label{subfig:unconstr-basic-time}
      }
      }
    }
    {Newsvendor problem: basic results and running time\label{fig:exp:unconstr-basic results}}
    {PO-L, KERNEL, and LDR are excluded in (a) for readability, whose costs are around 14, 11, and 9 respectively.}
  \end{figure}
Since the optimal DR in this case is a max-affine function with three pieces,  PADR$(K_1,K_2)$ with $K_1 \geq 3$ has no approximation error in principle. Figure~\ref{subfig:unconstr-basic} shows that the costs of PADR with $(K_1, K_2) = (3, 0)$ converge to SIMOPT as the training sample size $n$ increases from 50 to 1000.
RKHS-DR converges slower than PADR; additionally, its performance with $n=1000$ is worse than that of PADR with $n=400$. 
{NN-DR does not outperform PADR and RKHS-DR, possibly due to the instability of using gradient-based methods to optimize a nonconvex and nonsmooth objective function.}
Prescriptive methods, including  kNN, CART, RF, and SOF, are dominated by the PADR and RKHS-DR methods. %
For PO-PA, we observe that the prediction model provided by the PO-PA method is quite close to the true model $\overline{Y}(\cdot)$, so its substantial gap to SIMOPT could be mainly due to inconsistency between the prediction loss and Newsvendor cost with imbalanced $c_b$ and $c_h$. Numerical results with the parameter ratio $c_{b}:c_{h}$ varying from 1:9 to 9:1 exhibit similar patterns and thus are omitted.
Figure~\ref{subfig:unconstr-basic-time} shows the training time of each method.
The training time of RKHS-DR grows roughly exponentially when $n \geq 200$, whereas the time for PA-related methods remains stable and significantly shorter when $n=1000$, which is mainly credited to the enhanced and sampling techniques used in the ESMM algorithm. 
{NN-DR incurs a higher computational cost than most other methods on this problem.}%

\begin{figure}[tbp]
    \FIGURE
    {
      {
      \subfloat[$\varepsilon$-Strategies $\phantom{xxxxxx}$]{
          \includegraphics[width=0.33\textwidth]{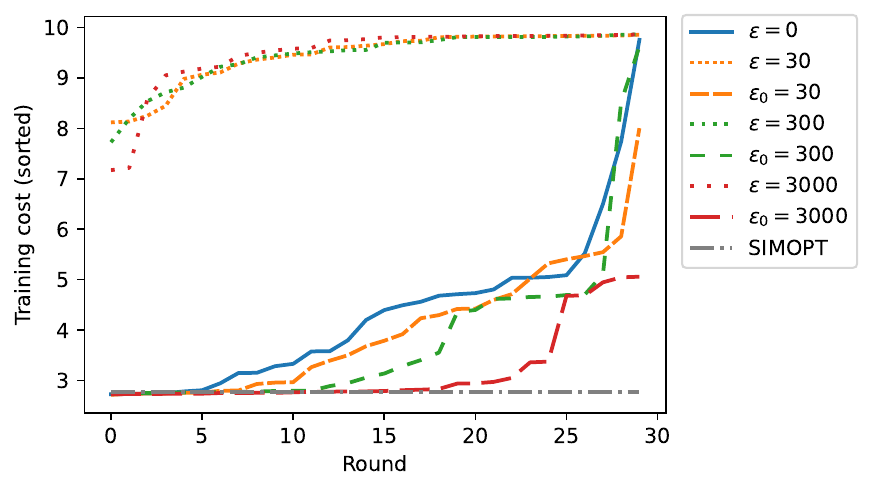}\label{subfig:shrinkingeps}
      }
      \subfloat[Sampling technique $\phantom{xxxxxxxxx}$]{
          \includegraphics[width=0.39\textwidth]{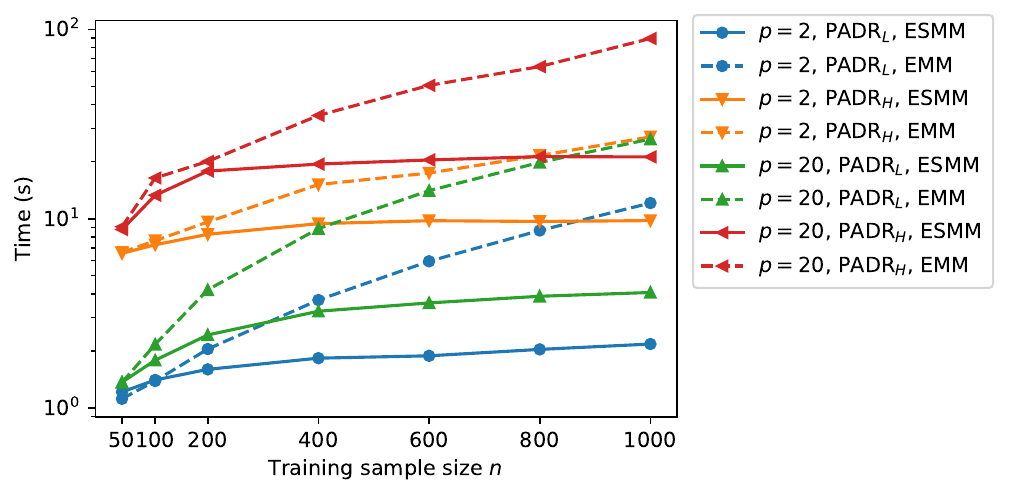}\label{subfig:sampling}
      }
      \subfloat[Time across $K_1$ and $K_2$]{
          \includegraphics[width=0.265\textwidth]{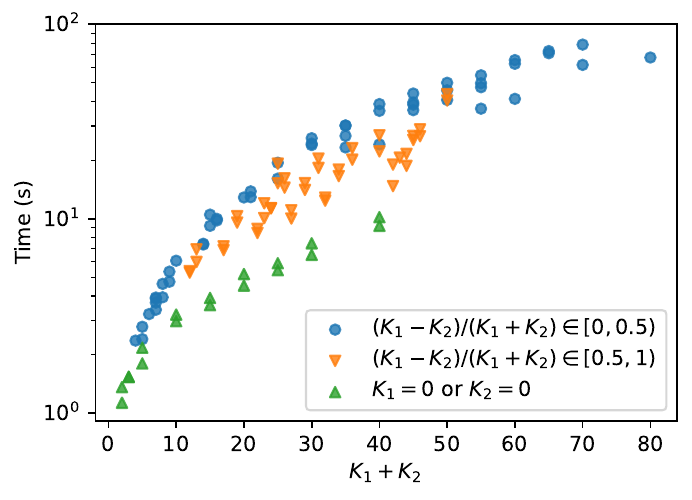}\label{subfig:timescalingwithK}
      }
      }
    }
    {Performance of ESMM under various settings\label{fig:settingsforESMM}}
    {In (a), solid and dotted lines represent constant-$\varepsilon$ strategies; dashed lines represent shrinking-$\varepsilon$ strategies with $T_0 = 3$ and $\varepsilon_1=0$. In (b), 
    experiments for $p=20$ are conducted on datasets using the sparse demand model introduced in Section~\ref{subsubsec:NV:varyingfeature}. The y-axis is set on a log scale. {In (c), each point represents the running time of PADR with $K_1$ and $K_2$ varying from 0 to 40 for the Newsvendor problem with $p = 2$ and $n = 1000$.}}
\end{figure}

In Figure~\ref{fig:settingsforESMM}, we provide results of the ESMM algorithm under different parameter settings that illustrate the effect of the enhanced sampling strategies. 
Figure~\ref{subfig:shrinkingeps} displays curves of sorted training cost of PADR provided by the ESMM algorithm under the constant-$\varepsilon$ and shrinking-$\varepsilon$ strategies over 30 rounds.
Under the constant-$\varepsilon$ strategy, the algorithm converges to stationary solutions which can be strengthened to global optimality with $\varepsilon$ being sufficiently large. However, Figure~\ref{subfig:shrinkingeps} shows that the constant-$\varepsilon$ strategy typically yields poor performance over all 30 rounds with 10 iterations. We conjecture that this is because when the sample size is large, the algorithm requires engaging with the random selection among a huge number of active indices and thus hardly makes a valid update at each iteration. In contrast, the shrinking-$\varepsilon$ strategy enables flexible adjustments of active indices for early iterations with a large positive value of $\varepsilon$ and produces iterate points serving as a warm start; afterward, it accelerates the convergence with a decreased value of $\varepsilon$. Figure~\ref{subfig:shrinkingeps} shows that the shrinking-$\varepsilon$ strategy has significantly better performance which attains the global optimum among roughly half of the total rounds.
 
Figure~\ref{subfig:sampling} reports the running time of PADR using the ESMM algorithm and the enhanced MM algorithm (EMM) without sampling (including all samples in the subproblem of each iteration). 
We use PADR(3,0) and PADR(6,4) to represent models with low and high complexity, denoted by PADR$_L$ and PADR$_H$, respectively.
The sampling technique significantly saves computational costs without compromising performance. The time consumption of ESMM is primarily determined by the sampling rate and remains relatively stable when $n \geq 400$, while the time required for EMM increases exponentially with the training sample size. The benefit of the sampling technique becomes more pronounced as the decision rule complexity increases. The time required for EMM is approximately 4 times longer than that for ESMM when $n=1000$.

{Figure~\ref{subfig:timescalingwithK} illustrates the running time of PADR for multiple settings of $(K_1, K_2)$ with each value ranging from 1 to 40. The running time for 10 rounds remains within 100 seconds when $\max\{K_1, K_2\} \leq 40$. 
An interesting observation is that, with the same total number of pieces (i.e., total number of parameters), the time consumption of PADR is relatively higher when its two components have more balanced piece numbers, and considerably lower when the decision rule consists of only one component (i.e., $K_1=0$ or $K_2=0$). This  can be explained by the complex structure of the surrogate function for the PADR with balanced piece numbers. 
}

\subsubsection{Varying feature dimensions.}\label{subsubsec:NV:varyingfeature}
We study the effect of the feature dimension of $X \in \mathbb R^p$ on the performance of PADR under sparse and dense demand models of $\overline{Y}$ following the settings in \citet{bertsimas2022data}. Specifically, in the sparse case, only the first two components of $X$ contribute to $\overline{Y}$ following the demand model~\eqref{eq:exp:max-affine model};
in the dense case, we replace the two feature variables in~\eqref{eq:exp:max-affine model} by the average of the first half of the features $\frac{2}{p} \sum_{i=1}^{p/2} X_i$ and the average of the other half respectively. 

\begin{figure}[htbp]
    \FIGURE
    {
    {
        \subfloat[Sparse demand model $\phantom{xxxxxx}$]{%
          \includegraphics[width=0.426\textwidth]{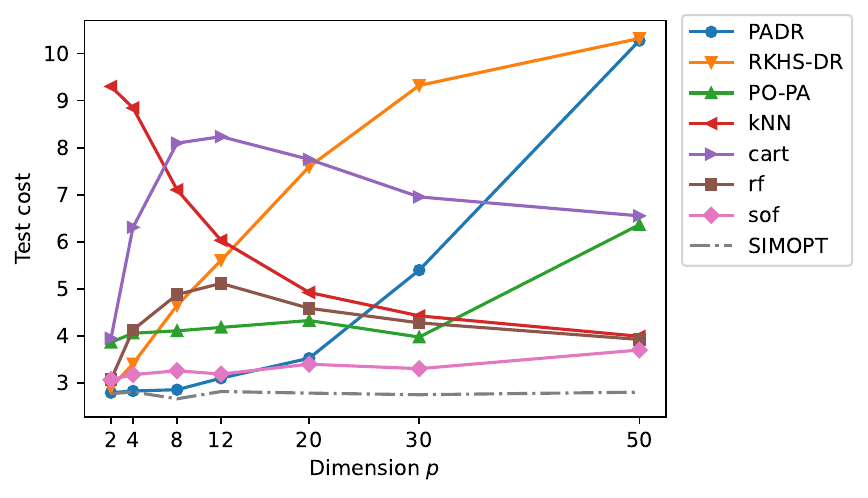}%
          \label{subfig:exp:UNC-MA-dim:sparse}
        }
        \subfloat[Dense demand model $\phantom{xxxxxxx}$]{%
          \includegraphics[width=0.41\textwidth]{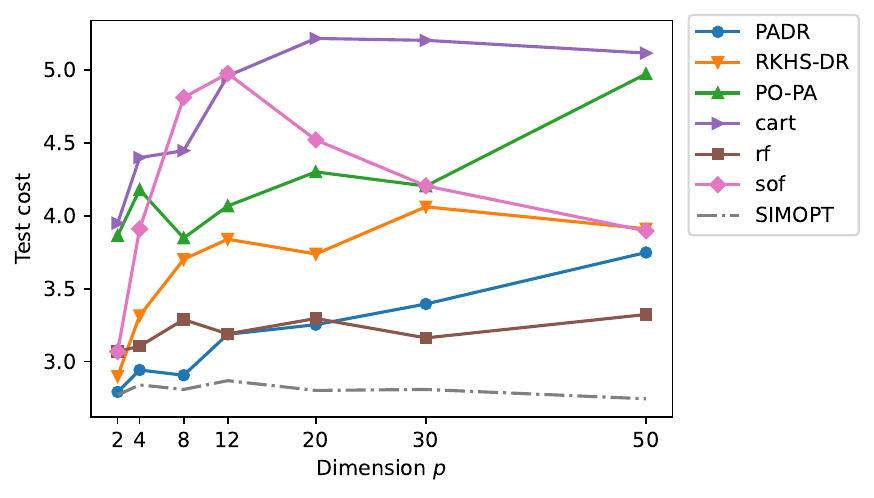}%
          \label{subfig:exp:UNC-MA-dim:dense}
        }
      }
    }
    {Newsvendor problem with varying feature dimensions\label{fig:exp:UNC-MA-dim}}
    {%
    For readability, kNN is excluded under the dense demand model, whose costs exceed 8 in most trials.
     }
  \end{figure}

\begin{figure}[htbp]
    \FIGURE
    {
    {
        \subfloat[Sparse demand model $\phantom{xxxxxx}$]{%
          \includegraphics[width=0.426\textwidth]{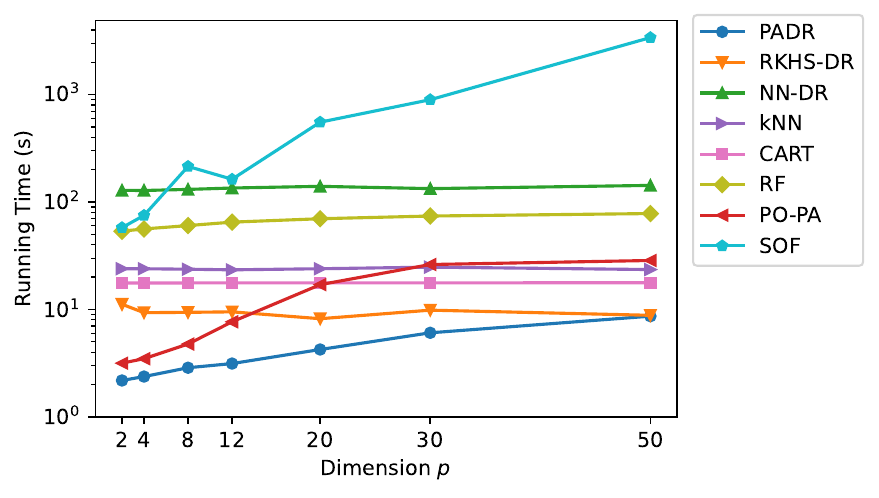}%
          \label{subfig:exp:UNC-MA-dim:time:sparse}
        }
        \subfloat[Dense demand model $\phantom{xxxxxxx}$]{%
          \includegraphics[width=0.41\textwidth]{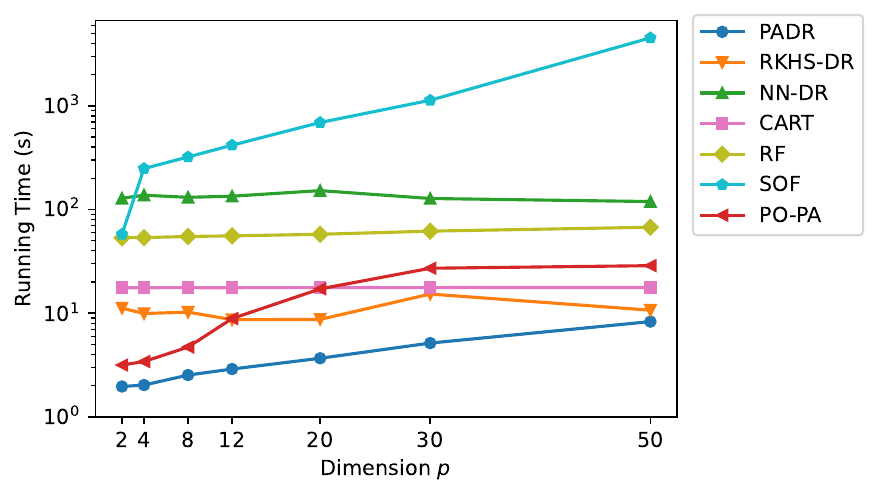}%
          \label{subfig:exp:UNC-MA-dim:time:dense}
        }
      }
    }
    {Running time in Newsvendor problem with varying feature dimensions\label{fig:exp:UNC-MA-dim:time}}
    {%
     }
  \end{figure}

{Figure~\ref{fig:exp:UNC-MA-dim} and Figure~\ref{fig:exp:UNC-MA-dim:time} report the test cost and running time, respectively.}
In the sparse case, forest-based prescriptive methods are robust to redundant features, as they only focus on informative features. It would be beneficial to employ an effective feature selection technique in the PADR method to improve the performance in the sparse case with a large feature dimension. 
In the dense case, PADR performs generally well, but its performance slightly deteriorates as $p$ increases due to growing estimation errors and optimization difficulties. We also observe that PADR with fewer pieces performs more stably in both cases despite the presence of approximation errors. 
{
When $p < 20$, PADR consistently achieves slightly better performance than benchmarks, while requiring substantially less training time.
For instance, at $p = 8$, PADR attains a 10\% lower test cost compared to NN-DR and SOF, with a 90\% less runtime than these methods.
This highlights the efficiency of PADR in learning decision rules in settings with moderate feature dimensions.
}

\subsubsection{Nonlinearity in the demand model.}\label{subsubsec:nonlinearity}
We study the effect of nonlinearity of the demand with respect to the feature variables by varying the coefficient $\nonlinearity$ in the demand model~\eqref{eq:exp:max-affine model} and report results in Figure~\ref{fig:exp:UNC-MA-nonlnr}. 
It indicates that the performance of PADR and PO-PA is irrelevant to the nonlinearity level of the demand model due to the strong approximation ability to the nonlinear structure. In contrast, the performance of RKHS-DR and prescriptive methods deteriorates significantly as the nonlinearity level increases. These methods require significantly more samples to achieve similar performance when the level of nonlinearity increases. 
{NN-DR shows the robustness towards non-linearity for the demand model with $p=20$; however, additional experiments in Figure~\ref{fig:exp:UNC-MA-nonlnr-NN} reporting relative cost gaps to SIMOPT for PADR, RKHS-DR, and NN-DR with multiple dimension settings, show that the performance of RKHS-DR deteriorates as $p$ increases, while NN-DR becomes less sensitive to nonlinearity as $p$ increases.
We speculate that when $p$ is small, the instability of gradient-based optimization becomes more pronounced, and this effect is further amplified under highly nonlinear demand functions.

To further test the effect of nonlinearity, we consider demand models with quadratic, cubic, and sine terms (in addition to a shared linear component). Results in Online Appendix~\ref{EC:sec:experiments} show consistent patterns across methods, validating PADR’s robustness to model misspecification under varying nonlinear structures.
}

\begin{figure}[htbp]
  \FIGURE
  {
    {
      \subfloat[$p=2 \phantom{xxxxxxx}$]{%
        \includegraphics[width=0.33\textwidth]{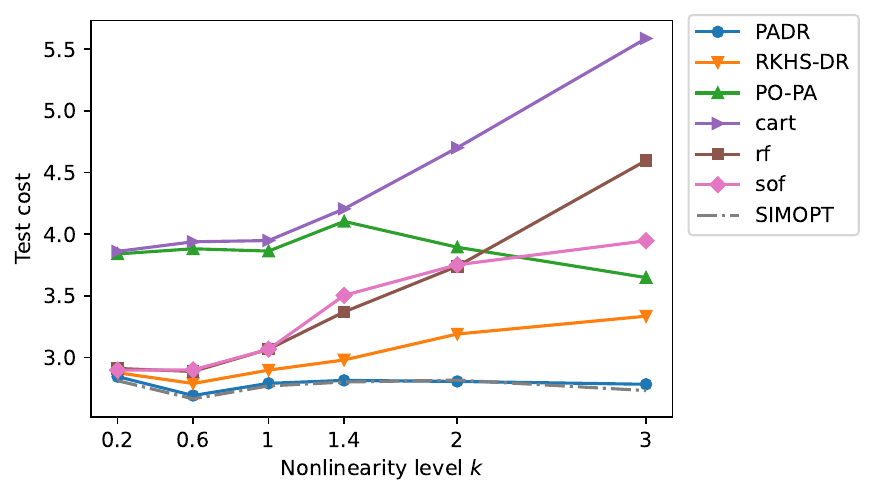}\label{fig:exp:UNC-MA-nonlnr-sub1}
      }
      \subfloat[$p=20$  $\phantom{xxxxxx}$]{%
        \includegraphics[width=0.33\textwidth]{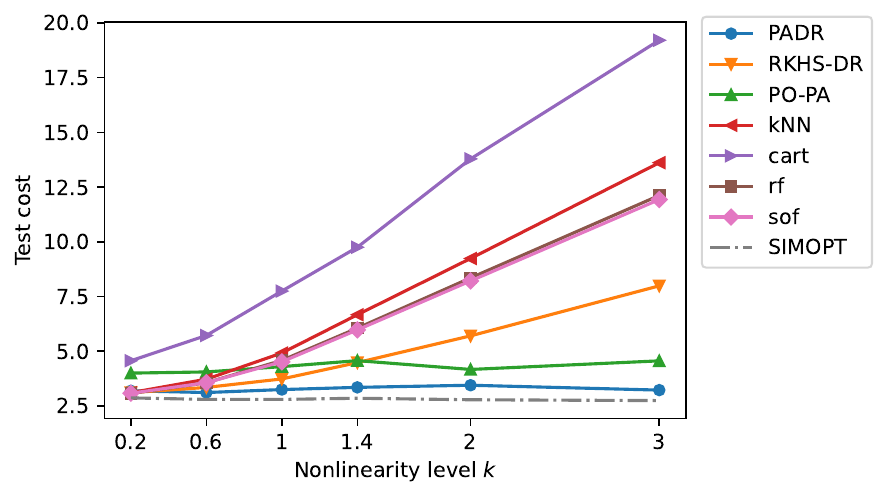}\label{fig:exp:UNC-MA-nonlnr-sub2}
      }
      \subfloat[Relative Cost Gap ($\omega=3$)]{%
        \includegraphics[width=0.265\textwidth]{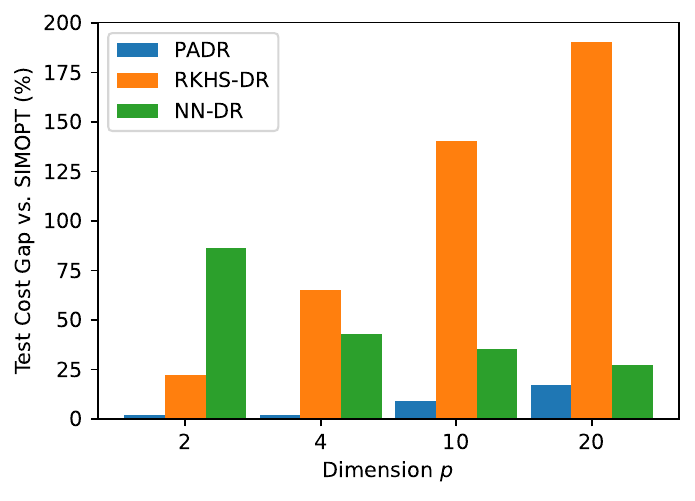}\label{fig:exp:UNC-MA-nonlnr-NN}
      }
    }
  }
  {Newsvendor problem with varying nonlinearity in demand model\label{fig:exp:UNC-MA-nonlnr}}
{In (b), the demand model is of the dense type. For readability, kNN, with costs increasing from 5 to 40, is excluded in (a).}
\end{figure}

\subsection{Newsvendor problem with nonconvex cost: real-world data}\label{subsec:exp:nonconvex Newsvendor}

In this section, we conduct experiments on a real-world dataset, which consists of order records for a specific model of mobile phone from a Chinese electronic product supplier. The order record, spanning from September 01, 2021 to November 30, 2023, contains order quantities and time stamps from retailers. 
In order to create a production plan, the supplier needs to determine the packing quantity in each day to fulfill the order requests.

We first organize the individual order records into daily order quantities. 
The aggregated daily order quantities have been rescaled for confidentiality purposes and also normalized for pretreatment.
The boxplots of daily orders in Figure~\ref{fig:exp:RW:demand distribution} show that order quantities are larger on weekdays and during the second quarter of each month, which aligns with retailers' ordering pattern.
\begin{figure}[htbp]
    \FIGURE
    {{
      \subfloat[Weekdays]{%
            \includegraphics[width=0.33\textwidth]{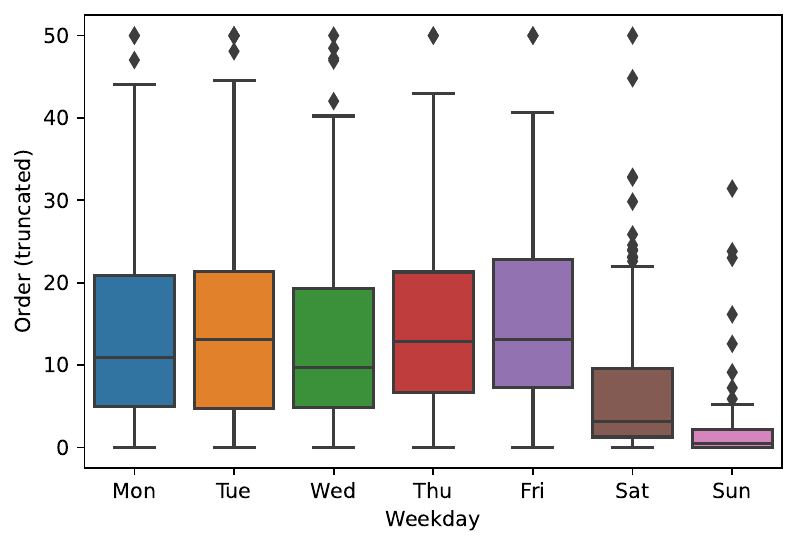}\label{subfig:rw-weekday}
      }
      \subfloat[Month days]{%
          \includegraphics[width=0.33\textwidth]{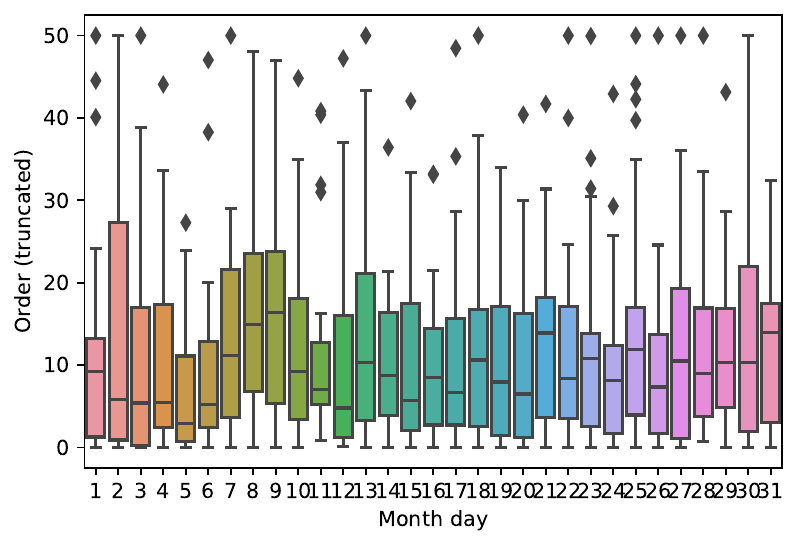}\label{subfig:rw-monthday}
      }
      \subfloat[Month quarters]{%
          \includegraphics[width=0.33\textwidth]{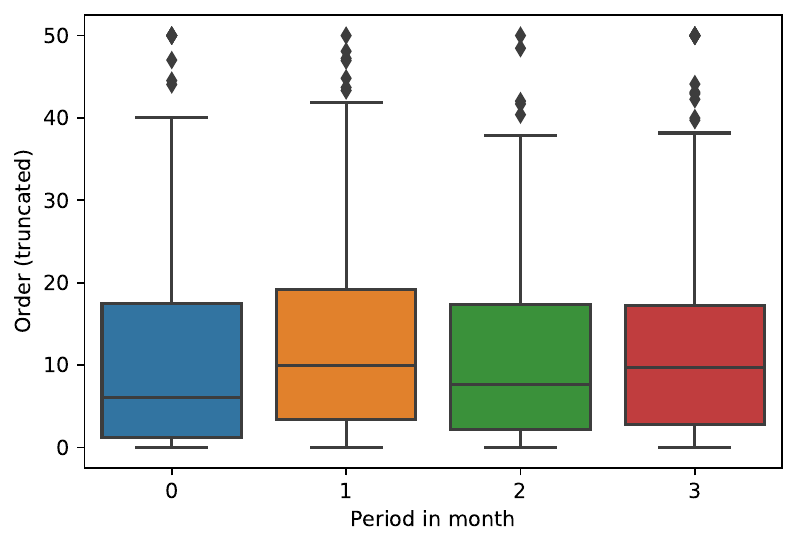}\label{subfig:rw-month quarters}
      }
    }}
    {Real-world dataset: boxplots of orders by time period\label{fig:exp:RW:demand distribution}}
    {%
    }
\end{figure}

We set observations from the first year as the training set and the remaining data as the test set.
Here we consider three types of features, namely weekday, month quarter, and order quantities in past days. 
The first two types of features are converted into dummy variables (7 for weekdays and 4 for month quarters) and set as basic features. 
For the third type of features, we incorporate order quantities in past days of different lengths as different instances in conducting the experiments. {Note that incorporating retailer's past order quantities into the feature vector introduces temporal dependence, violating the i.i.d. assumption in Assumption~\ref{assum:A1:iid dataset S}. 
This issue can be addressed using generalization bounds for non-stationary $\beta$-mixing processes, as established by \citet{kuznetsov2017generalization}.
Since only a fixed-length history is used in the feature, the resulting data sequence is $\beta$-mixing with finite memory, i.e., there exists a lag $a_0$ such that $\beta(a) = 0$ for all $a > a_0$. 
Under such settings, Theorem~1 and Corollary~1 from \citet{kuznetsov2017generalization} ensure that the generalization bound still holds, and thus the excess risk analysis remains valid.
}

To determine the packing
quantity in each day to fulfill the order requests, we formulate a Newsvendor problem with features listed above and report the numerical results in Table~\ref{tab:exp:RWresults-NV}. {The results show that PADR and NN-DR exhibit competitive performance. PADR outperforms NN-DR when $p \leq 15$, whereas NN-DR demonstrates greater robustness as the feature dimension increases.}

\begin{table}[htbp]
    \TABLE
    {Test costs for Newsvendor on a real-world dataset\label{tab:exp:RWresults-NV}}
    {
    \small
    \begin{tabular}{cccccccccc}
        \toprule
        \multicolumn{1}{c}{Added features} & \multirow{2}[2]{*}{$p$} & \multicolumn{3}{c}{DR methods} & \multicolumn{4}{c}{Prescriptive methods} & \multicolumn{1}{l}{PO methods} \\
        \cmidrule{3-10}
        (all include basic features) &       & PADR  & RKHS-DR  & NN-DR & kNN   & CART  & RF    & SOF   & PO-PA \\
        \midrule
         None & 11    & \redlog{55.95} & 56.10 & 56.74 & 56.06 & \bluelog{56.00} & \bluelog{56.00} & \bluelog{56.00} & 79.02 \\
         +  order quantities in past 1 day & 13    & \redlog{42.67} & 51.53 & \bluelog{45.13} & 50.25 & 61.72 & 49.54 & \bluelog{48.70} & 68.18 \\
         +  order quantities in past 3 days & 15    & \redlog{41.65} & 55.59 & 42.46 & 48.25 & 54.64 & \bluelog{41.74} & \bluelog{42.19} & 59.73 \\
         +  order quantities in past 5 days & 17    & \bluelog{41.40} & 65.17 & \bluelog{41.78} & 49.59 & 56.74 & \redlog{40.90} & 41.83 & 52.39 \\
         + order quantities in past 7 days & 19    & \bluelog{40.93} & 70.60 & \redlog{40.27} & 49.02 & 65.62 & 43.62 & \bluelog{43.51} & 54.52 \\
         + order quantities in past 14 days & 26    & \redlog{41.03} & 78.37 & \bluelog{41.04} & 46.81 & 58.39 & \bluelog{42.45} & 42.48 & 52.10 \\
         + order quantities in past 21 days & 33    & {42.43} & 80.24 & \redlog{40.12} & 46.61 & 59.45 & \bluelog{41.59} & \redlog{41.33} & 53.44 \\
        \bottomrule
        \end{tabular}
      }
      {The lowest three costs in each row are highlighted in bold and color, with the lowest one being underlined and colored red.}
\end{table}

We next consider a practical setting in which the supplier incurs an additional cost due to packaging requirements in actual production operations, aiming to balance daily inventory and packing costs.
{The packing cost follows an incremental discount scheme~\citep{tersine1991economic}, where the unit price decreases for quantities exceeding certain thresholds.
To align with the rescaled order quantities, we specify the following piecewise affine concave cost function:}
\begin{equation}\label{eq:exp:nonconvex cost}
C(z) = \min\{z, 0.6z+8, 0.4z+15.6 \}.
\end{equation}
Accordingly, we consider a Newsvendor problem with an additional term $C(z)$ in the objective function.
The inclusion of this concave term introduces nonconvexity into the optimization problem.
we use the deterministic MM algorithm to solve the ERM problem for learning RKHS-based decision rule, referred to as RKHS-DR(MM). 
For prescriptive methods, we solve one-dimensional nonconvex problems by line search. 
{As SOF from \cite{kallus2022stochastic} requires a new development of second-order perturbation analysis of stochastic optimization to accommodate the nonconvex term, we instead replace $C(z)$ by its affine piece $0.6z + 8$ when computing weights for training samples.
We refer to this linearized version as SOF(L).}
\begin{table}[htbp]
    \TABLE
    {Test costs for nonconvex Newsvendor on a real-world dataset\label{tab:exp:RWresults-NcvxNV}}
    {
        
    \small
    \begin{tabular}{cccccccccc}
        \toprule
        \multicolumn{1}{c}{Added features} & \multirow{2}[2]{*}{$p$} & \multicolumn{3}{c}{DR methods} & \multicolumn{4}{c}{Prescriptive methods} & \multicolumn{1}{l}{PO methods} \\
        \cmidrule{3-10}
        (all include basic features) &       & PADR  & RKHS-DR  & NN-DR & kNN   & CART  & RF    & SOF(L)   & PO-PA \\
        \midrule
         None & 11    & \redlog{74.28} & 81.89 & \bluelog{80.67} & 84.52 & \bluelog{81.56} & \bluelog{81.56} & 82.50 & 90.48 \\
         + order quantities in past 1 day & 13    & \redlog{64.75} & \bluelog{72.11} & \bluelog{68.19} & 77.86 & 86.58 & 73.78 & 74.45& 77.89 \\
         + order quantities in past 3 days & 15    & \redlog{64.10} & 68.29 & 66.62 & 75.58 & 78.31 & \bluelog{64.99} & \bluelog{66.54} & 73.61 \\
         + order quantities in past 5 days & 17    & \redlog{63.20} & \bluelog{64.72} & \bluelog{64.07} & 76.81 & 79.54 & 65.74 & 67.66 & 71.56 \\
         + order quantities in past 7 days & 19    & \redlog{62.01} & \bluelog{65.35} & \bluelog{62.25} & 75.38 & 85.95 & 68.35 & 69.23 & 69.97 \\
         + order quantities in past 14 days & 26    & \bluelog{64.86} & 65.15 & \redlog{63.45} & 73.81 & 75.52 & \bluelog{64.71} & 65.58 & 70.80 \\
         + order quantities in past 21 days & 33    & \bluelog{65.65} & \bluelog{65.15} & \redlog{63.33} & 73.81 & 75.52 & 65.73 & 66.01 & 71.02 \\
        \bottomrule
        \end{tabular}
    }
      {The lowest three  costs in each row are highlighted in bold and color, with the lowest one being underlined and colored red.}
\end{table}

{Table~\ref{tab:exp:RWresults-NcvxNV} demonstrates that PADR consistently outperforms other methods when $p \leq 17$, with performance gains exceeding 5\% for $p \leq 13$, and achieves the lowest test cost at $p = 19$. 
These findings highlight PADR's effectiveness in addressing nonconvex contextual SP problems. On the other hand, NN-DR maintains superior performance when $p \geq 26$.
Forest-based methods, however, perform relatively poorly in this nonconvex setting. 
The performance gap between SOF(L) and RF becomes more pronounced than in the convex case, since linearizing the concave cost introduces considerable model misspecification.
Moreover, the persistently weak performance of PO-PA across both convex and nonconvex scenarios indicates that its sequential prediction and decision-making process may not be suitable for this real-world dataset.}

\subsection{Product placement problem}\label{subsec:exp:pp}
We consider a product placement problem on a graph $\mathcal G=(\mathcal N, \mathcal A)$ with the aim to determine the quantities of products $z \in \RR^{\abs{\mathcal N}}$ to place at different nodes. 
After placing the products, demands $y \in \RR^{\abs{\mathcal N}}$ are observed at nodes $\mathcal N$, and then products are shipped across arcs $\mathcal A$ to fulfill these demands. 
We denote the node-arc matrix of $\mathcal G$ by $W$, the product placement cost by $c\in \RR^{\abs{\mathcal N}}$, the shipment cost along edges by $g \in \RR^{\abs{\mathcal A}}$, and the penalty for unmet demand by $b\in \RR^{\abs{\mathcal N}}$. The problem is formulated as follows,
\[
  \min_{z\geq 0}\quad \EE\left[
    c^\top z + V(z, Y)
    \,\middle\vert\, X = x
  \right],
\]
where $V(z, Y) = \min_{f \geq 0} \left\{g^\top f + b^\top (Y + Wf - z)^+ \right\}$. 

This is a two-stage stochastic linear program, and its ERM problem can be converted into a large linear program by adding the second-stage decision variables for each scenario. 
Following the graph generation settings in \citet{bertsimas2022data}, we set $\abs{\mathcal N}=5$, $\abs{\mathcal A} = 10$, $c = (3,\dots,3)$, $g = (1,\dots,1)$, and $b = (5, \dots, 5)$.
To evaluate PADR's performance under more complex conditions, we employed the following demand model:
\begin{equation}\label{eq:exp:sin+ma model}
\overline{Y}'(X) = \nonlinearity\alpha_1\sin(\alpha_2 X_1) + \nonlinearity\max\{\alpha_3 X_2, \alpha_4 X_2\} + \alpha_5,
\end{equation}
where the sine function represents the seasonal dependence of the demand on $X_1$, and the max-affine function represents the nonlinear trend with respect to $X_2$. We set $\nonlinearity=1$ unless otherwise specified.
We randomly generate parameters $\{\alpha_i\}_{i=1}^5$ and the node-arc matrix $W$ for each trial and report performance of methods over 10 trials for each experimental setting. 
{For NN-DR, incorporating the second-stage problem into the loss function within the standard training framework is nontrivial.
To address this difficulty, we treat the second-stage decision variables as additional model parameters and jointly optimize them along with the NN-DR parameters. 
It is worth noting that this implementation is tailored specifically to the problem at hand, and it underscores the limitation of NN–based approaches for multi-stage stochastic programs.}
The SOF method of \citet{kallus2022stochastic} is not included in this experiment, as applying it to a two-stage stochastic program would require substantial extensions involving second-order perturbation analysis, which are beyond the scope of this paper.

\subsubsection{Varying feature dimensions.}
We first study the effect of the feature dimension $p$ on the performance of PADR under sparse and dense demand models similar to the discussion in Section~\ref{subsubsec:NV:varyingfeature}. 
The results are reported in Figure~\ref{fig:exp:PP:p}. 
Similar to the results in Section~\ref{subsubsec:NV:varyingfeature} for the Newsvendor problem, with sparse demand model, PADR performs well in low-dimensional cases but deteriorates as $p$ increases which may necessitate the proper parameter regularization. 
NN-DR performs poorly in the sparse case, also suffering from overfitting.
In contrast, forest-based prescriptive methods are robust to redundant features. 
With dense demand model, PADR and NN-DR have the best performance, with the former one slightly outperforming when $p$ is small.
{We also report the running time results in Figure~\ref{fig:exp:PP:p-time}.
PADR requires substantially more computational time in this problem, as the underlying optimization becomes more complex and involves multiple decision variables. 
In the current implementation, the decision rule is constructed separately for each dimension.
Given so, we believe the computational efficiency of PADR could be further improved through vectorized multi-dimensional decision rule implementation and engineering optimizations in the ESMM surrogate construction.
}

Despite the increased complexity and the challenge in accurately approximating the true model, PADR continues to demonstrate superior performance for instances with dense demand model, high nonlinearity, and sparse demand models with small feature dimensions. %

\begin{figure}[htbp]
    \FIGURE
    {{
      \subfloat[Sparse demand model $\phantom{xxxxxx}$]{%
          \includegraphics[width=0.426\textwidth]{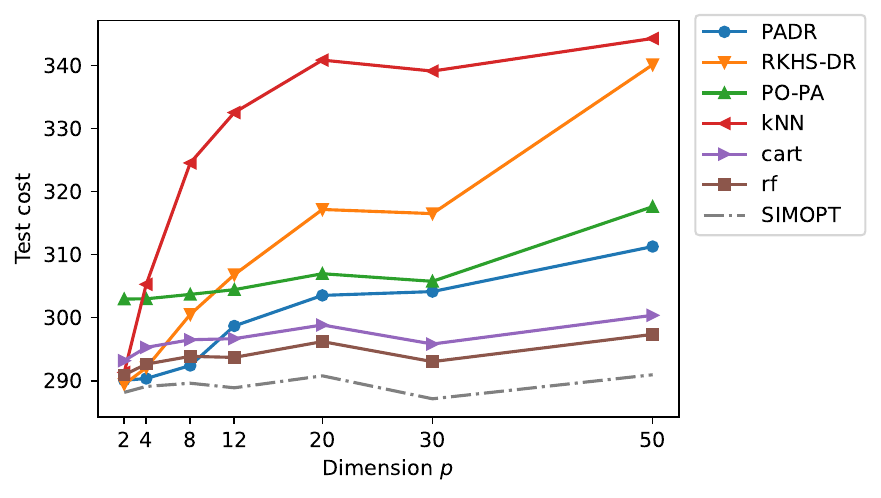}\label{subfig:pp-sinpa-k5-p-sparse}
      }
      \subfloat[Dense demand model $\phantom{xxxxxx}$]{%
            \includegraphics[width=0.426\textwidth]{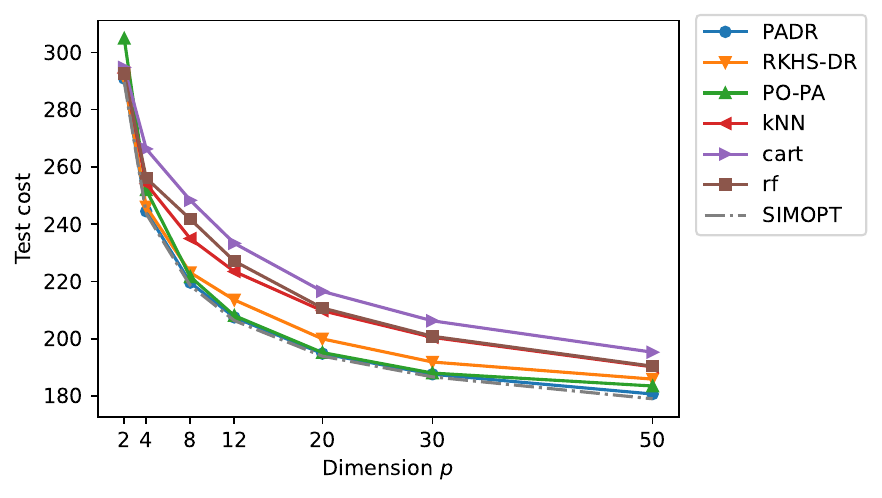}\label{subfig:pp-sinpa-k5-p-dense}
      }
    }}
    {Product placement problem with varying feature dimensions\label{fig:exp:PP:p}}
    {%
    }
\end{figure}

\begin{figure}[htbp]
    \FIGURE
    {{
      \subfloat[Sparse demand model $\phantom{xxxxxx}$]{%
          \includegraphics[width=0.426\textwidth]{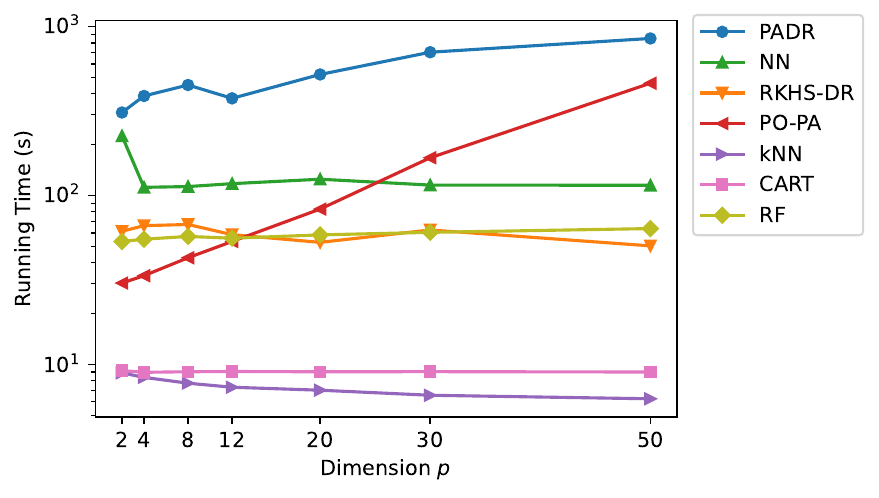}\label{subfig:pp-sinpa-k5-p-sparse-time}
      }
      \subfloat[Dense demand model $\phantom{xxxxxx}$]{%
            \includegraphics[width=0.426\textwidth]{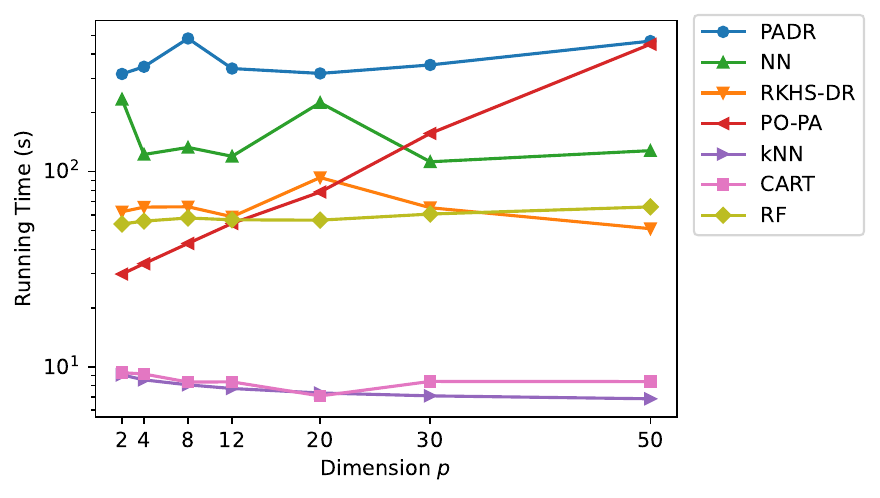}\label{subfig:pp-sinpa-k5-p-dense-time}
      }
    }}
    {Running time in product placement problem with varying feature dimensions\label{fig:exp:PP:p-time}}
    {%
    }
\end{figure}

\subsubsection{Nonlinearity in demand model.}
We also study the effect of nonlinearity in the demand model by varying the coefficient $\nonlinearity$. Notice that the nonlinear structure of the demand model here is different from the one in Section \ref{subsec:exp:Newsvendor} with the sine component function. 
The results for $p=20$ are provided in Figure~\ref{fig:exp:PP:k}.
The increasing gap in Figure~\ref{subfig:pp-k} indicates that except for PADR, NN-DR and PO-PA, the performance of all the other methods deteriorates as the nonlinearity of the demand model increases. Figure~\ref{subfig:pp-k-relative} indicates that PADR and NN-DR are slightly better than PO-PA and outperforms the other methods, and the relative gap among them remains approximately stable as $\nonlinearity$ increases.

\begin{figure}[htbp]
    \FIGURE
    {{
      \subfloat[$\phantom{xxxx}$]{%
          \includegraphics[width=0.426\textwidth]{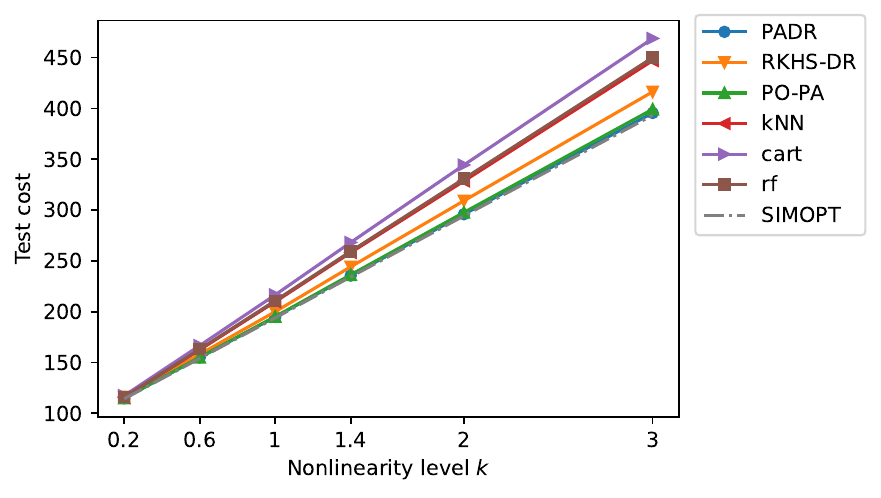}\label{subfig:pp-k}
      }
      \subfloat[$\phantom{xxxxx}$]{%
            \includegraphics[width=0.426\textwidth]{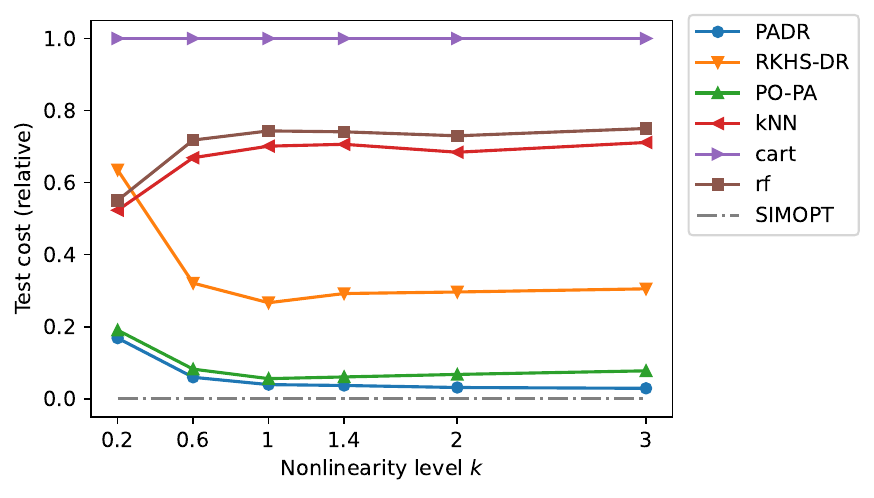}\label{subfig:pp-k-relative}
      }
    }}
    {Product placement problem: varying nonlinearity level $\nonlinearity$ (dense model, $p=20$) \label{fig:exp:PP:k}}
    {In (b), with the cost of the CART scaled to 1 and the cost of SIMOPT set to 0, we plot the relative cost of other methods compared to SIMOPT.
    }
\end{figure}

\section{Conclusion}\label{sec:Conclusion}
We propose the PADR-based ERM method together with the ESMM algorithm to solve a broad class of SP problems with covariate information.
For our PADR-based ERM framework, we provide the approximation bound of the PA function class and further establish the nonasymptotic consistency result for unconstrained problems.
To solve the ERM problem with a composite DC structure, we develop the ESMM algorithm by combining enhanced $\varepsilon$-active surrogation with sequential sampling. We provide the asymptotic convergence of the ESMM algorithm to composite strong d-stationarity, and the nonasymptotic convergence to d-stationarity via the error-bound analysis.  
Numerical experiments indicate the advantages of our method over PO methods and prominent DR-based methods for various instances regarding the test cost, computational time, and robustness to nonlinearity of the underlying dependency. 
In conclusion, the PADR is a favorable hypothesis class for DR-based ERM methods, supported by theoretical consistency, algorithmic convergence, and computational efficiency for solving SP with covariate information. 
Extensions of the PADR-based ERM framework for more complicated SP problems deserve to be explored in the future, such as problems with chance constraints or endogenous random variables, which necessitate a further investigation of model construction, consistency analysis, and data-driven algorithmic development for problems with more complex composite structures.

\bibliographystyle{informs2014}
\bibliography{refs}

\newpage

\setcounter{page}{1}\def\thepage{ec\arabic{page}}%
\par\noindent{\raggedright\fs.15.18.\bf 
\begin{center}
  $ $\\
  Supplementary Material for\\[8pt]
  \myfulltitle
\end{center}\endgraf}\vspace{8pt}
\begin{APPENDICES}
\section{Proofs of Theoretical Consistency Results in Section~\ref{sec:Model}}
\subsection{Proof of Proposition~\ref{prop:model:universal approximation of PADR}}
\proof{Proof of Proposition~\ref{prop:model:universal approximation of PADR}.}
    \label{proof:prop:universal approximation}
    \newcommand{\Coverrr}{\widehat{\cC}_{2}(\epsilon; \bar{\cX})}
    \newcommand{\Aone}[1]{\mathcal A_1(#1)}
    \newcommand{\Atwo}[1]{\mathcal A_2(#1)}
    It suffices to justify the bound when $d = 1$. 
    For $\epsilon > 0$ satisfying $(\sqrt{p}\bar{X} / \epsilon) \geq 1$, we first construct a $(2\epsilon/\sqrt{p})$-grid of the cube $\bar{\cX}\coloneqq [-\bar{X}, \bar{X}]^p$, which is also an $\epsilon\text{-cover}$ of $\bar{\cX}$ with $\ell_2$-norm, denoted by $\Coverrr$. Let
        $d_{\min}\coloneqq \underset{\widehat{x}_1, \widehat{x}_2 \in \Coverrr,\, \widehat{x}_1 \neq \widehat{x}_2}{\min}
        \norm{\widehat{x}_1 - \widehat{x}_2} = \frac{2\epsilon}{\sqrt{p}}$.
    The cardinality of $\Coverrr$ satisfies
    \[
        \left(\frac{\sqrt{p}\bar{X}}{\epsilon}\right)^p
        \leq \abs{\Coverrr} 
        \leq \left(\frac{\sqrt{p}\bar{X}}{\epsilon} + 1\right)^p
        \leq \left(\frac{2\sqrt{p}\bar{X}}{\epsilon}\right)^p.
    \]
    
    Given $f_0 \in \cF_{L_0, M_0}$, by \citet[Theorem 1]{mcshane1934extension}, we can extend $f_0$ to an $L_0$-Lipschitz function $\f: \Xcube \rightarrow \RR$ satisfying $\norm{\f}_\infty \leq M_0$ and $\f(x) = f_0(x)$ for any $x \in \cX$.
    Next, we construct a piecewise affine interpolation function to a set of points $\{(\widehat{x}, \bar{f}_0(\widehat{x})) \mid \widehat{x} \in \Coverrr\}$:
    \begin{align}\label{eq:EC:2-1:interpolation PA function}
        f_{\widehat{\cC}}(x) 
        \coloneqq &\max_{\widehat{x} \in \Coverrr} \left\{
            C\widehat{x}^\top(x - \widehat{x}) + \frac{C}{2} \norm{\widehat{x}}^2 + \frac{1}{2} \f(\widehat{x})
        \right\}
        - \max_{\widehat{x} \in \Coverrr} \left\{
            C\widehat{x}^\top(x - \widehat{x}) + \frac{C}{2} \norm{\widehat{x}}^2 - \frac{1}{2} \f(\widehat{x})
        \right\},
    \end{align}
    where $C \coloneqq L_0 / d_{\min} = L_0 \sqrt{p} / (2\epsilon)$.
    The function $f_{\widehat{\cC}}$ indeed interpolates $\{(\widehat{x},  \bar{f}_0(\widehat{x})) : \widehat{x} \in \Coverrr\}$.
    Specifically, for any $\widehat{y} \in \Coverrr$,
    \begin{align*}
        \max_{\widehat{x} \in \Coverrr} 
        \left\{
            C \widehat{x}^\top (\widehat{y} - \widehat{x}) + \frac{1}{2} C \norm{\widehat{x}}^2 + \frac{1}{2}\f(\widehat{x})\right\}
        &= \max_{\widehat{x} \in \Coverrr} 
        \left\{
            \frac{1}{2} C \norm{\widehat{y}}^2 - \frac{1}{2} C \norm{\widehat{x} - \widehat{y}}^2 + \frac{1}{2}\f(\widehat{x})
        \right\}\\
        &\leq \max_{\widehat{x} \in \Coverrr} 
        \left\{
            \frac{1}{2} C \norm{\widehat{y}}^2 - \frac{1}{2} L_0 \norm{\widehat{x} - \widehat{y}} + \frac{1}{2}\f(\widehat{x})
        \right\}\\
        &\leq \max_{\widehat{x} \in \Coverrr} 
        \left\{
            \frac{1}{2} C \norm{\widehat{y}}^2 - \frac{1}{2} \abs{\f(\widehat{x}) - \f(\widehat{y})} + \frac{1}{2}\f(\widehat{x})
        \right\}\\
        &\leq \frac{1}{2} C \norm{\widehat{y}}^2 + \frac{1}{2}\f(\widehat{y}),
    \end{align*}
    where the first inequality holds since $C = L_0 / d_{\min} \geq L_0 / \norm{\widehat{x} - \widehat{y}}$. %
    Moreover, this upper bound is attained by taking $\widehat{x} = \widehat{y}$ in the max operator.
    Similarly, the value of the second max-affine component of $f_{\widehat{\cC}}$ at $\widehat{y}$ is $\frac{1}{2} C \norm{\widehat{y}}^2 - \frac{1}{2}\f(\widehat{y})$. 
    Hence, we obtain $f_{\widehat{\cC}}(\widehat{y}) = \f(\widehat{y})$, verifying the interpolation.
    With $K \coloneqq \abs{\Coverrr}$, for any $\widehat{x} \in \Coverrr$, we obtain $
        \norm{C \widehat{x}}_\infty 
        \leq \frac{L_0 \sqrt{p}}{2\epsilon}\bar{X} 
        \leq \frac{L_0}{2}K^{1/p}$
    and
    \[\abs{-C \norm{\widehat{x}}^2 + \frac{1}{2} C \norm{\widehat{x}}^2 \pm \frac{1}{2}\f(\widehat{x})}
            \leq \frac{1}{2} C  p \norm{\widehat{x}}_{\infty}^2 %
            + \frac{1}{2}M_0
            \leq \frac{1}{4}p L_0 \bar{X} K^{1/p} + \frac{1}{2}M_0.
    \]
    Therefore, with $\mu \geq \max\left\{\frac{1}{2}L_0 K^{1/p}, \frac{1}{4}p L_0 \bar{X} K^{1/p} + \frac{1}{2}M_0\right\}$, we have $f_{\widehat{\cC}} \in \HK$.
    
    Note that for any $w \in \cX$, there exists $\widehat{x}(w)  \in \Coverrr$ such that $\norm{w - \widehat{x}(w) } \leq \epsilon$ and $f_{\widehat{\cC}}(\widehat{x}(w)) = \f(\widehat{x}(w))$.
    We next show that $f_{\widehat{\cC}}$ is Lipschitz continuous with modulus $(\sqrt{p} + 2)L_0$, with which 
    we obtain the desired inequality,
    \[
        \begin{array}{ll}
            \displaystyle{\min_{f \in \HK}} \max_{w \in \cX} \abs{f(w) - f_0(w)} 
            &\leq \displaystyle{\max_{w \in \cX}} \abs{f_{\widehat{\cC}}(w) - f_0(w)}  \leq \displaystyle{\max_{w \in \cX}} \left\{\abs{f_{\widehat{\cC}}(w) - f_{\widehat{\cC}}(\widehat{x}(w))} + \abs{\bar{f}_0(\widehat{x}(w)) - \bar{f}_0(w)}\right\}\\
            &\leq 2(p^{1/2} + 3)p^{1/2} L_0\bar{X} K^{-1/p},
        \end{array}
    \]
    where the last inequality follows from the Lipschitz continuity of $f_{\widehat{\cC}}$ and $\bar{f}_0$ and the relation $K \leq (2\sqrt{p}\bar{X} / \epsilon)^p$.
    For any $x \in \Xcube$, let
    \[
    \begin{array}{ll}
        & \mathcal A_1(x) \coloneqq \displaystyle{\mathop{\arg\max}_{\widehat{x} \in \Coverrr}} \left\{
            C \widehat{x}^\top (x - \widehat{x}) + \frac{C}{2} \norm{\widehat{x}}^2 + \frac{1}{2}\f(\widehat{x})
        \right\},   \\
        & \mathcal A_2(x) \coloneqq \displaystyle{\mathop{\arg\max}_{\widehat{x} \in \Coverrr}} \left\{
            C \widehat{x}^\top (x - \widehat{x}) + \frac{C}{2} \norm{\widehat{x}}^2 - \frac{1}{2}\f(\widehat{x})
        \right\}.
        \end{array}
    \]
    It suffices to show that $\norm{C(\widehat{x}_1^* - \widehat{x}_2^*)} \leq (\sqrt{p}+2)L_0$ for any $\widehat{x}_1^* \in \Aone{x}$ and $\widehat{x}_2^* \in \Atwo{x}$. 
    Let $\widehat{x} \in \Coverrr$ satisfy $\norm{x - \widehat{x}} \leq \epsilon$.
    For any $\widehat{y} \in \Coverrr$ satisfying $\norm{x - \widehat{y}} > \frac{\sqrt{p} + 2}{\sqrt{p}}\epsilon$, we have that
    \begin{align*}
        C \widehat{y}^\top (x - \widehat{y}) + \frac{1}{2} C \norm{\widehat{y}}^2 + \frac{1}{2}\f(\widehat{y}) & = \frac{1}{2} \left[C \norm{x}^2 - C \norm{x - \widehat{y}}^2 + \f(\widehat{y})\right]\\
        &= \frac{1}{2} \left[\f(\widehat{y}) - \f(x) - C\norm{x - \widehat{y}}^2\right] + \frac{1}{2}\left[C \norm{x}^2 + \f(x)\right]\\
        &\leq \frac{1}{2} \left[
            L_0 \norm{x - \widehat{y}} - \frac{\sqrt{p}L_0}{2 \epsilon} \norm{x - \widehat{y}}^2
        \right] + \frac{1}{2} \left[C \norm{x}^2 + \f(x)\right]\\
        &< \frac{1}{2} \left[
            L_0 \cdot \frac{\sqrt{p}+2}{\sqrt{p}}\epsilon - \frac{\sqrt{p}L_0}{2\epsilon}\left(\frac{\sqrt{p}+2}{\sqrt{p}}\epsilon\right)^2
        \right] + \frac{1}{2}\left[C \norm{x}^2 + \f(x)\right]\\
        &= \frac{1}{2}\left[-L_0 \cdot \epsilon - \frac{\sqrt{p}L_0}{2\epsilon}\epsilon^2\right] + \frac{1}{2}\left[C \norm{x}^2 + \f(x)\right]\\
        &\leq \frac{1}{2}\left[\f(\widehat{x}) - \f(x) - C\norm{x - \widehat{x}}^2\right] + \frac{1}{2}\left[C \norm{x}^2 + \f(x)\right]\\
        &= C\widehat{x}^\top (x - \widehat{x}) + \frac{1}{2} C\norm{\widehat{x}}^2 + \frac{1}{2}\f(\widehat{x}),
    \end{align*}
    where the second inequality is derived by the strict monotonicity of the quadratic function $t \mapsto L_0 t - \frac{\sqrt{p}L_0}{2\epsilon} t^2$ on $[\epsilon/\sqrt{p}, +\infty)$. 
    The above result shows $\norm{\widehat{x}_1^* - x} \leq \frac{\sqrt{p} + 2}{\sqrt{p}}\epsilon$ for all $\widehat{x}_1^* \in \Aone{x}$. Similarly, we have $\norm{\widehat{x}_2^* - x} \leq \frac{\sqrt{p}+2}{\sqrt{p}}\epsilon$ for all $\widehat{x}_2^* \in \Atwo{x}$. It follows that
    \[
        \norm{C(\widehat{x}_1^* - \widehat{x}_2^*)}
        \leq C\norm{\widehat{x}_1^* - x} + C \norm{\widehat{x}_2^* - x}
        \leq 2\frac{\sqrt{p}L_0}{2\epsilon}\cdot \frac{\sqrt{p}+2}{\sqrt{p}} \epsilon
        \leq (\sqrt{p} + 2)L_0,
    \]
    which completes the proof.\qedsymbol
\endproof

\subsection{Proofs of Proposition~\ref{prop:model:Uniform generalization error bound} and Theorem~\ref{thm:excess-wo-constr}}\label{EC:subsec:proofs for model:nonasymptotic consistency}

\subsubsection{Preliminaries on Rademacher complexity.}
We start with a brief introduction to the Rademacher complexity, which is utilized to quantify the generalization error.
The empirical Rademacher complexity of a function class is defined as follows.
\begin{definition}
    Let $\Phi$ be a class of functions mapping from $\cU$ to $\RR$ and $\Xi_N = \{u^s\}_{s=1}^N \subseteq \cU$ be a sample set of size $N$. Then the empirical Rademacher complexity of $\Phi$ with respect to $\Xi_N$ is defined as
    \[
        \widehat{\cR}_{\Xi_N}(\Phi) \coloneqq \EE_{\boldsymbol{\sigma}} 
        \left[
        \sup_{\phi \in \Phi} \frac{1}{N}\sum_{s = 1}^N \sigma_s \phi(u^s)
        \right],
    \]
    where ${\boldsymbol\sigma} = (\sigma_1,\dots,\sigma_N)^\top$ is a uniform random vector taking values in $\{-1, +1\}^N$.
\end{definition}
The Rademacher complexity of $\Phi$ is defined as the expectation of the empirical Rademacher complexity over all samples of size $N$ drawn from the underlying distribution:
\[
    \cR_N(\Phi) \coloneqq \EE_{\Xi_N}[\widehat{\cR}_{\Xi_N}(\Phi)].
\]
It captures the richness of the class $\Phi$ by measuring the degree to which a hypothesis set can fit random noise.
The following lemma quantifies the Rademacher complexity of a class of parametric functions, which is generalized from \citet[Lemma B.2]{ermoliev2013sample} with minor modifications made to accommodate the developments in this paper.
\begin{lemma}
    \label{lemma:EC:Rademacher of Param functions}
    Let $\widehat{\Phi} \coloneqq \{\phi(\cdot; \theta):\cU \rightarrow \RR \mid \theta \in \Theta\}$ be a class of functions parameterized by vectors in a compact set $\Theta \subseteq \RR^q$ with $\max_{\theta\in \Theta}\norm{\theta}_\infty = \mu$. Suppose that there exist $M_{\widehat{\Phi}} > 0$ and $L_{\widehat{\Phi}} > 0$ such that $\max_{\theta \in \Theta}\norm{\phi(\cdot; \theta)}_\infty \leq M_{\widehat{\Phi}}$ and $\norm{\phi(\cdot; \theta_1) - \phi(\cdot; \theta_2)}_\infty \leq L_{\widehat{\Phi}} \norm{\theta_1 - \theta_2}$ for any $\theta_1, \theta_2\in \Theta$. Then
    \[
        \cR_N(\widehat\Phi) \leq \frac{2\mu L_{\widehat{\Phi}} \sqrt{q}}{N} + M_{\widehat{\Phi}}\sqrt{\frac{q\log N}{N}}.
    \]
\end{lemma}
\proof{Proof of Lemma~\ref{lemma:EC:Rademacher of Param functions}}
    Given a sample set $\Xi_N = \{u^s\}_{s=1}^N$ of size $N$, following the proof of \citet[Lemma B.2]{ermoliev2013sample}, we derive that for any $\epsilon \in (0, 2\mu]$,
    \[
        \widehat{\cR}_{\Xi_N}(\widehat{\Phi}) \leq L_{\widehat{\Phi}}\sqrt{q}\epsilon + M_{\widehat{\Phi}}\sqrt{\frac{2q\log(2\mu / \epsilon)}{N}}.
    \]
    By setting $\epsilon = 2\mu / N$ and taking expectation over $\Xi_N$, we obtain the desired upper bound.\qedsymbol
\endproof

The following lemma gives the generalization bound based on the Rademacher complexity.

\begin{lemma}
    [\citet{ECmohri2018foundations}, Theorem 3.3]\label{lemma:EC:Gen Error bound with Rademacher complex}
    Let $\Phi$ be a function class mapping from $\cU$ to $[-M, M]$.
    Then for any $\delta > 0$, with probability at least $1-\delta$ over the draw of an independently and identically distributed (i.i.d.) sample $\Xi_{N}$ of size $N$, the following holds true for all $\phi \in \Phi$:
    \[
        \EE[\phi(\tilde{u})]
        \leq 
        \frac{1}{N}\sum_{s=1}^N \phi(u^s) + 2{\cR}_{N}(\Phi) + \sqrt{2} M \sqrt{\frac{\log(1/\delta)}{N}}.
    \]
\end{lemma}

\subsubsection{Proofs of Proposition~\ref{prop:model:Uniform generalization error bound} and Theorem~\ref{thm:excess-wo-constr}.}
\proof{Proof of Proposition~\ref{prop:model:Uniform generalization error bound}.}
    \label{proof:prop:uniform gen error bound}
    \newcommand{\SU}{S_u}
    \newcommand{\SV}{S_v}
    \newcommand{\EEg}{\EE[\phi]}
    \newcommand{\ES}[1]{\widehat{\EE}_{#1}[\phi]}
    We first quantify the Rademacher complexity of the composite function class $\Phi_{\HK} \coloneqq \left\{\phi:\cX\times \cY \rightarrow \RR \mid \phi(x, y) = \cost{f(x)}{x}{y},\,f\in \HK\right\}$. 
    Notice that any function in $\Phi_{\HK}$ is parameterized by some $\theta \in [-\mu, \mu]^{2dK(p+1)}$. For any $(x, y)\in \cX \times \cY$ and $\theta_1, \theta_2 \in [-\mu, \mu]^{2dK(p+1)}$,
    \[
        \abs{\varphi(f(x, \theta_1), y) - \varphi(f(x, \theta_2), y)}
        \leq L_\varphi \norm{f(x, \theta_1) - f(x, \theta_2)} \leq 2 L_\varphi \sqrt{p \bar{X}^2 + 1} \norm{\theta_1 - \theta_2},
    \]
    where the second inequality follows from Assumption~\ref{assum:A2:compactness of X}.
    Moreover, for any $f \in \HK$,
    \begin{equation}\label{eq:EC:boundedness of varphi}
        \abs{\varphi(f(x), y)} 
        \leq \abs{\varphi(z_0, y)} + \abs{\varphi(f(x), y) - \varphi(z_0, y)}
        \leq C_{\varphi} + L_\varphi \norm{f(x)-z_0} \leq M_{p,d}(\mu).
    \end{equation}
    Therefore, by Lemma~\ref{lemma:EC:Rademacher of Param functions}, we deduce
    \begin{equation}\label{eq:EC:Rad. complexity of PhiHK}
        \cR_n(\Phi_{\HK}) \leq \frac{4\mu L_\varphi\sqrt{2dK(p+1)(p\bar{X}^2 + 1)}}{n} + M_{p,d}(\mu)\sqrt{\frac{2dK(p+1)\log n}{n}}.
    \end{equation}
    
    Now we prove the uniform generalization bound. By Lemma~\ref{lemma:EC:Gen Error bound with Rademacher complex} and the inequality~\eqref{eq:EC:boundedness of varphi}, with probability at least $1-\delta/2$,
    \[
        \max_{\phi \in \Phi_{\HK}}\left(\EE_{X,Y}\phi(X, Y) - \frac{1}{n}\sum_{s=1}^n \phi(x^s, y^s)\right) 
        \leq 2 {\cR}_{n}(\Phi_{\HK}) + \sqrt{2} M_{p,d}(\mu) \sqrt{\frac{\log(2/\delta)}{n}}.
    \]
    For $\Phi_{\HK}^- \coloneqq \{-\phi:\phi \in \Phi_{\HK}\}$, we have ${\cR}_{n}(\Phi_{\HK}^-) = {\cR}_{n}(\Phi_{\HK})$.
    Therefore, by applying Lemma~\ref{lemma:EC:Gen Error bound with Rademacher complex} again to $\Phi_{\HK}^-$ and using the union bound, we derive that with probability at least $1 - \delta$,
    \begin{equation}\label{eq:proof:gen bound 2nd}
        \max_{f \in \HK}\abs{R(f) - \widehat{R}_{\Xi}(f)} \leq 2 \cR_n(\Phi_{\HK}) + \sqrt{2}M_{p,d}(\mu)\sqrt{\frac{\log(2/\delta)}{n}}.
    \end{equation}
    Finally, by combining~\eqref{eq:EC:Rad. complexity of PhiHK} and~\eqref{eq:proof:gen bound 2nd}, we complete the proof.\qedsymbol
\endproof

\proof{Proof of Theorem~\ref{thm:excess-wo-constr}.}
    Let $f^*_{\HK}$ denote the optimal solution to problem~\eqref{eq:intro:contextalSP:Full-E} restricted to $\HK$.
    Following the excess risk decomposition~\eqref{eq:model-3-excess-risk:excess-risk-decomposition-1}, we derive the excess risk bound by bounding the approximation and estimation errors respectively.
    With $f_{app} \in \mathop{\arg\min}_{f \in \HK} \norm{f - f^*}_\infty$, the approximation error can be bounded as follows:
    \begin{align*}
            R(f^*_{\HK}) - \min_{f\in \cF} R(f)
            &\leq R(f_{app}) - R(f^*)\\
            &= \EE_{X,Y} \left[\cost{f_{app}(X)}{X}{Y} - \cost{f^*(X)}{X}{Y}\right]\\
            &\leq L_\varphi \sup_{x\in \cX}\norm{f_{app}(x) - f^*(x)}_2\\
            &\leq 2L_\varphi \sqrt{d} (p^{1/2} + 3)p^{1/2}L_z \bar{X}K^{-1/p},
    \end{align*}
    where the last inequality follows from Proposition~\ref{prop:model:universal approximation of PADR}. 
    For the estimation error, we have
    \begin{align*}
    R(f^*_{\Xi, \HK}) - R(f^*_{\HK}) 
    &= R(f^*_{\Xi, \HK}) - \widehat{R}_{\Xi}(f^*_{\Xi, \HK}) + \widehat{R}_{\Xi}(f^*_{\Xi, \HK})  - R(f^*_{\HK})  \notag\\
    & \leq R(f^*_{\Xi, \HK}) - \widehat{R}_{\Xi}(f^*_{\Xi, \HK}) + \widehat{R}_{\Xi}(f^*_{\HK})  - R(f^*_{\HK})  \notag\\
    & \leq 2 \left( \max_{ f\in \HK} \, | R(f) - \widehat{R}_{\Xi}(f) |  \right).
    \end{align*}
    By Proposition~\ref{prop:model:Uniform generalization error bound}, with probability at least $1 - \delta$, we have
    \begin{align*}
    R(f^*_{\Xi, \HK}) - R(f^*_{\HK})
    \leq &\ \frac{16\mu L_\varphi\sqrt{2dK(p+1)(p\bar{X}^2 + 1)}}{n} + 4M_{p,d}(\mu)\sqrt{\frac{2dK(p+1)\log n}{n}}\\
    &\ + 2\sqrt{2}M_{p,d}(\mu) \sqrt{\frac{\log(2 / \delta)}{n}}.
    \end{align*}
    By combining these two error bounds, we obtain the desired inequality. \qedsymbol
\endproof

\section{Surrogation Construction of $F(\cdot, \xi)$ for~\ref{enum:algo:phi1} to~\ref{enum:algo:phi3}}
\label{EC:sec:surrogation for sec3}
We provide the construction of surrogation of $F(\cdot, \xi)\coloneqq \varphi(f(\cdot,\xi), \xi)$ for~\ref{enum:algo:phi1} to~\ref{enum:algo:phi3} in the proof of Proposition~\ref{prop:verify Cases1-3}.
In order to verify Assumption~\ref{assum:B4:minSurr dd-consistency}, 
we prove that the minimum surrogate function $\widehat{F}^{\varepsilon}_{\min}(\cdot, \xi; \bar{\theta})$ at a reference point $\bar{\theta}$ satisfies the following condition for any $\xi \in \Xi$:
\begin{enumerate}[label = (P\arabic*), resume*=P]
    \item\label{enum:algo:P5:quadratic surrogation gap} 
       Quadratic surrogation gap: for some $\kappa > 0$, $F(\theta, \xi) + \displaystyle{\frac{\kappa}{2}}\norm{\theta - \bar{\theta}}^2 \geq \widehat{F}^\varepsilon_{\min}(\theta, \xi; \bar{\theta})$ for any $\theta \in \Theta$.
\end{enumerate}
For directionally differentiable surrogate functions, Condition~\ref{enum:algo:P5:quadratic surrogation gap} is stronger than directional derivative consistency condition because it also limits the growth rate of the surrogate function with respect to the original function, which is essential to develop the nonasymptotic convergence analysis in Section~\ref{EC:sec:nonasymptotic convergence for general cases}.

\proof{Proof of Proposition~\ref{prop:verify Cases1-3}}

We first construct the surrogate functions for the three cases and verify Assumptions~\ref{assum:B1:Touching}-\ref{assum:B3:Continuous given fixed I}, followed by the verification of Assumption~\ref{assum:B4:minSurr dd-consistency}.

\underline{{\it Case 1:}} 
Suppose that $\varphi(z, \xi) = \varphi_1(z, \xi) - \varphi_2(z, \xi)$, where $\varphi_1(z, \xi) = \max_{i \in [K_\varphi]}\varphi_{1i}(z, \xi)$ and $\varphi_2(z, \xi) = \max_{i \in [K_\varphi]}\varphi_{2i}(z, \xi)$ are both max-affine functions. We denote the affine functions by $\{\varphi_{1i}(z, \xi) = a_{1i}^\top z + b_{1i}\}_{i \in [K_\varphi]}$ and $\{\varphi_{2i}(z, \xi) = a_{2i}^\top z + b_{2i}\}_{i \in [K_\varphi]}$, where $\xi$ represents the random components of $\{a_{1i}, b_{1i}, a_{2i}, b_{2i}\}_{i \in [K_\varphi]}$. Then the surrogate function of $F(\cdot, \xi)$ at $\bar{\theta}$ is constructed as follows,
\begin{align*}
    \widehat{F}(\theta, \xi; \bar{\theta}, I) 
    &= \max_{i \in [K_\varphi]} \left\{
    \max\{a_{1i}, 0\}^\top \widehat{f}(\theta, \xi; I_2) - \max\{-a_{1i}, 0\}^\top \widecheck{f}(\theta, \xi; I_1) + b_{1i}\right\}\\
    &\phantom{=} - \left(\max\{a_{2I_3}, 0\}^\top \widecheck{f}(\theta, \xi; I_1) - \max\{-a_{2I_3}, 0\}^\top \widehat{f}(\theta, \xi; I_2) + b_{2I_3}\right),
\end{align*}
where $I_1 \in \mathcal I_h^\varepsilon(\bar{\theta}, \xi)$, $I_2 \in \mathcal I_g^\varepsilon(\bar{\theta}, \xi)$, and $I_3 \in \mathcal I_{\varphi_2}^\varepsilon(f(\bar{\theta}, \xi), \xi)$. With a slight abuse of notation, we denote $I = (I_1, I_2, I_3) \in \mathcal I^\varepsilon(\bar{\theta}, \xi) \coloneqq \mathcal I_g^\varepsilon(\bar{\theta}, \xi) \times \mathcal I_h^\varepsilon(\bar{\theta}, \xi) \times \mathcal I_{\varphi_2}^\varepsilon(f(\bar{\theta}, \xi), \xi)$. Assumptions~\ref{assum:B1:Touching} and~\ref{assum:B2:Majorizing and Convex} hold in this case by construction.
Assumption~\ref{assum:B3:Continuous given fixed I} is justified since given a fixed surrogation index $I$, the surrogate function is independent of the reference point $\bar{\theta}$.

\underline{{\it Case 2:}} by assumption, $\varphi(\cdot, \xi)$ can be written as
$\varphi(\cdot, \xi) = \varphi^{\uparrow}(\cdot, \xi) + \varphi^{\downarrow}(\cdot, \xi)$,
where $\varphi^{\uparrow}(\cdot, \xi)$ and $\varphi^{\downarrow}(\cdot, \xi)$ are nondecreasing and nonincreasing convex functions respectively, that is, for any $\theta_1, \theta_2 \in \Theta$ satisfying $\theta_1 \leq \theta_2$ componentwisely,
we have $\varphi^{\uparrow}(\theta_2, \xi) \geq \varphi^{\uparrow}(\theta_1, \xi)$ and $\varphi^{\downarrow}(\theta_2, \xi) \leq \varphi^{\downarrow}(\theta_1, \xi)$.
With $I = (I_{1},I_{2}) \in \cI^\varepsilon(\bar{\theta}, \xi)$, the  surrogate function of $F(\cdot, \xi)$ at $\wtheta$
is constructed as
\[
    \widehat{F}(\theta, \xi; \wtheta, I) 
    = \varphi^{\uparrow}(\widehat{f}(\theta, \xi; I_{2}), \xi) 
    + \varphi^{\downarrow}(\widecheck{f}(\theta, \xi; I_{1}), \xi),
\]
where $\widehat{f}(\theta, \xi; I_{2})$ and $\widecheck{f}(\theta, \xi; I_{1})$ are the upper and lower  surrogate functions of $f(\theta, \xi)$ respectively.
By monotonicity of $\varphi^{\uparrow}$ and $\varphi^{\downarrow}$, $
        \widehat{F}(\theta, \xi; \wtheta, I)\geq \varphi^{\uparrow}(f(\theta, \xi), \xi) + \varphi^{\downarrow}(f(\theta, \xi), \xi) = F(\theta, \xi)$ for any $\theta \in \Theta$,
and thus the majorizing condition~\ref{enum:algo:P2:majorization} holds. 
Moreover, $\widehat{F}(\cdot, \xi; \bar{\theta}, I)$ satisfies the convex condition~\ref{enum:algo:P3:convexity} since a nondecreasing convex function composite with a convex function and a nonincreasing convex function composite with a concave function are both convex. 
The touching condition~\ref{enum:algo:P1:touching} clearly holds when $\varepsilon=0$. 
With a fixed index~$I$, $\widehat{F}(\cdot, \xi; \bar{\theta}, I)$ is also independent of $\bar{\theta}$ in this case, thus validating the bi-variate continuity assumption~\ref{assum:B3:Continuous given fixed I}.

\underline{{\it Case 3:}} since $\Xi$ is finite, for any $\xi \in \Xi$, $\varphi(\cdot, \xi)$ is Lipschitz smooth with uniform modulus denoted by $L_{\nabla \varphi}$, and $f(\cdot, \xi)$ is Lipschitz continuous with uniform modulus $L_f \coloneqq 2\sqrt{\max_{s\in[n]}\norm{x^s}^2 + 1}$.
Then the  surrogate function of $F(\cdot, \xi)$ at $\wtheta$ with $I = (I_{1}, I_{2}) \in \cI^\varepsilon(\wtheta, \xi)$ is constructed by surrogating the inner PA function according to componentwise signs of the gradient $\nabla\varphi(f(\bar{\theta}, \xi), \xi)$ as follows,
\begin{align*}
        \widehat{F}(\theta, \xi; \wtheta, I)
        &= \max\{\nabla{\varphi}(f(\wtheta, \xi), \xi), 0\}^\top \widehat{f}(\theta, \xi; I_{2})
        - \max\{-\nabla{\varphi}(f(\wtheta, \xi), \xi), 0\}^\top \widecheck{f}(\theta, \xi; I_{1})\\
        &\phantom{=,} - \nabla{\varphi}(f(\wtheta, \xi), \xi)^\top f(\wtheta, \xi) + \varphi(f(\wtheta, \xi), \xi) + \frac{1}{2}L_{\nabla\varphi} L_f^2 \norm{\theta - \wtheta}^2.
\end{align*}
The convex condition~\ref{enum:algo:P3:convexity} holds evidently; the touching condition~\ref{enum:algo:P1:touching} also holds when $\varepsilon=0$.
The surrogate $\widehat{F}(\cdot, \xi; \cdot, I)$ is continuous on $\Theta\times\Theta$ since $\nabla\varphi(\cdot, \xi)$ is continuous, validating Assumption~\ref{assum:B3:Continuous given fixed I}.
Based on the quadratic upper bound of a Lipschitz smooth function, the majorizing condition~\ref{enum:algo:P2:majorization} is verified as follows. 
Since $\varphi(\cdot, \xi)$ is $L_{\nabla\varphi}$-smooth, we have for any $z_1,z_2 \in \RR^d$ that 
\[
    \varphi(z_1, \xi) \leq \varphi(z_2, \xi) + \nabla\varphi(z_2, \xi)^\top(z_1 - z_2) + \frac{L_{\nabla\varphi}}{2}\norm{z_1 - z_2}^2.
\]
By setting $z_1 = f(\theta, \xi)$ and $z_2 = f(\wtheta, \xi)$, we obtain that for any $I = (I_{1}, I_{2}) \in \cI^\varepsilon(\wtheta, \xi)$,
\begin{align*}
    \label{eq:algo:surrogate for the lsmooth case}
    F(\theta, \xi) 
    &\leq F(\wtheta, \xi) + \nabla\varphi(f(\wtheta, \xi), \xi)^\top(f(\theta, \xi) - f(\wtheta, \xi)) 
        + \frac{L_{\nabla\varphi}}{2}\norm{f(\theta, \xi) - f(\wtheta, \xi)}^2\\
    &\leq F(\wtheta, \xi) 
          + \max\{\nabla{\varphi}(f(\wtheta, \xi), \xi), 0\}^\top \widehat{f}(\theta, \xi; I_{2})\\
    &\phantom{\leq,} - \max\{-\nabla{\varphi}(f(\wtheta, \xi), \xi), 0\}^\top \widecheck{f}(\theta, \xi; I_{1}) - \nabla\varphi(f(\wtheta, \xi), \xi)^\top f(\wtheta, \xi) + \frac{L_{\nabla\varphi}L_f^2}{2}\norm{\theta - \wtheta}^2\\
    &= \widehat{F}(\theta, \xi; \wtheta, I).
\end{align*}

Now we prove Assumption~\ref{assum:B4:minSurr dd-consistency} for~\ref{enum:algo:phi1} to~\ref{enum:algo:phi3}.
We prove that the minimum surrogate function $\widehat{F}^{\varepsilon}_{\min}(\cdot, \xi; \bar{\theta})$ satisfies Condition~\ref{enum:algo:P5:quadratic surrogation gap}.
Assumption~\ref{assum:B4:minSurr dd-consistency} is then verified since $\widehat{F}^{\varepsilon}_{\min}(\cdot, \xi; \bar{\theta})$ by construction is directionally differentiable and satisfies~\ref{enum:algo:P2:majorization}.
Since $\HatFmin(\cdot, \xi; \bar{\theta})$ also satisfies the majorizing condition~\ref{enum:algo:P2:majorization}, by the compactness of $\Theta$, it suffices to show that 
\[
    \widehat{\kappa} \coloneqq \limsup_{\theta \rightarrow \bar{\theta}, \theta \neq \bar{\theta}} \frac{
        2(\HatFmin(\theta, \xi; \bar{\theta}) - F(\theta, \xi))
    }{\norm{\theta - \bar{\theta}}^2} < +\infty.
\]
For the piecewise affine function $f(\theta, \xi) \coloneqq g(\theta, \xi) - h(\theta, \xi)$, we have $\cI_h^0(\theta, \xi) \subseteq \cI_h^\varepsilon(\bar{\theta}, \xi)$ when $\theta \in \Theta$ is sufficiently close to $\bar{\theta}$.
It follows from~\ref{enum:algo:P2:majorization} that for any $\theta\in \Theta$ sufficiently close to $\bar{\theta}$,
\[
    f(\theta, \xi) 
    \leq \min_{I_2 \in \cI^\varepsilon_h(\bar{\theta}, \xi)} \widehat{f}(\theta, \xi; I_2)
    \leq \min_{I_2 \in \cI_h^0(\theta, \xi)}\widehat{f}(\theta, \xi; I_2) = f(\theta, \xi).
\]
Therefore, by similar arguments for lower  surrogate functions of $f(\cdot, \xi)$, we obtain that for any $\theta\in \Theta$ sufficiently close to $\bar{\theta}$,
\[
    \min_{I_2 \in \cI^\varepsilon_h(\bar{\theta}, \xi)}\widehat{f}(\theta, \xi; I_2) 
    = f(\theta, \xi) 
    = \max_{I_1 \in \cI^\varepsilon_g(\bar{\theta}, \xi)}\widecheck{f}(\theta, \xi; I_1).
\]

For~\ref{enum:algo:phi1}, when $\theta$ is sufficiently close to $\bar{\theta}$, we have
\begin{align}\label{eq:=F}
    \HatFmin(\theta, \xi; \bar{\theta})
    &\coloneqq \min_{(I_1, I_2, I_3) \in \mathcal I^\varepsilon(\bar{\theta}, \xi)}\widehat{F}(\theta, \xi; \bar{\theta}, I) \notag\\
    &= \underbrace{\min_{I_1 \in \cI^\varepsilon_g(\bar{\theta}, \xi), I_2 \in \cI^\varepsilon_h(\bar{\theta}, \xi)}\max_{i \in [K_\varphi]} \left\{
    \max\{a_{1i}, 0\}^\top \widehat{f}(\theta, \xi; I_2) - \max\{-a_{1i}, 0\}^\top \widecheck{f}(\theta, \xi; I_1) + b_{1i}\right\}}_{\mathcal V_1}\notag\\
    &\phantom{==} - \underbrace{\max_{(I_1, I_2, I_3) \in \mathcal I^\varepsilon(\bar{\theta}, \xi)}\left(\max\{a_{2I_3}, 0\}^\top \widecheck{f}(\theta, \xi; I_1) - \max\{-a_{2I_3}, 0\}^\top \widehat{f}(\theta, \xi; I_2) + b_{2I_3}\right)}_{\mathcal V_2}\notag\\
    &= F(\theta, \xi).
\end{align}
The third equality is validated as follows. For the first component $\mathcal V_1$, we have
\[
\begin{aligned}
    \mathcal V_1 
    &\geq \max_{i \in [K_\varphi]}\min_{I_1 \in \cI^\varepsilon_g(\bar{\theta}, \xi), I_2 \in \cI^\varepsilon_h(\bar{\theta}, \xi)} \left\{
    \max\{a_{1i}, 0\}^\top \widehat{f}(\theta, \xi; I_2) - \max\{-a_{1i}, 0\}^\top \widecheck{f}(\theta, \xi; I_1) + b_{1i}\right\}\\
    &= \max_{i \in [K_\varphi]} \left\{\max\{a_{1i}, 0\}^\top \min_{I_2 \in \cI^\varepsilon_h(\bar{\theta}, \xi)}\widehat{f}(\theta, \xi; I_2) - \max\{-a_{1i}, 0\}^\top \max_{I_1 \in \cI^\varepsilon_g(\bar{\theta}, \xi)}\widecheck{f}(\theta, \xi; I_1) + b_{1i}\right\}\\
    &= \max_{i \in [K_\varphi]}\left\{a_{1i}^\top f(\theta, \xi) + b_{1i}\right\} \eqqcolon \varphi_{1}(f(\theta, \xi), \xi).
\end{aligned}
\]
Since there exists $\bar{I}_1 \in \mathcal I^\varepsilon_g(\bar{\theta}, \xi)$ and $\bar{I}_2 \in \mathcal I^\varepsilon_h(\bar{\theta}, \xi)$ such that $\widehat{f}(\theta, \xi; \bar{I}_2) = f(\theta, \xi)$ and $\widecheck{f}(\theta, \xi; \bar{I}_1) = f(\theta, \xi)$, we also have
\[
\begin{aligned}
    \mathcal V_1 
    &\leq \max_{i \in [K_\varphi]} \left\{\max\{a_{1i}, 0\}^\top \widehat{f}(\theta, \xi; \bar{I}_2) - \max\{-a_{1i}, 0\}^\top \widecheck{f}(\theta, \xi; \bar{I}_1) + b_{1i}\right\}\\
    &= \max_{i \in [K_\varphi]}\left\{a_{1i}^\top f(\theta, \xi) + b_{1i}\right\} \eqqcolon \varphi_{1}(f(\theta, \xi), \xi),
\end{aligned}
\]
indicating $\mathcal V_1 = \varphi_1(f(\theta, \xi), \xi)$. Moreover, for the second component $\mathcal V_2$, when $\theta$ is sufficiently close to $\bar{\theta}$, there exists $\bar{I}_3 \in \cI_{\varphi_2}^\varepsilon(f(\bar{\theta}, \xi), \xi)$ such that $\varphi_{2\bar{I}_3}(f(\theta, \xi), \xi) = \varphi_2(f(\theta, \xi), \xi)$.
Then, we have 
\begin{align*}
\mathcal V_2 
&\geq \max\{a_{2\bar{I}_3}, 0\}^\top \widecheck{f}(\theta, \xi; \bar{I}_1) 
- \max\{-a_{2\bar{I}_3}, 0\}^\top \widehat{f}(\theta, \xi; \bar{I}_2) + b_{2\bar{I}_3}
= \varphi_2(f(\theta, \xi), \xi)
\end{align*}
and
\begin{align*}
\mathcal V_2
&\leq \max_{
(I_1, I_2) \in \mathcal I_g^\varepsilon(\bar{\theta}, \xi) \times \cI_h^\varepsilon(\bar{\theta}, \xi), i\in [K_\varphi]}
\left(\max\{a_{2i}, 0\}^\top \widecheck{f}(\theta, \xi; I_1) - \max\{-a_{2i}, 0\}^\top \widehat{f}(\theta, \xi; I_2) + b_{2i}\right)\\
&\leq \max_{
(I_1, I_2) \in \mathcal I_g^\varepsilon(\bar{\theta}, \xi) \times \cI_h^\varepsilon(\bar{\theta}, \xi), i\in [K_\varphi]} \left\{
a_{2i}^\top f(\theta, \xi) + b_{2i}
\right\} = \varphi_2(f(\theta, \xi), \xi),
\end{align*}
indicating $\mathcal V_2 = \varphi_2(f(\theta, \xi), \xi)$.
Therefore, we derive~\eqref{eq:=F} and directly obtain $\widehat{\kappa} = 0$. 
For~\ref{enum:algo:phi2}, when $\theta$ is sufficiently close to $\bar{\theta}$, we have 
\begin{align*}
    \HatFmin(\theta, \xi; \bar{\theta})
    &\coloneqq \min_{(I_1, I_2) \in \cI^\varepsilon(\bar{\theta}, \xi)}\widehat{F}(\theta, \xi; \bar{\theta}, (I_1, I_2))\\
    &= \min_{I_2 \in \cI^\varepsilon_h(\bar{\theta}, \xi)}\varphi^{\uparrow}(\widehat{f}(\theta, \xi; I_2), \xi) 
        + \min_{I_1\in \cI^\varepsilon_g(\bar{\theta}, \xi)}\varphi^{\downarrow}(\widecheck{f}(\theta, \xi; I_1), \xi)\\
    &= \varphi^{\uparrow}\left(\min_{I_2 \in \cI^\varepsilon_h(\bar{\theta}, \xi)} \widehat{f}(\theta, \xi; I_2), \xi\right)
        + \varphi^{\downarrow}\left(\max_{I_1\in \cI^\varepsilon_g(\bar{\theta}, \xi)} \widecheck{f}(\theta, \xi; I_1), \xi\right)\\
    &= F(\theta, \xi),
\end{align*}
also implying $\widehat{\kappa} = 0$.
For~\ref{enum:algo:phi3}, when $\theta$ is sufficiently close to $\bar{\theta}$, we obtain
\[
    \HatFmin(\theta, \xi; \bar{\theta})
    = \nabla \varphi(f(\bar{\theta}, \xi), \xi)^\top (f(\theta, \xi) - f(\bar{\theta}, \xi)) + \varphi(f(\bar{\theta}, \xi), \xi)
    + \frac{1}{2}L_{\nabla\varphi}L_f^2 \norm{\theta - \bar{\theta}}^2.
\]
Thus, by the second-order Taylor expansion of $\varphi(\cdot, \xi)$ at $f(\bar{\theta}, \xi)$,
we obtain
\[
    \HatFmin(\theta, \xi; \bar{\theta}) - F(\theta, \xi) 
    = \frac{1}{2}L_{\nabla\varphi}L_f^2 \norm{\theta - \bar{\theta}}^2 - \frac{1}{2}v^\top\nabla^2\varphi(f(\bar{\theta}, \xi), \xi)v + o(\norm{v}^2),
\]
where $\norm{v} = \norm{f(\theta, \xi) - f(\bar{\theta}, \xi)}\leq L_f\norm{\theta - \bar{\theta}}$ by the Lipschitz continuity of $f$. 
By dividing the both sides by $\norm{\theta - \bar{\theta}}^2$ and taking $\limsup$ as $\theta\rightarrow\bar{\theta}$ and $\theta\neq \bar{\theta}$, we obtain that $\widehat{\kappa}$ is finite.\qedsymbol

\endproof

\section{Details of Convergence Analysis in Section~\ref{subsec:the surrogation and the algorithm}}\label{EC:sec:proofs for Section Convergence} %

\subsection{Proof of Proposition~\ref{prop:convergence:relationship-S*-Sds} and equivalence between $\cS^d_{\varepsilon, \rho}$ and $\cS^*$} %
\label{EC:subsec:proofs for Section4-1}

\proof{Proof of Proposition~\ref{prop:convergence:relationship-S*-Sds}.}\label{proof:EC:proof of S*Sd relationships}
Given $\varepsilon \geq 0$ and $\rho > 0$, for any $\theta^* \in \cS^*$
it holds that $\HatFmin(\theta; \theta^*) + \frac{\rho}{2}\norm{\theta - \theta^*}^2 > F(\theta) \geq F(\theta^*) = \HatFmin(\theta^*; \theta^*)$ for any $\theta \neq \theta^*$, showing $\theta^* \in \cS^d_{\varepsilon, \rho}$ and thus $\cS^* \subseteq \cS^d_{\varepsilon, \rho}$.
For any $\theta^d \in \cS^d_{\varepsilon, \rho}$, from Condition~\ref{assum:B4:minSurr dd-consistency} 
and the first-order optimality condition, it follows that $F'(\theta^d; \theta - \theta^d) = \HatFmin(\cdot; \theta^d)'(\theta^d; \theta - \theta^d) \geq 0$ for all $\theta \in \Theta$, yielding $\theta^d \in \cS^d$ and thus $\cS^d_{\varepsilon, \rho} \subseteq \cS^d$.

It remains to show $\cS^d \subseteq \cS^d_{0, \rho}$ for any $\rho > 0$.
Assume by contradiction that there exists $\theta^d \in \cS^d\backslash\cS^d_{0,\widehat{\rho}}$ for some $\widehat{\rho} > 0$.
Then there exist $\{\widehat{I}_s \in \cI^0(\theta^d, \xi^s) \mid s=1, \dots, n\}$ and $\widehat{\theta} (\neq \theta^d) \in \Theta$ such that $\widehat{\theta} \in \mathop{\arg\min}_{\theta \in \Theta}\left\{
    \frac{1}{n}\sum_{s=1}^n\widehat{F}(\theta, \xi; \theta^d, \widehat{I}_s\,) + \frac{\widehat{\rho}}{2}\norm{\theta - \theta^d}^2
\right\}$.
By the strong convexity of the function $\theta \mapsto \frac{1}{n}\sum_{s=1}^n\widehat{F}(\theta, \xi; \theta^d, \widehat{I}_s) + \frac{\widehat{\rho}}{2}\, \norm{\theta - \theta^d}^2$, we have 
$\frac{1}{n}\sum_{s=1}^n\widehat{F}(\cdot, \xi; \theta^d, \widehat{I}_s)'(\theta^d; \widehat{\theta} - \theta^d) < 0$. Thus, there exists $t_0 > 0$ such that
\begin{align*}
    F(\theta^d + t(\widehat{\theta} - \theta^d))
    \leq 
    \frac{1}{n}\sum_{s=1}^n\widehat{F}(\theta^d + t(\widehat{\theta} - \theta^d), \xi^s; \theta^d, \widehat{I}_s) %
    < \frac{1}{n}\sum_{s=1}^n\widehat{F}(\theta^d, \xi^s; \theta^d, \widehat{I}_s) = F(\theta^d),\quad \forall\, t \in (0, t_0],
\end{align*}
where 
the equality holds by the touching condition~\ref{enum:algo:P1:touching}.
This contradicts that $\theta^d$ is a d-stationary point.
Therefore, we have $\cS^d \subseteq \cS^d_{0, \rho}$, which completes the proof.\qedsymbol
\endproof

As an extension, we show the equivalence between $\cS^d_{\varepsilon, \rho}$ and $\cS^*$ with sufficiently large $\varepsilon$ and suitable $\rho$ in Proposition~\ref{prop:EC:relationship-S* equiv Sder}.
The following assumption is needed:  %
\begin{enumerate}
    [label=(EC\arabic*), series=EC, labelwidth=!, labelindent=0pt]
    \item\label{assum:B4:F-gap-global-d-stat} 
    There exists $M_\Delta > 0$ such that $F(\theta^d) \geq \min_{\theta\in \Theta} F(\theta) + M_\Delta$ for any $\theta^d \in \cS^d\backslash \cS^*$.
\end{enumerate}
The assumption is naturally satisfied when the outer function $\varphi$ of $F$ is in~\ref{enum:algo:phi1} or~\ref{enum:algo:phi2}. 
Specifically, from the proof of Claim~\ref{claim:EC:generalized QG for strongly convex varphi} afterward, we see that in these two cases, $\Theta$ can be decomposed into finitely many convex regions, on each of which any d-stationary point $\theta^d \in \cS^d$ is a minimizer of $F$ over that region.
Since the number of such regions is finite, the function $F$ takes finitely many values on $\cS^d$, verifying Assumption~\ref{assum:B4:F-gap-global-d-stat}.
Discussions on the finite number of d-stationary values for more general nonconvex problems with piecewise and composite structures can be seen in \citet{cui2018finite}.
\begin{proposition}\label{prop:EC:relationship-S* equiv Sder}
    Suppose Assumptions~\ref{assum:B1:Touching}-\ref{assum:B4:minSurr dd-consistency} hold.
    For problem~\eqref{eq:algo:formulation:finite sum primal problem} with $\max_{\theta\in\Theta}\norm{\theta}_\infty = \mu$,
    the following statements hold: 
    \begin{enumerate}[label={(\alph*)}, labelwidth=!, labelindent=0pt]
        \item With $\varphi$ in~\ref{enum:algo:phi1} or~\ref{enum:algo:phi2}, $\cS^d_{\varepsilon, \rho} = \cS^*$ for any $\varepsilon \in [\varepsilon_0, +\infty)$ and $\rho \in (0, \frac{M_{\Delta}}{2 q\mu^2}]$, where $\varepsilon_0\coloneqq 2\mu\sqrt{q}\sqrt{\max_{s\in [n]}\norm{x^s}^2 + 1}$, and $M_\Delta$ is the constant defined in~\ref{assum:B4:F-gap-global-d-stat}.
        \item Suppose that Assumption~\ref{assum:B4:F-gap-global-d-stat} holds with $M_{\Delta} > 2\kappa q\mu^2$, and that the minimum surrogate $\widehat{F}^\varepsilon_{\min}(\theta, \xi; \bar{\theta})$ satisfies Condition~\ref{enum:algo:P5:quadratic surrogation gap} for any $\varepsilon \geq 0$ and $\xi \in \Xi$.
        It holds that
        $\cS^d_{\varepsilon, \rho} = \cS^*$ for any $\varepsilon \in [\varepsilon_0, +\infty)$ and $\rho \in (0, \frac{M_{\Delta}}{2 q\mu^2} - \kappa]$.
    \end{enumerate}
\end{proposition}
\proof{Proof of Proposition~\ref{prop:EC:relationship-S* equiv Sder}.}
By Proposition~\ref{prop:convergence:relationship-S*-Sds}, it suffices to prove the inclusion $\cS^d_{\varepsilon, \rho} \subseteq \cS^*$.
Since the two max-affine components of $f(\cdot, \xi)$ for any $\xi \in \Xi$ are Lipschitz continuous with modulus $L_{g,h} \coloneqq \sqrt{\max_{s\in [n]}\norm{x^s}^2 + 1}$, we see that with $\varepsilon\geq\varepsilon_0 = 2 L_{g,h} \sqrt{q}\mu$, the set $\cI^\varepsilon(\bar{\theta}) \coloneqq \{\{I^\varepsilon_s\}_{s=1}^n \mid I^\varepsilon_s \in \cI(\bar\theta, \xi^s), \,s=1,\dots,n\}$ contains all index combinations from the $\varepsilon$-active index sets $\{\cI(\bar\theta, \xi^s)\}_{s=1}^n$, that is, $\cI^\varepsilon(\bar{\theta}) = [K]^{2n}$ for any $\bar{\theta} \in \Theta$.

For the first statement,
assume by contradiction that there exists $\widehat{\theta}^d \in \cS^d_{\widehat{\varepsilon}, \widehat{\rho}} \backslash \cS^*$ for some $\widehat{\varepsilon} \in [\varepsilon_0, +\infty)$ and $\widehat{\rho} \in (0, \frac{M_\Delta}{2q\mu^2}]$.
Since the set $\cI^{\widehat{\varepsilon}}(\widehat{\theta}^d)$ contains all index combinations for the inner PA functions, by the surrogation construction for~\ref{enum:algo:phi1} and~\ref{enum:algo:phi2} in Online Appendix~\ref{EC:sec:surrogation for sec3}, for any $\theta^* \in \cS^*$ there exists $\{I^{\widehat{\varepsilon}}_s\}_{s=1}^n \in \cI^{\widehat{\varepsilon}}(\widehat{\theta}^d)$ such that $\widehat{F}^{\widehat{\varepsilon}}_{\min}(\theta^*; \widehat{\theta}^d) = \frac{1}{n}\sum_{s=1}^n\widehat{F}(\theta^*,\xi^s; \widehat{\theta}^d, I^{\widehat{\varepsilon}}_s) = \min_{\theta\in\Theta}F(\theta)$.
It follows that
\[
    \widehat{F}^{\widehat{\varepsilon}}_{\min}(\theta^*; \widehat{\theta}^d) + \frac{\widehat{\rho}}{2}\, \big\lVert \,{\theta^* - \widehat{\theta}^d}\, \big\rVert^2 \leq \min_{\theta\in\Theta}F(\theta) + M_\Delta \leq F(\widehat{\theta}^d)  = \widehat{F}^{\widehat{\varepsilon}}_{\min}(\widehat{\theta}^d; \widehat{\theta}^d),
\]
which contradicts $\widehat{\theta}^d \in \cS^d_{\widehat{\varepsilon}, \widehat{\rho}}$.
For the second statement, in a similar way, assume that there exists $\widehat{\theta}^d \in \cS^d_{\widehat{\varepsilon}, \widehat{\rho}} \backslash \cS^*$ for some $\widehat{\varepsilon} \in [\varepsilon_0, +\infty)$ and $\widehat{\rho} \in (0, \frac{M_\Delta}{2 q\mu^2} - \kappa ]$. By Condition~\ref{enum:algo:P5:quadratic surrogation gap}, we have that for $\theta^* \in \cS^*$,
\[
    \widehat{F}^{\widehat{\varepsilon}}_{\min}(\theta^*; \widehat{\theta}^d) + \frac{\widehat{\rho}}{2}\, \big\lVert \,{\theta^* - \widehat{\theta}^d}\, \big\rVert ^2 
    \leq F(\theta^*) + \frac{\widehat{\rho} + \kappa}{2} \, \big\lVert \,{\theta^* - \widehat{\theta}^d}\, \big\rVert^2
    \leq \min_{\theta\in\Theta}F(\theta) + M_\Delta 
    \leq F(\widehat{\theta}^d).
\]
Again, it leads to a contradiction, which completes the proof.\qedsymbol
\endproof

\subsection{Proof of Proposition~\ref{prop:convergence:convergence rate of bar-r}}\label{EC:subsec:proofs for Sec4-2}

The following lemma from \citet{deng2021minibatch} is rephrased in our context for the proof of Proposition~\ref{prop:convergence:convergence rate of bar-r}. 
\begin{lemma}
    [\citet{deng2021minibatch}, Theorem 3.2]
    \label{lemma:stability:deng-minibatch}
    Suppose that a stochastic function $\phi(\cdot, \xi)$ is convex and Lipschitz continuous on $\Theta$ with modulus $L_{\phi}$ for any $\xi \in \Xi$ and $\bar{\theta} \in \Theta$, and that samples $\Xi_N \coloneqq \{\widetilde{\xi}^s\}_{s = 1}^N$ are i.i.d. generated from $\Xi$.
    For $\eta > 0$, denote $\bar{\theta}_+ \coloneqq \mathop{\arg\min}_{\theta\in \Theta}\left\{
        \frac{1}{N}\sum_{s=1}^N \phi(\theta, \widetilde{\xi}^s) + \frac{\eta}{2}\norm{\theta - \bar{\theta}}^2
    \right\}$. Then it holds that
    \[
        \abs{\EE_{\Xi_N}\left[\frac{1}{N}\sum_{s=1}^N \phi(\bar{\theta}_+, \widetilde{\xi}^s) - \EE_{\widetilde\xi_n}\left[\phi(\bar{\theta}_+, \widetilde\xi_n)\right]\right]}
        \leq \frac{2 L_{\phi}^2}{\eta N}.
    \]
\end{lemma}
\proof{Proof of Proposition~\ref{prop:convergence:convergence rate of bar-r}.}
    \label{proof:prop::algo:unconstrained:convergence rate of M mapping}
    \newcommand{\tIJen}{\widetilde{\bI}^{\varepsilon}_{\theta^\nu}}
    \newcommand{\tIJn}{\widetilde{\bI}^0_{\theta_\nu}}
    \newcommand{\vthetan}{P^\rho_{\wtbI^\varepsilon_\nu}(\theta^\nu)}
    \newcommand{\vthetanRV}{P^\rho_{\widetilde{\bI}^\varepsilon_\nu}(\theta^\nu)}
    \newcommand{\vthetanI}{P^\rho_{\wtbI^0_\nu}(\theta^\nu)}
    \newcommand{\vthetanIRV}{P^\rho_{\widetilde{\bI}^0_\nu}(\theta^\nu)}
    Let $\tIJen = \{\widetilde{I}^\varepsilon_\xi\}_{\xi \in \Xi}$ be a random index combination uniformly drawn from $\cI^\varepsilon(\theta^\nu)$, $\{\widetilde{\xi}^s_\nu\}_{s=1}^{N_{\nu}}$ be random vectors independently and uniformly drawn from $\Xi$. 
    Let $\cG_{\nu-1}$ denote the $\sigma$-algebra generated by all past samples $\{\{\widetilde{\xi}_{t}^s\}_{s = 1}^{N_t}\}_{t = 0}^{\nu-1}$ and combinations $\{\widetilde{\bI}^\varepsilon_{t}\}_{t=0}^{\nu - 1}$, and $\cG_{\nu - 1}^+$ denote the $\sigma$-algebra generated by $\{\{\widetilde{\xi}_{t}^s\}_{s = 1}^{N_t}\}_{t = 0}^{\nu-1}$ and $\{\widetilde{\bI}^\varepsilon_{t}\}_{t=0}^{\nu}$. Given the current iterate $\theta^\nu$ and the index combination $\IJen \in \cI^\varepsilon(\theta^\nu)$, by the definition of $P_{\IJen}^\rho(\cdot)$ in~\eqref{eq:mapping:barP}, we have
    \begin{equation}
        \label{eq:proof:algo:unconstr:descent}
        \min\left\{
            \widehat{F}(\vthetan; \theta^\nu, \IJen), F(\theta^{\nu})
        \right\} + \frac{\rho}{2}\, \big\lVert \,{\vthetan - \theta^\nu}\, \big\rVert^2 
        \leq %
         F(\theta^\nu).
    \end{equation}
    With the notation $\EE_{N_{\nu}}$ representing the expectation with respect to $\{\widetilde{\xi}_\nu^s\}_{s = 1}^{N_{\nu}}$, by rearranging terms of~\eqref{eq:proof:algo:unconstr:descent} and taking the conditional expectation on both sides, we obtain that
    \begin{align}\label{eq:EC:rate:two terms}
            \frac{\rho}{2}\EE\left[\, \big\lVert \,{\vthetanRV - \theta^\nu}\, \big\rVert ^2 \,\middle\vert\, \cG_{\nu - 1} \right] 
            &\leq F(\theta^\nu) - \EE\left[\min\left\{\widehat{F}(\vthetanRV; \theta^\nu, \tIJen), F(\theta^\nu)\right\} \,\middle\vert\, \cG_{\nu - 1} \right]\notag\\
            &= \underbrace{F(\theta^\nu) - W(\theta^{\nu})}_{\cE_1^\nu}
            + \underbrace{W(\theta^{\nu}) - \EE\left[\min\left\{\widehat{F}(\vthetanRV; \theta^\nu, \tIJen), F(\theta^\nu)\right\}\,\middle\vert\, \cG_{\nu - 1}
            \right]}_{\cE_2^\nu},
    \end{align}
    where $W(\theta^\nu)\coloneqq
            \EE_{\tIJen}\left[\min\left\{
                \EE_{N_{\nu}}\left[\widehat{F}_{\nu}(\bthetan; \theta^\nu, \widetilde{\bI}^\varepsilon_{\nu}) + \frac{\eta}{2}\, \big\lVert \,{\bthetan - \theta^\nu}\, \big\rVert ^2 \,\middle\vert\, \cG_{\nu - 1}^+ \right], F(\theta^\nu)
            \right\} \,\middle\vert\, \cG_{\nu - 1} \right]$.
            
    From the update rule in Step~5 of Algorithm~\ref{ALGO:unconstrained case}, we have $\widehat{F}_{\nu}(\bthetan; \theta^\nu, \wtbI^\varepsilon_{\nu}) + \frac{\eta}{2}\, \big\lVert \,{\bthetan - \theta^\nu}\, \big\rVert^2 \geq F_{\nu}(\theta^{\nu + 1})$ and $F_{\nu}(\theta^\nu)\geq F_{\nu}(\theta^{\nu + 1})$.
    It follows that
    \begin{align*}
        W(\theta^\nu)
        &\geq 
        \EE_{\tIJen}\left[\min\left\{
            \EE_{N_{\nu}} \left[F_{\nu}(\theta^{\nu + 1}) \mid \cG_{\nu - 1}^+\right], F(\theta^\nu)
        \right\} \mid \cG_{\nu - 1}\right]
        = \EE\left[F_{\nu}(\theta^{\nu + 1})\mid \cG_{\nu - 1}\right].
    \end{align*}
    With $\Theta$ being a compact set, let $\max_{\theta\in\Theta} \norm{\theta}_\infty = \mu$. Define 
    \[
        \widehat{\Phi}_{\HK} \coloneqq \{\phi(\cdot, \theta):\Xi \rightarrow \RR \mid \phi(\xi, \theta) = \varphi(f(\theta, \xi), \xi),\, \xi \in \Xi, f(\theta, \cdot)\in\HK, \theta\in\Theta\}.
    \]
    From finiteness of $\Xi$ and continuity of $\varphi(\cdot, \xi)$, it follows that for any $\phi \in \widehat{\Phi}_{\HK}$, $\phi(\xi,\cdot)$ is Lipschitz continuous on the compact set $\Theta$ with uniform modulus, and that $\phi(\cdot, \theta)$ is uniformly bounded.
    Therefore, for the term $\cE_1^\nu$, there exists a constant $C_1 > 0$ such that
    \begin{align}
    \label{eq:proof:algo:convergence:T1bound}
        \cE_1^\nu
        &\leq F(\theta^\nu) - \EE\left[F_{\nu}(\theta^{\nu+1})\mid \cG_{\nu - 1}\right]\notag\\
        & = F(\theta^\nu) - \EE\left[F(\theta^{\nu+1})\mid \cG_{\nu - 1}\right] 
        + \EE\left[F(\theta^{\nu+1})\mid \cG_{\nu - 1}\right] - \EE\left[F_{\nu}(\theta^{\nu + 1})\mid \cG_{\nu - 1}\right]\notag\\
        &\leq F(\theta^\nu) - \EE\left[F(\theta^{\nu+1})\mid \cG_{\nu - 1}\right] + \EE\left[\sup_{\theta\in \Theta}\abs{F(\theta) - F_{\nu}(\theta)}\,\middle\vert\, \cG_{\nu - 1}\right]\notag\\
        &\leq F(\theta^\nu) - \EE\left[F(\theta^{\nu+1})\mid \cG_{\nu - 1}\right] + 2\cR_{N_{\nu}}(\widehat{\Phi}_{\HK})\notag\\
        &\leq F(\theta^\nu) - \EE\left[F(\theta^{\nu+1})\mid \cG_{\nu - 1}\right] + C_1\sqrt{\frac{\log N_{\nu}}{N_{\nu}}},
    \end{align}
    where the third inequality follows from inequalities (3.8)-(3.13) in \citet{ECmohri2018foundations} that are part of the proof of Lemma~\ref{lemma:EC:Gen Error bound with Rademacher complex}, and the last inequality is derived by Lemma~\ref{lemma:EC:Rademacher of Param functions}.

    Next, we bound the term $\cE_2^\nu$. Given $\tIJen$ and $\{\xi_{\nu}^s\}_{s = 1}^{N_{\nu}}$, by the optimality of $\bthetan$, we have
    \[
        \widehat{F}_{\nu}(\bthetan; \theta^\nu, \wtbI^\varepsilon_{\nu}) + \frac{\eta}{2}\, \big\lVert \,{\bthetan - \theta^{\nu}}\, \big\rVert^2 
        \leq 
        \widehat{F}_{\nu}(\vthetan; \theta^{\nu}, \wtbI^\varepsilon_{\nu}) + \frac{\eta}{2}\, \big\lVert \,\vthetan - \theta^{\nu}\, \big\rVert^2.
    \]
    It follows that
    \begin{align*}
            W(\theta^\nu)
            &\leq
            \EE_{\tIJen}\left[\min\left\{
                \EE_{N_{\nu}}\left[\widehat{F}_{\nu}(\vthetanRV; \theta^\nu, \widetilde{\bI}^\varepsilon_{\nu}) + \frac{\eta}{2}\, \big\lVert \,\vthetanRV - \theta^\nu\, \big\rVert ^2 \,\middle\vert\, \cG_{\nu - 1}^+\right], F(\theta^\nu)
            \right\}\,\middle\vert\, \cG_{\nu - 1}\right]\\
            &\leq
            \EE_{\tIJen}\left[\min\left\{
                \EE_{N_{\nu}}\left[\widehat{F}_{\nu}(\vthetanRV; \theta^\nu, \widetilde{\bI}^\varepsilon_{\nu})\,\middle\vert\, \cG_{\nu - 1}^+\right], F(\theta^\nu)
            \right\} + \frac{\eta}{2}\, \big\lVert \,\vthetanRV - \theta^\nu\, \big\rVert ^2\,\middle\vert\, \cG_{\nu - 1}\right]\\
            &=
            \EE\left[\min\left\{
                \widehat{F}(\vthetanRV; \theta^\nu, \tIJen), F(\theta^\nu)
            \right\}\,\middle\vert\, \cG_{\nu - 1}\right] + \frac{\eta}{2}\EE\left[\, \big\lVert \,\vthetanRV - \theta^\nu\, \big\rVert ^2 \,\middle\vert\, \cG_{\nu - 1} \right],
    \end{align*}
    where the equality holds since both $\theta^\nu$ and $\vthetanRV$ are independent of $\{\widetilde{\xi}_\nu^s\}_{s = 1}^{N_{\nu}}$. 
    Then we obtain
    \begin{align}
    \label{eq:proof:algo:convergence:T2bound}
        \cE_2^\nu
        & \leq \frac{\eta}{2}\EE\left[\, \big\lVert \,\vthetanRV - \theta^\nu\, \big\rVert ^2 \,\middle\vert\, \cG_{\nu - 1} \right].
    \end{align}
    
    Combining bounds~\eqref{eq:proof:algo:convergence:T1bound} and~\eqref{eq:proof:algo:convergence:T2bound} of $\cE_1^\nu$ and $\cE_2^\nu$ and rearranging terms, we derive that
    \[
        \frac{\rho - \eta}{2}\EE\left[\, \big\lVert \,{\vthetanRV - \theta^{\nu}}\, \big\rVert^2\,\middle\vert\, \cG_{\nu - 1}\right] 
        \leq F(\theta^{\nu}) - \EE\left[F(\theta^{\nu + 1})\,\middle\vert\, \cG_{\nu - 1}\right] + C_1\sqrt{\frac{\log N_{\nu}}{N_{\nu}}}.
    \]
    By applying Jensen's inequality, taking full expectation on both sides, and summing inequalities up over $\nu = 0, \dots, T-1$, we have
    \[
        \frac{\rho - \eta}{2}\sum_{\nu = 0}^{T-1} \big( \, \EE\,[\bar{r}^{\varepsilon, \rho}(\theta^\nu)] \, \big)^2
        \leq \frac{\rho - \eta}{2} \sum_{\nu = 0}^{T-1} \EE \, \big\lVert \,{\vthetanRV - \theta^{\nu}}\, \big\rVert^2
        \leq F(\theta^0) - \EE\left[F(\theta^{T})\right] + C_1\cdot\sum_{\nu = 0}^{T-1}\sqrt{\frac{\log N_{\nu}}{N_{\nu}}}.
    \]
    Dividing both sides by $T$, applying Jensen's inequality again, 
    and using the fact that $\EE\left[F(\theta^T)\right] \geq \min_{\theta\in\Theta} F(\theta)$, we finally obtain that
    \[
        \left(\EE\,\big[\,\bar{r}^{\varepsilon, \rho}(\widetilde{\theta}^T)\,\big]\right)^2 
        \leq \frac{F(\theta^0) - \min_{\theta\in\Theta} F(\theta) + C_1 \cdot \sum_{\nu = 0}^{T-1}\sqrt{\log N_{\nu}/N_{\nu}}}{\frac{\rho - \eta}{2}\cdot T}.
    \]

    Now we consider the special case with $\varepsilon = 0$ and $\eta > 0$.
    In such a case, $\bthetan = \theta^{\nu + 1}$, and the intermediate term $W(\theta^\nu)$ degenerates into
    \[
        W(\theta^\nu) 
        = \EE\left[\widehat{F}_{\nu}(\theta^{\nu + 1}; \theta^\nu, \widetilde{\bI}^0_{\nu}) 
        + \frac{\eta}{2}\norm{\theta^{\nu + 1} - \theta^\nu}^2\,\middle\vert\, \cG_{\nu - 1}\right].
    \]
    The term $\cE_2^\nu$ is still bounded by~\eqref{eq:proof:algo:convergence:T2bound}.
    For $\cE_1^\nu$, we can derive a tighter bound in the complexity of the sampling size $N_{\nu}$.
    Specifically, there exists $C_2 > 0$ such that
    \begin{align*}
            \cE_1^\nu 
            &= F(\theta^\nu) - \EE\left[\widehat{F}_{\nu}(\theta^{\nu + 1}; \theta^\nu, \widetilde{\bI}^0_{\nu}) 
            + \frac{\eta}{2}\norm{\theta^{\nu + 1} - \theta^\nu}^2 \,\middle\vert\, \cG_{\nu - 1} \right]\\
            &\leq F(\theta^\nu) - \EE\left[\widehat{F}_{\nu}(\theta^{\nu + 1}; \theta^\nu, \widetilde{\bI}^0_{\nu})\,\middle\vert\, \cG_{\nu - 1}\right]\\
            &= F(\theta^\nu) - \EE\left[\widehat{F}(\theta^{\nu + 1}; \theta^\nu, \tIJn) \,\middle\vert\, \cG_{\nu - 1}\right]
            + \EE\left[\widehat{F}(\theta^{\nu + 1}; \theta^\nu, \tIJn) \,\middle\vert\, \cG_{\nu - 1}\right]
            - \EE\left[\widehat{F}_{\nu}(\theta^{\nu + 1}; \theta^\nu, \tIJn)\,\middle\vert\, \cG_{\nu - 1}\right]\\
            &\leq F(\theta^\nu) - \EE\left[F(\theta^{\nu + 1})\,\middle\vert\, \cG_{\nu - 1}\right] + \EE\left[
                    \widehat{F}(\theta^{\nu + 1}; \theta^\nu, \tIJn) - \widehat{F}_{\nu}(\theta^{\nu + 1}; \theta^\nu, \widetilde{\bI}^0_{\nu})
                \,\middle\vert\, \cG_{\nu - 1}
                \right]\\
            &\leq F(\theta^\nu) - \EE\left[F(\theta^{\nu + 1})\,\middle\vert\, \cG_{\nu - 1}\right] + \frac{C_2}{N_{\nu}},
    \end{align*}
    where the last inequality holds by Lemma~\ref{lemma:stability:deng-minibatch}, and the Lipschitz continuity of surrogates follows from the bi-variate continuity of $\widehat F(\cdot, \xi;\cdot, I)$, compactness of $\Theta$, and finiteness of $\Xi$.
    Then the second inequality in the proposition can be proved with arguments similar to those above.\qedsymbol
\endproof

\subsection{Proofs of Theorems~\ref{thm:convergence:convergence results} and~\ref{thm:convergence:convergence results:PA}}\label{EC:subsec:analytical details on convergence}
In this section, we prove the convergence results in Theorems~\ref{thm:convergence:convergence results} and~\ref{thm:convergence:convergence results:PA}.
As in Section~\ref{subsec:the surrogation and the algorithm}, we suppose throughout this section that 
$\varphi$ is directionally differentiable and locally Lipschitz continuous, $\Theta$ is a compact convex set, and the surrogate family $\{\widehat F(\cdot, \xi; \bar{\theta}, I)\}_{I \in \cI^\varepsilon(\bar{\theta}, \xi)}$ of $F(\cdot, \xi)$ satisfies Assumptions~\ref{assum:B1:Touching}-\ref{assum:B4:minSurr dd-consistency} for each $\xi \in \Xi$.

The convergence analysis requires the following two residual functions.
The first residual of a reference point $\bar{\theta} \in \Theta$ with $\varepsilon \geq 0$ and $\rho > 0$ measures the proximal mapping gap with respect to the  surrogate function $\HatFmin(\theta; \bar{\theta}) \coloneqq \min_{\bI^\varepsilon \in \cI^\varepsilon(\bar{\theta})}\widehat{F}(\theta; \bar{\theta}, \bI^\varepsilon)$, which is defined as
\begin{equation}\label{eq:definition of the residual}
    r^{\varepsilon, \rho}(\bar{\theta}) \coloneqq \max_{\theta \in \cP_{\min}^{\varepsilon, \rho}(\bar{\theta})}\norm{\theta - \bar{\theta}},
    \text{\ where\ }
    \cP_{\min}^{\varepsilon, \rho}(\bar{\theta}) \coloneqq 
    \mathop{\arg\min}_{\theta \in \Theta} \left\{
        \HatFmin(\theta; \bar{\theta}) + \frac{\rho}{2}\norm{\theta - \bar{\theta}}^2
    \right\}.
\end{equation}
Clearly, $\theta^d \in \Theta$ is a composite $(\varepsilon, \rho)$-strongly d-stationary point of problem~\eqref{eq:algo:formulation:finite sum primal problem} if and only if $r^{\varepsilon, \rho}(\theta^d) = 0$.
The residual $r^{\varepsilon, \rho}$ works as a criterion for the subsequential convergence to composite strong d-stationarity and is also exploited to establish the error-bound analysis.
The other useful residual is defined as $r_d(\theta) \coloneqq \underset{d \in \cT(\theta; \Theta), \norm{d}\leq 1}{\max} \left\{- F'(\theta; d)\right\}$,
where $\cT(\theta; \Theta)$ is the tangent cone of $\Theta$ at $\theta \in \Theta$.
This residual is utilized to derive the global error bound in terms of $r^{\varepsilon, \rho}$ in the second statement of Proposition~\ref{prop:EC:limiting stationarity and error bounds}.
In Figure~\ref{fig:proof:convergence roadmap}, we summarize a roadmap for the convergence analysis of the ESMM algorithm.

\begin{figure}[htbp]
    \FIGURE
    {\includegraphics[width=1\textwidth]{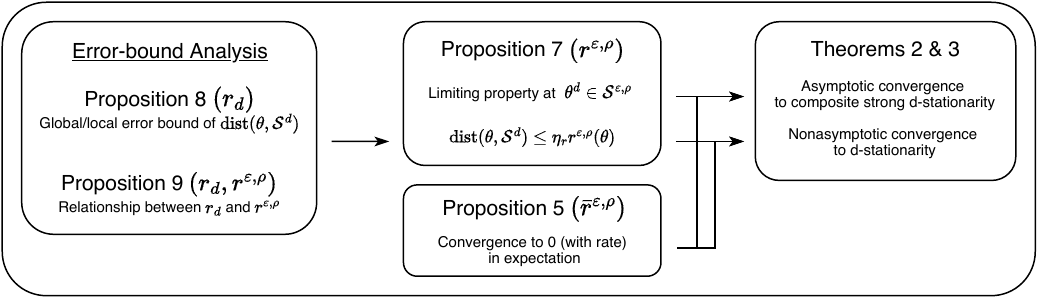}}
    {A roadmap for convergence analysis of ESMM\label{fig:proof:convergence roadmap}}
    {}
\end{figure}

We first provide in the following proposition the limiting property of the residual $r^{\varepsilon, \rho}$ and the global error bound of distance to $\cS^d$ in terms of $r^{\varepsilon, \rho}$. 

\begin{proposition}\label{prop:EC:limiting stationarity and error bounds}
    For problem~\eqref{eq:algo:formulation:finite sum primal problem}, the following statements hold:
    \begin{enumerate}[label={(\alph*)}, labelwidth=!, labelindent=0pt]
        \item For ${\varepsilon} > 0$ and ${\rho} > 0$, let $\theta_\infty$ be a limit point of a sequence $\{\theta_t\}_t$ in $\Theta$ satisfying $\lim_{t\rightarrow \infty} r^{{\varepsilon}, {\rho}}(\theta_t) = 0$. 
        Then $\theta_\infty$ is a composite $(\varepsilon', \rho')$-strongly d-stationary point for any $\varepsilon' \in [0, \varepsilon)$ and $\rho' \in (\rho, +\infty)$.
        
        \item For $\varepsilon \geq 0$ and $\rho > 0$, if $\varphi$ is piecewise affine and $\Theta$ is a polytope, there exists a constant $\eta_r > 0$ such that $\distSd{\theta} \leq \eta_r r^{\varepsilon, \rho}(\theta)$ for any $\theta \in \Theta$.
    \end{enumerate}
\end{proposition}

Before proving this proposition, we show that Theorem~\ref{thm:convergence:convergence results} could be derived straightforwardly by combining Propositions~\ref{prop:convergence:convergence rate of bar-r} and~\ref{prop:EC:limiting stationarity and error bounds}, similarly for Theorem~\ref{thm:convergence:convergence results:PA}. 
\proof{Proof of Theorem~\ref{thm:convergence:convergence results}.}
    According to the definition in~\eqref{eq:definition of the residual}, for any $\theta \in \cP^{\varepsilon, \rho}_{\min}(\bar{\theta})$, there exists $\bI^\varepsilon \in \cI^\varepsilon(\bar{\theta})$ such that $\theta = P^\rho_{\bI^\varepsilon}(\bar{\theta})$.
    Due to the finite-sum structure of the objective $F$, there exists a uniform constant $C_{\cI} > 0$ such that the cardinality $\abs{\cI^\varepsilon(\bar{\theta})}\leq C_{\cI}$ for any $\bar{\theta} \in \Theta$, and thus ${r}^{\varepsilon, \rho}(\bar{\theta}) \leq \abs{\cI^\varepsilon(\bar{\theta})} \bar{r}^{\varepsilon, \rho}(\bar{\theta})\leq C_{\cI} \bar{r}^{\varepsilon, \rho}(\bar{\theta})$, which implies that
    $\EE\,[{r}^{\varepsilon,\rho}(\widetilde{\theta}^{\,T})]$ converges to 0 at least at the same rate as $\EE\,[\bar{r}^{\varepsilon, \rho}(\widetilde{\theta}^{\,T})]$ in terms of $\{N_{\nu}\}_\nu$ and $T$.
    Therefore, the statements of the theorem follow from Propositions~\ref{prop:convergence:convergence rate of bar-r} and~\ref{prop:EC:limiting stationarity and error bounds}.\qedsymbol
\endproof

We next present the proof of Proposition~\ref{prop:EC:limiting stationarity and error bounds}. 
The overall scheme of the error-bound analysis follows the error-bound theory developed in~\citet{ECliu2022solving} for the post-convergence analysis of the SMM algorithm.
To justify the global error bound in the second statement of Proposition~\ref{prop:EC:limiting stationarity and error bounds}, we first give the error bounds of distance to $\cS^d$ in terms of the residual $r_d$ in Proposition~\ref{prop:proof:error bounds by r-d} and then relate $r_d$ with $r^{\varepsilon, \rho}$ in Proposition~\ref{prop:EC:connecting r-d and r:liu2022}.
\begin{proposition}\label{prop:proof:error bounds by r-d}
    For problem~\eqref{eq:algo:formulation:finite sum primal problem}, if $\varphi$ is piecewise affine and $\Theta$ is a polyhedral set, 
    then there exists $\eta_d > 0$ such that
    \[
        \distSd{\theta} \leq \eta_d r_d(\theta),\quad \forall\, \theta\in \Theta.
    \]
\end{proposition}

 \proof{Proof of Proposition~\ref{prop:proof:error bounds by r-d}.}\label{proof:prop:proof:error bounds by r-d}
    Since the objective $F$ is a piecewise affine function on a polyhedral set $\Theta$, $\Theta$ can be decomposed into finitely many polyhedral sets such that $F$ is affine on each of them, showing that $r_d(\cdot)$ takes finitely many values on $\Theta$.
    Notice that $\theta \in \cS^d$ if and only if $r_d(\theta) = \max_{\dir} \{-F'(\theta; d)\} = 0$.
    Hence, there exists $l_0 > 0$ such that $r_d(\theta)\geq l_0$ for any $\theta\in \Theta\backslash\cS^d$.
    By setting $\eta_d \coloneqq \frac{2}{l_0}\max_{\theta \in \Theta}\norm{\theta}$, we deduce $\distSd{\theta} \leq 2\max_{\theta \in \Theta}\norm{\theta} \leq \eta_d r_d(\theta)$.
    \qedsymbol
\endproof

\begin{proposition}%
    \label{prop:EC:connecting r-d and r:liu2022}
    Suppose that the minimum surrogate $\widehat{F}^\varepsilon_{\min}(\theta, \xi; \bar{\theta})$ satisfies the quadratic surrogation gap condition~\ref{enum:algo:P5:quadratic surrogation gap} for any $\varepsilon \geq 0$ and $\xi \in \Xi$.
    For problem~\eqref{eq:algo:formulation:finite sum primal problem} with $\varepsilon \geq 0$ and $\rho > 0$, it holds that for any $\theta \in \Theta$, there exists $\theta' \in \Theta$ such that $\norm{\theta' - \theta} \leq 2\,r^{\varepsilon, \rho}(\theta)$ 
    and $r_d(\theta') \leq \left(2\rho + \frac{5}{2}\kappa\right)\,r^{\varepsilon, \rho}(\theta)$.
\end{proposition}

A similar result to Proposition~\ref{prop:EC:connecting r-d and r:liu2022} is provided in \citet[Lemma 4]{ECliu2022solving} under the assumption that the relevant  surrogate function satisfies Condition~\ref{enum:algo:P5:quadratic surrogation gap}.
With the  surrogate function $\HatFmin(\cdot; \bar{\theta})$ associated with $r^{\varepsilon, \rho}(\bar{\theta})$ satisfying~\ref{enum:algo:P5:quadratic surrogation gap}, we can prove Proposition~\ref{prop:EC:connecting r-d and r:liu2022} following the lines of~\citet[Lemma 4]{ECliu2022solving}, and thus the details are omitted here.

\proof{Proof of Proposition~\ref{prop:EC:limiting stationarity and error bounds}.}
    For the first statement, let $\{\theta_t\}_{t\in \cK}$ be a convergent subsequence with the limit point $\theta_\infty$.
    Given $\varepsilon' \in [0, \varepsilon)$, it follows from the Lipschitz continuity of $f$ that $\cI^{\varepsilon'}(\theta_\infty) \subseteq \cI^\varepsilon(\theta_t)$ for sufficiently large $t \in \cK$.
    With $\widehat{P}^\rho_{\bI}(\theta_t) \coloneqq %
    \underset{\theta \in \Theta}{\mathop{\arg\min}}\left\{
        \widehat{F}(\theta; \theta_t, \bI) + \frac{\rho}{2}\norm{\theta - \theta_t}^2
    \right\}$, by definition of the residual $r^{\varepsilon, \rho}$ there exists $\bI^*_t \in \underset{\bI\in \cI^\varepsilon(\theta_t)}{\arg\min} \left\{
        \widehat{F}(\widehat{P}^\rho_{\bI}(\theta_t); \theta_t, \bI) + \frac{{\rho}}{2}\, \big\lVert \,{\widehat{P}^\rho_{\bI}(\theta_t) - \theta_t}\, \big\rVert ^2\right\}$ such that $r^{\varepsilon, \rho}(\theta_t) = \, \big\lVert \,{\widehat{P}^\rho_{\bI^*_t}(\theta_t) - \theta_t}\, \big\rVert$. Then for sufficiently large $t \in \cK$, we have that
    \begin{align}\label{eq:EC:prf-prop-EC4-1}
        \widehat{F}^*(\theta_t) &\coloneqq 
        \min_{\bI \in \cI^{\varepsilon'}(\theta_\infty)}\left\{
        \widehat{F}(\widehat{P}^\rho_{\bI}(\theta_t); \theta_t, \bI) + \frac{{\rho}}{2}\, \big\lVert \,\widehat{P}^\rho_{\bI}(\theta_t) - \theta_t\, \big\rVert ^2\right\}\notag\\
        &\geq \min_{\bI \in \cI^\varepsilon(\theta_t)}\left\{
            \widehat{F}(\widehat{P}^\rho_{\bI}(\theta_t); \theta_t, \bI) + \frac{{\rho}}{2}\, \big\lVert \,\widehat{P}^\rho_{\bI}(\theta_t) - \theta_t\, \big\rVert ^2\right\}\notag\\
        &= \widehat{F}(\widehat{P}^\rho_{\bI^*_t}(\theta_t); \theta_t, \bI^*_t) + \frac{{\rho}}{2}\, \big\lVert \,\widehat{P}^\rho_{\bI^*_t}(\theta_t) - \theta_t\, \big\rVert^2\notag\\
        &\geq F(\widehat{P}^\rho_{\bI^*_t}(\theta_t)) + \frac{{\rho}}{2}\, \big\lVert \,\widehat{P}^\rho_{\bI^*_t}(\theta_t) - \theta_t\, \big\rVert^2.
    \end{align}

    Since $\lim_{t(\in \cK)\rightarrow \infty}r^{{\varepsilon}, {\rho}}(\theta_t) = 0$, we deduce
    $\lim_{t(\in\cK)\rightarrow \infty} \widehat{P}^\rho_{\bI^*_t}(\theta_t) = \theta_\infty$. 
    Additionally, for a fixed index combination $\bI\in \cI^{\varepsilon'}(\theta_\infty)$, the  surrogate function $\widehat{F}(\cdot; \cdot, \bI)$ is continuous on $\Theta\times \Theta$ by assumption, and thus by \citet[Theorem 7.23]{ECshapiro2021lectures}, $\widehat{P}^\rho_{\bI}(\cdot)$ is continuous at $\theta_\infty$.
    It follows from~\eqref{eq:EC:prf-prop-EC4-1} that
    \[
        F(\theta_\infty) 
        \leq \lim_{t(\in \cK)\rightarrow \infty}\widehat{F}^*(\theta_t) 
        = \min_{\bI \in \cI^{\varepsilon'}(\theta_\infty)} \left\{\widehat{F}(\widehat{P}^\rho_{\bI}(\theta_\infty); \theta_\infty, \bI) + \frac{{\rho}}{2} \, \big\lVert \, \widehat{P}^\rho_{\bI}(\theta_\infty) - \theta_\infty \, \big\rVert^2\right\},
    \]
    which leads to $F(\theta_\infty) <\displaystyle{ \min_{\bI \in \cI^{\varepsilon'}(\theta_\infty)\backslash\cI^0(\theta_\infty)} } 
        \left\{\widehat{F}(\widehat{P}^{\rho'}_{\bI}(\theta_\infty); \theta_\infty, \bI) + \frac{{\rho'}}{2}\, \big\lVert \,\widehat{P}^{\rho'}_{\bI}(\theta_\infty) - \theta_\infty\, \big\rVert^2\right\}$
    for $\rho' > \rho$, yielding $r^{\varepsilon', \rho'}(\theta_\infty) = 0$.

    We next prove the second statement.
    If $\varphi$ is piecewise affine, by Proposition~\ref{prop:verify Cases1-3} in Appendix~\ref{EC:sec:surrogation for sec3}, the minimum surrogate function $\widehat{F}^\varepsilon_{\min}(\theta, \xi; \bar{\theta})$ satisfies Condition~\ref{enum:algo:P5:quadratic surrogation gap} for any $\varepsilon\geq 0$ and $\xi \in \Xi$.
    Given $\varepsilon \geq 0$ and $\rho > 0$, by Proposition~\ref{prop:proof:error bounds by r-d} there exists $\eta_d > 0$ such that $\distSd{\theta} \leq \eta_d r_d(\theta)$ for any $\theta \in \Theta$; by Proposition~\ref{prop:EC:connecting r-d and r:liu2022}, for any $\theta \in \Theta$, there exists $\theta' \in \Theta$ such that 
    $\norm{\theta' - \theta} \leq 2 r^{\varepsilon, \rho}(\theta)$ and $r_d(\theta') \leq (2\rho + \frac{5\kappa}{2})r^{\varepsilon, \rho}(\theta)$. Setting $\eta_r = 2 + \eta_d\left(2\rho + \frac{5\kappa}{2}\right)$, we have that
    \begin{equation*}
        \distSd{\theta} 
        \leq \mathrm{dist}(\theta', \cS^d) + \norm{\theta - \theta'}
        \leq \eta_d r_d(\theta') + 2 r^{\varepsilon, \rho}(\theta)
        \leq \eta_r r^{\varepsilon, \rho}(\theta),
    \end{equation*}
    which completes the proof.\qedsymbol
\endproof

\subsection{Nonasymptotic convergence result for general cases}\label{EC:sec:nonasymptotic convergence for general cases}

In this section, we extend the nonasymptotic convergence result for the general case where the outer function $\varphi$ is directionally differentiable and locally Lipschitz continuous, and the surrogate functions of $F$ satisfy Assumptions~\ref{assum:B1:Touching}-\ref{assum:B4:minSurr dd-consistency}.
The error-bound analysis in this case requires the following two assumptions. We denote the closure operator of a set by $\CL{\cdot}$ and $\cP_{\min}^{\varepsilon, \rho}(\bar{\theta}) \coloneqq 
\mathop{\arg\min}_{\theta \in \Theta} \left\{
    \HatFmin(\theta; \bar{\theta}) + \frac{\rho}{2}\norm{\theta - \bar{\theta}}^2
\right\}$.
\begin{enumerate}
    [label=(B\arabic*), resume*=B]
    \item\label{assum:B5:surrogation quadratic gap} The minimum surrogate $\widehat{F}^\varepsilon_{\min}(\theta, \xi; \bar{\theta})$ satisfies the quadratic surrogation gap condition~\ref{enum:algo:P5:quadratic surrogation gap} for any $\varepsilon \geq 0$ and $\xi \in \Xi$;
    \item\label{assum:B6:generalized quadratic growth} 
    generalized quadratic growth: given $\varepsilon>0$ and $\rho > 0$, there exist positive constants $\delta_q$, $\lambda_q$, and $c_q$ such that 
    for any $\theta \in \Theta$ satisfying $\distSd{\theta} \leq \delta_q$ and $\max_{\theta' \in \cP_{\min}^{\varepsilon, \rho}(\theta)}\norm{\theta' - \theta} \leq \lambda_q$,
    there exists $\theta^d \in \CL{\cS^d}$
    such that $\bar{F}(t; \theta^d, \theta)\coloneqq F(\theta^d + t(\theta - \theta^d))$ is convex on $[0,1]$ and
    \[
        \bar{F}(t; \theta^d, \theta) \geq \bar{F}(0; \theta^d, \theta) + \frac{c_q}{2} t^2\norm{\theta - \theta^d}^2, \quad \forall\,t \in [0,1].
    \]
\end{enumerate}

Assumption~\ref{assum:B5:surrogation quadratic gap} has appeared in Proposition~\ref{prop:EC:connecting r-d and r:liu2022}, and we state it here officially.
Condition~\ref{enum:algo:P5:quadratic surrogation gap} originating from \citet{liu2022solving} is a technical condition for the nonasymptotic convergence analysis of the ESMM algorithm. Note that Condition~\ref{enum:algo:P5:quadratic surrogation gap} is verified for~\ref{enum:algo:phi1}-\ref{enum:algo:phi3} in Online Appendix~\ref{EC:sec:surrogation for sec3}. 
Assumption~\ref{assum:B6:generalized quadratic growth} is generalized from~\citet{drusvyatskiy2018error}. 
The condition holds, for example, when the outer function $\varphi$ is strongly convex and $\Theta$ is a polytope, which is verified in Online Appendix~\ref{EC:sec:verification}.
It can be expected to hold in more general cases since it is a local property only requiring the convexity and quadratic growth of $F$ near $\cS^d$.

With those assumptions, the nonasymptotic convergence result of ESMM for general cases is stated as follows.
\begin{theorem}[Nonasymptotic convergence of ESMM for general cases]\label{thm:EC:nonasymptotic convergence for general cases}
    Under Assumptions~\ref{assum:B1:Touching}-\ref{assum:B6:generalized quadratic growth}, for problem~\eqref{eq:algo:formulation:finite sum primal problem} under the same settings of Algorithm~\ref{ALGO:unconstrained case} specified in Theorem~\ref{thm:convergence:convergence results},
    $\EE\,[\distSd{\widetilde \theta^{\, T}}] = \widetilde{\cO}\Big(T^{-\frac{1}{2}\min\left\{\frac{\alpha}{2},1\right\}}\Big)$.
\end{theorem}

Following the error-bound analysis developed in Section~\ref{EC:subsec:analytical details on convergence}, the convergence result can be established with the following proposition.
\begin{proposition}\label{prop:EC:limiting stationarity and error bounds plus}
    Under Assumptions~\ref{assum:B1:Touching}-\ref{assum:B6:generalized quadratic growth}, for $\varepsilon > 0$ and $\rho > 0$, there exists a constant $\eta_r > 0$ such that $\distSd{\theta} \leq \eta_r r^{\varepsilon, \rho}(\theta)$ for any $\theta \in \Theta$.
\end{proposition}

Similar to the proof of the second statement of Proposition~\ref{prop:EC:limiting stationarity and error bounds}, to prove the above error bound, we require the following proposition, which is a generalization of Proposition~\ref{prop:proof:error bounds by r-d}.
\begin{proposition}\label{prop:proof:error bounds by r-d plus}
    Under Assumptions~\ref{assum:B1:Touching}-\ref{assum:B6:generalized quadratic growth}, for problem~\eqref{eq:algo:formulation:finite sum primal problem}with $\varepsilon>0$, $\rho > 0$, and the corresponding constants $\delta_q$, $\lambda_q$, and $c_q$ in Assumption~\ref{assum:B6:generalized quadratic growth},
    \[
        \distSd{\theta} \leq \frac{2}{c_q} r_d(\theta),\quad \text{$\forall\, \theta\in \Theta$ satisfying $\distSd{\theta} \leq \delta_q$ and $r^{\varepsilon, \rho}(\theta) \leq \lambda_q$.}
    \]
\end{proposition}
\proof{Proof of Proposition~\ref{prop:proof:error bounds by r-d plus}}
    The inequality holds trivially if $\theta \in \CL{\cS^d}$.
    By Assumption~\ref{assum:B6:generalized quadratic growth}, 
    for any $\theta \in \Theta\backslash \CL{\cS^d}$ satisfying $r^{\varepsilon, \rho}(\theta) \leq \lambda_q$ and $\distSd{\theta} \leq \delta_q$,
    there exists $\theta^d \in \CL{\cS^d}$ such that $\bar{F}(\cdot; \theta^d, \theta)$ is convex and $\bar{F}(t; \theta^d, \theta) \geq \bar{F}(0; \theta^d, \theta) + \frac{c_q}{2}t^2\norm{\theta - \theta^d}^2$ for all $t\in [0,1].$
    According to \citet[Theorem 3.3]{ECdrusvyatskiy2018error}, it follows that
    \[
        \norm{\theta - \theta^d}
        \leq \frac{2}{c_q} \mathrm{dist}(0, \partial \bar{F}(1; \theta^d, \theta)) \leq \frac{2}{c_q} \left[- \bar{F}(\cdot; \theta^d, \theta)'(1; -1)\right]  =    \frac{2}{c_q} \left[- F'\left(\theta; \frac{\theta^d - \theta}{\norm{\theta^d - \theta}}\right)\right] 
        \leq \frac{2}{c_q}r_d(\theta).
    \] 
    Thus, we obtain $\mathrm{dist}(\theta, \cS^d) \leq \norm{\theta - \theta^d} \leq \frac{2}{c_q} r_d(\theta)$, which completes the proof.\qedsymbol
\endproof

\proof{Proof of Proposition~\ref{prop:EC:limiting stationarity and error bounds plus}.}
    We first show that given $\varepsilon > \varepsilon' \geq 0$ and $\rho' > \rho > 0$, 
    for any $\lambda'_d > 0$ there exists $\lambda_d > 0$ such that 
    \begin{align}
    \label{eq:prop:EC4_residual_neighbor}
    \Big\{\theta' \in \Theta: r^{\varepsilon, \rho}(\theta) \leq \lambda_d, \ \lVert\theta' - \theta\rVert \leq 2\lambda_d,\  \theta \in \Theta\Big\}
    \subseteq  \left\{\theta \in \Theta : r^{\varepsilon', \rho'}(\theta) \leq \lambda_d' \right\}.
    \end{align}
    By contradiction, suppose that there exists $\bar{\lambda}'_d > 0$ such that for any $\lambda > 0$, 
    there exist points $\theta, \bar{\theta} \in \Theta$ 
    satisfying $r^{\varepsilon, \rho}(\theta) \leq \lambda$, $\norm{\theta - \bar{\theta}} \leq 2\lambda$, and $r^{\varepsilon', \rho'}(\bar{\theta}) > \bar{\lambda}'_d$.
    Then there exist two sequences $\{\theta_k\}_k$ and $\{\bar{\theta}_k\}_k$ in $\Theta$ with $\inf_{k \in \mathbb{N}}r^{\varepsilon', \rho'}(\bar{\theta}_k) \geq \bar{\lambda}'_d$ such that both $r^{\varepsilon, \rho}(\theta_k)$ and $\norm{\theta_k - \bar{\theta}_k}$ converge to 0.
    Since $\Theta$ is compact, by passing to a subsequence if necessary, we can assume that both $\{\theta_k\}_k$ and $\{\bar{\theta}_k\}_k$ converge to $\theta_\infty \in \Theta$. 
    By the first statement of the proposition, we have $r^{\varepsilon', \rho'}(\theta_\infty) = 0$.
    Let %
    $\bI_k \in {\arg\min}_{\bI \in \cI^{\varepsilon'}(\bar{\theta}_k)}\left\{
        \widehat{F}(\widehat{P}^{\rho'}_{\bI}(\bar{\theta}_k); \bar{\theta}_k, \bI) + \frac{{\rho'}}{2}\, \big\lVert \,{\widehat{P}^{\rho'}_{\bI}(\bar{\theta}_k) - \bar{\theta}_k}\, \big\rVert ^2\right\}$ such that $r^{\varepsilon', \rho'}(\bar{\theta}_k) = \, \big\lVert \,{\widehat{P}^{\rho'}_{\bI_k}(\bar{\theta}_k) - \bar{\theta}_k}\, \big\rVert  > \bar{\lambda}'_d$.  
    Since 
    the total number of index combination is finite, by passing to a subsequence if necessary,
    we can assume that $\bI_k \equiv \widehat{\bI}$ for any $k \in \mathbb{N}$.
    By continuity of $\widehat{P}^{\rho'}_{\widehat{\bI}}(\cdot)$ at $\theta_\infty$, we derive
    $\, \big\lVert \,{\widehat{P}^{\rho'}_{\widehat{\bI}}(\theta_\infty) - \theta_\infty}\, \big\rVert \geq \bar{\lambda}'_d$ and $ \widehat{F}(\widehat{P}^{\rho'}_{\widehat{\bI}}(\theta_\infty); \theta_\infty, \widehat{\bI}) + \frac{\rho'}{2}\, \big\lVert \,{\widehat{P}^{\rho'}_{\widehat{\bI}}(\theta_\infty) - \theta_\infty}\, \big\rVert^2 \leq F(\theta_\infty)$, implying $r^{\varepsilon', \rho'}(\theta_\infty) > 0$, a contradiction.

    Given $\varepsilon' \in (0, \varepsilon)$ and $\rho' \in (\rho, +\infty)$,
    let $\delta_q', \lambda_q'$, and $c_q'$ be the corresponding constants in Assumption~\ref{assum:B6:generalized quadratic growth}. 
    There exists $\lambda_q > 0$ such that~\eqref{eq:prop:EC4_residual_neighbor} holds for $\lambda_q'$. %
    For $\theta \in \Theta$ satisfying $\distSd{\theta} \leq \frac{\delta_q'}{3}$ and $r^{\varepsilon, \rho}(\theta) \leq \lambda_0 \coloneqq \min\left\{ \lambda_q, \frac{\delta_q'}{3} \right\}$, by Proposition~\ref{prop:EC:connecting r-d and r:liu2022} let $\theta'$ be the point satisfying $\norm{\theta' - {\theta}} \leq 2\,r^{\varepsilon, \rho}({\theta}) \leq 2 \lambda_q$
    and $r_d(\theta') \leq \left(2\rho + \frac{5}{2}\kappa\right)\,r^{\varepsilon, \rho}({\theta})$. Then we have $r^{\varepsilon', \rho'}(\theta') \leq \lambda_q'$ from~\eqref{eq:prop:EC4_residual_neighbor} and
    \[
        \mathrm{dist}(\theta', \cS^d) 
        \leq \distSd{\theta} + \norm{\theta - \theta'}
        \leq \frac{\delta_q'}{3} + 2 r^{\varepsilon, \rho}(\theta) \leq \delta_q'.
    \]
    It follows from Proposition~\ref{prop:proof:error bounds by r-d} that $ \distSd{\theta'} \leq \frac{2}{c_q'} r_d(\theta')$. 
    We further obtain that
    \[
        \distSd{\theta} 
        \leq \mathrm{dist}(\theta', \cS^d) + \norm{\theta - \theta'} 
        \leq \frac{2}{c_q'} r_d(\theta') + 2 r^{\varepsilon, \rho}(\theta)
        \leq \left(\frac{4\rho + 5\kappa}{c_q'} + 2\right) r^{\varepsilon, \rho}(\theta).
    \]
    Next, for $\theta \in \Theta$ satisfying $\distSd{\theta} > \frac{\delta_q'}{3}$, by the first statement of the proposition, there exists $\lambda_1 > 0$ such that $r^{\varepsilon, \rho}(\theta) \geq \lambda_1$. 
    With $\Theta$ being a compact set, let $\max_{\theta\in\Theta}\norm{\theta}_\infty = \mu$.
    Then we obtain $ \distSd{\theta} \leq 2\sqrt{q}\mu \leq \frac{2\sqrt{q}\mu}{\lambda_1}r^{\varepsilon, \rho}(\theta)$.
    Finally, for $\theta \in \Theta$ satisfying $r^{\varepsilon, \rho}(\theta) > \lambda_0$, 
    we simply have $\distSd{\theta} \leq 2\sqrt{q}\mu < \frac{2 \sqrt{q}\mu}{\lambda_0} r^{\varepsilon, \rho}(\theta)$.
    Therefore, we obtain $\distSd{\theta} \leq \eta_r {r}^{\varepsilon, \rho}(\theta)$ for any $\theta \in \Theta$ with $\eta_r
    = \max\left\{\frac{4\rho + 5\kappa}{c_q'} + 2, \frac{2 \sqrt{q}\mu}{\lambda_1}, \frac{2 \sqrt{q}\mu}
    {\min\{\lambda_q, \delta_q'/3\}}\right\}$, which completes the proof.\qedsymbol
\endproof

\subsection{Verification of Assumption~\ref{assum:B6:generalized quadratic growth}}
\label{EC:sec:verification}

We verify that Assumption~\ref{assum:B6:generalized quadratic growth} holds when the outer function $\varphi$ is strongly convex and $\Theta$ is a polytope by proving a stronger result in Claim~\ref{claim:EC:generalized QG for strongly convex varphi} with the following lemma.
\begin{lemma}[\citet{ECcui2018composite}, Proposition 2 (ii)]\label{lemma:EC:cui-strongconvexity-localminimizer}
    Let $\Theta \subseteq \RR^{q}$ be a closed convex set. Let $f:\RR^{q} \rightarrow \RR^d$ be piecewise affine and $\varphi : \RR^d \rightarrow \RR$ be convex. It holds that every d-stationary point of the composite function $\varphi \circ f$ on $\Theta$ is a local minimizer.
\end{lemma}
Recall the notations $\cP_{\min}^{\varepsilon, \rho}(\theta) \coloneqq 
    \mathop{\arg\min}_{\theta' \in \Theta} \left\{
        \HatFmin(\theta'; \theta) + \frac{\rho}{2}\norm{\theta' - \theta}^2
    \right\}$ and $\bar{F}(t;\theta_1, \theta_2) \coloneqq F(\theta_1 + t(\theta_2 - \theta_1))$ for $\theta_1, \theta_2 \in \Theta$. In the analysis below, let $\partial(S)$ denote the boundary of $\Xi$, $\INT{S}$ denote the interior of $\Xi$, and $\RI{S}$ denote the relative interior of $\Xi$.
    \begin{claim}\label{claim:EC:generalized QG for strongly convex varphi}
        Let  $\Theta \subseteq \RR^q$ be a polytope %
        and $\varphi(\cdot, \xi)$ be a strongly convex function with modulus $c_\varphi$ for any $\xi \in \Xi$.
        Given $\varepsilon > 0$ and $\rho > 0$, there exist $\lambda_q > 0$ and $c_q > 0$ such that for any $\theta \in \Theta$ satisfying $r^{\varepsilon, \rho}(\theta) \coloneqq \underset{\theta' \in \cP_{\min}^{\varepsilon, \rho}(\theta)}{\max} \, \norm{\theta - \theta'}  \leq \lambda_q$, there exists $\theta^d \in \CL{\cS^d}$ such that $\bar{F}(\cdot; \theta^d, \theta)$ is convex on $[0,1]$ and
        \[
            \bar{F}(t; \theta^d, \theta) \geq \bar{F}(0; \theta^d, \theta) + \frac{c_q}{2} t^2\norm{\theta - \theta^d}^2, \quad \forall\,t \in [0,1].
        \]
    \end{claim}
    \proof{Proof of Claim~\ref{claim:EC:generalized QG for strongly convex varphi}.}\label{proof:EC:claim:GQG condition for strongly cvx case}
    \newcommand{\SdcC}{\cS^d\bigcap\cC}
    \newcommand{\clSdcC}{\CL{\cS^d}\bigcap\cC}
    \newcommand{\StarC}{\cS^*_\cC}
    \newcommand{\NAcC}[1]{[\{#1\} + \cN_A]\bigcap \cC}
    \newcommand{\fC}{\mathfrak{C}}
    \newcommand{\fCC}{\widehat{\mathfrak{C}}}
    \newcommand{\cCC}{\widehat{\cC}}
    
    For $\theta \in \CL{\cS^d}$, the claim holds clearly with $\theta^d = \theta$.
    For $\theta \in \Theta\backslash\CL{\cS^d}$, it suffices to consider the case with $n=1$,
    for which we denote  $F(\theta) = \varphi(f(\theta))$, where $f:\RR^{q} \rightarrow \RR^d$ is piecewise linear and $\varphi$ is strongly convex. 
    The polytope $\Theta$ can be partitioned into a finite number of $q$-polytopes, denoted by $\fC \coloneqq \{\cC_i\}_{i=1}^N$, such that $f$ is linear on each of them and the intersection of any two polytopes in $\fC$ is their common face.
    Given $\cC \in \fC$, $f$ restricted to $\cC$ can be represented as $f(\theta) = A \theta$ for some $A \in \RR^{d \times q}$. %
    Since the strongly convex function $\varphi$ has a unique minimizer on $\{A\theta\in \RR^d \mid \theta \in \cC\}$, the set of minimizers of $F$ on $\cC$ denoted by $\StarC$ can be represented as $\StarC = \NAcC{\theta'}$ for any $\theta' \in \StarC$, where $\cN_A$ denotes the null space of $A$.
    
    We first show that $\StarC = \clSdcC$ if $\cS^d \bigcap \cC \neq \emptyset$. %
    By Lemma~\ref{lemma:EC:cui-strongconvexity-localminimizer}, 
    $\clSdcC\subseteq \StarC$. %
    We show the reverse inclusion by considering the two cases $\cS^d \bigcap \INT{\cC} \neq \emptyset$ and $\cS^d \bigcap \INT{\cC} = \emptyset$ respectively.
    If there exists $\widehat{\theta}_1^d \in \cS^d \bigcap \INT{\cC}$, then for any $ {\theta}^{\, \prime} \in [\{\widehat{\theta}_1^d\} + \cN_{A}]\bigcap \INT{\cC}$, 
    \[F'(\theta'; \theta - {\theta}')  
        = \varphi'(A \widehat{\theta}_1^d + A ({\theta'} - \widehat{\theta}_1^d); A(\theta - \widehat{\theta}_1^d) + A(\widehat{\theta}_1^d - \theta'))
        = F'(\widehat{\theta}_1^d; \theta - \widehat{\theta}_1^d) \geq 0,\quad\forall\,\theta \in \Theta,
    \]
    which indicates $\StarC = \NAcC{\widehat{\theta}_1^d} \subseteq \clSdcC$.
    If $\cS^d \bigcap \INT{\cC} = \emptyset$, we must have $\StarC \subseteq \partial(\cC)$. 
    Let $\widehat{\theta}_2^d \in \cS^d \bigcap {\cC} \subseteq \StarC$.
    Since $\StarC = \NAcC{\widehat{\theta}_2^d}$ is a polyhedral set, it is contained in a face of~$\cC$.
    Therefore, we have either $\RI\StarC \subseteq \cC'$ or $\RI\StarC \bigcap \cC' = \emptyset$ for any $\cC' \in \fC$.
    Then for any $\theta'\in \RI\StarC$, by restricting to the common polytope and noticing $ \widehat\theta_2^d- \theta' \in \cN_A$, we obtain
    $F'({\theta}^{\prime}; \theta - {\theta}^{\prime}) 
        = F'(\widehat{\theta}_2^d; \theta - \widehat{\theta}_2^d) \geq 0$ for any $\theta$ sufficiently close to $\theta'$.
    This implies $\RI{\StarC} \subseteq \SdcC$ and further $\StarC \subseteq \clSdcC$.
    Combining the above two cases, we obtain that  $\StarC %
    = \clSdcC$ if $\cS^d \bigcap \cC \neq \emptyset$. %
    
    We next show that for any $\cC \in \fC$ with $\SdcC = \emptyset$,  given $\varepsilon > 0$ and $\rho > 0$, there exists $\widehat{\lambda}_q > 0$ such that $r^{\varepsilon, \rho}(\theta) > \widehat{\lambda}_q$ for any $\theta \in \cC$.
    For such $\cC \in \fC$, %
    we have $\StarC \subseteq \partial{(\cC)}$ and $\StarC \bigcap \cS^d = \emptyset$.
    Since $\cS^d \bigcap \INT\cC = \emptyset$,
    by arguments above there exists a polytope $\cC' \in \fC$ such that $\StarC \subseteq \cC'$ and
    $\min_{\theta \in \StarC}\dist{\theta}{\cS^*_{\cC'}} > 0$. 
    Let $\bI_{\cC'}$ %
    denote the index combination of $f$ corresponding to $\cC'$. %
    Due to the Lipschitz continuity of $F$ on $\Theta$, for any $\varepsilon > 0$ there exist $\delta > 0$ and $\Delta > 0$ such that any $\theta \in \cC$ with $\dist{\theta}{\StarC}\leq \delta$ satisfies $\bI_{\cC'} \in \cI^\varepsilon(\theta)$ and $F(\theta) > \min_{\theta' \in \cC'}F(\theta') + \Delta$;
    consequently, for any given $\rho > 0$, there exists $\lambda_1 > 0$ such that $r^{\varepsilon, \rho}(\theta) > \lambda_1$. 
    On the other hand, for any point $\theta \in \cC$ satisfying $\dist{\theta}{\cS^*_{\cC}} > \delta$, by the first statement of Proposition~\ref{prop:EC:limiting stationarity and error bounds}, there exists $\lambda_2 > 0$ such that $r^{\varepsilon, \rho}(\theta) > \lambda_2$.
    Combining the above two cases, we obtain the desired statement.  %
    It yields that there exists $\lambda_q > 0$ such that for any $\theta \in \Theta$ satisfying $r^{\varepsilon, \rho}(\theta) \leq {\lambda}_q$, we have $\SdcC \neq \emptyset$ and $\clSdcC = \StarC$ for any $\mathcal C \in \fC_{\theta}\coloneqq \{\cC' \in \fC \mid \theta \in \cC'\}$. 
    
    Now, consider a point $\theta \in \Theta \backslash \CL{\cS^d}$ satisfying $r^{\varepsilon, \rho}(\theta) \leq {\lambda}_q$.  
    Given a polytope $\cC \in \fC_\theta$ denoted by $\cC \coloneqq \{\theta \in \RR^{q}\mid B \theta \leq b\}$, since $\clSdcC = \StarC$,
    we also have $\theta \in \cC\backslash\StarC$. 
    With $f(\theta) \coloneqq A \theta$ over $\cC$ for some matrix $A$ and $a^* \coloneqq A\theta^*$ for $\theta^* \in \StarC$, we obtain $\StarC = \{\theta \in \RR^{q}\mid B \theta \leq b, A\theta \leq a^*, -A\theta \leq -a^*\}$.
    By compactness of $\cC$ and Hoffman's lemma \citep{hoffman1952approximate}, there exists $c > 0$ depending on $\cC$ such that for any $\theta \in \cC$, there exists $\theta^d \in \StarC = \clSdcC$ satisfying
    \[
        \norm{\theta - \theta^d} = \dist{\theta}{\StarC} \leq c \, \big(  \norm{(B \theta - b)_+} + \norm{A\theta - a^*}  \big) = c \norm{A(\theta - \theta^d)}.
    \]
    Clearly, $\bar{F}(\cdot; \theta^d, \theta)$ is convex on $[0,1]$.
    Moreover, since $\varphi$ is strongly convex, for any $t \in [0,1]$,
    \begin{align*}
       \bar{F}(t; \theta^d, \theta) =  \varphi(f(\theta^d + t(\theta - \theta^d))) 
        &\geq \varphi(f(\theta^d)) + \frac{c_\varphi}{2}\norm{f(\theta^d + t(\theta - \theta^d)) - f(\theta^d)}^2\\
        &= \varphi(A \theta^d) + \frac{c_\varphi}{2} t^2\norm{A(\theta - \theta^d)}^2  \\
        &\geq \bar{F}(0; \theta^d, \theta) + \frac{c_\varphi}{2c^2} t^2\norm{\theta - \theta^d}^2.
    \end{align*}
    By the finiteness of $\fC$, %
    we thus prove the desired claim.\qedsymbol
    \endproof

\section{Constrained SP with Covariate Information}\label{section:Extension:constrained}
In this section, we study problem~\eqref{eq:intro:contextalSP} with deterministic constraints $\psi_j(z)\leq 0$ for $j \in [J]$.
We first provide the formulation of the PADR-based ERM method. Then, we establish the asymptotic consistency of the constrained PADR-based ERM model in Section~\ref{EC:subsec:asymptotic consistency of constrained ERM with PADR} and develop  ESMM algorithm with exact penalization in Section~\ref{EC:subsec:extension:convergence-to-b-stat}. We present the numerical results on the Newsvendor problem with convex and nonconvex constraints in Section~\ref{subsec:exp:cvx-constrained-nv}.

The optimization problem for learning the optimal decision rule can be formulated as
\begin{miniinf}
    {\substack{f\in \mathcal F}}{\EE_{X, Y}\left[\varphi(f(X), Y)\right]}{}{}
    \label{eq:model:constrained contextualSP-full E}
    \addConstraint{\psi_j(f(x))}{\leq 0,\quad}{\text{a.e.}\ x \in \cX, \  \forall\, j \in [J].}
\end{miniinf}
Constraints bring up the feasibility issue for data-driven decision rule methods. We show that in the following three special cases,  probabilistic constraints can be equivalently converted into deterministic ones or embedded into the structure of PADR, with which we can construct the PADR-based ERM problem while retaining the feasibility of the obtained PADR:
\begin{enumerate}[label = (C\arabic*), series=C]
    \item\label{enum:C1:finite scenarios} Finite scenarios: when the support set $\cX$ of the feature vector $x$ contains a finite number of known scenarios, feasibility can be guaranteed by imposing constraints on the decision rule under all scenarios.
    \item\label{enum:C2:box constraints} Box constraints: when deterministic constraints lead to a box-constrained set $\{z \in \mathbb R^d: \norm{z}_\infty \leq M\}$, a truncated PADR class $\barHK \coloneqq \{\bar{f}\mid \bar{f} = \max\{\min\{f, M\}, -M\},\,f\in \HK\}$ can be used, where $\HK$ represents the PA hypothesis class, and the max/min are element-wise operators. This ensures feasibility of the ERM solution.
\end{enumerate}
Following the excess risk decomposition in Section~\ref{subsec:nonasymptotic consistency of ERM with PADR}, it is straightforward to see that the PADR-based ERM model for the three special cases \ref{enum:C1:finite scenarios} and \ref{enum:C2:box constraints} is nonasymptotically consistent with the original SP  problem~\eqref{eq:model:constrained contextualSP-full E}. {For Case~\ref{enum:C1:finite scenarios}, problems with finitely many scenarios can be addressed using the exact penalization formulation~\eqref{eq:algo:penalized problem} developed in Section~\ref{EC:subsec:extension:convergence-to-b-stat}, and the feasibility of the solution is guaranteed by Lemma~\ref{lemma:cui2021modern}.}

For the case with general constraints, it is difficult to convert them into deterministic ones, so we impose constraints only on finite scenarios in $\Xi$ and we formulate a restricted PADR-based ERM problem with an additional constant $\gamma >0$ as below,
\begin{mini}
    {\substack{f\in \HK}}{\frac{1}{n}\sum_{s=1}^n \cost{f(x^s)}{x^s}{y^s}}{}{}
    \addConstraint{\psi_j(f(x^s)) + \gamma}{\leq 0,\quad}{\forall\,s\in [n],\ \forall\,j\in[J].}
    \label{eq:model:constrained ERM-H}
\end{mini}
Let ${\widetilde{\cH}}^K_{\mu, \gamma} \coloneqq  \HK \bigcap \left\{f\in \cF\mid \psi_j(f(x^s)) + \gamma \leq 0, \, s \in [n], \, j \in [J] \right\}$ , and let $f^*_{\Xi,\widetilde{\cH}^K_{\mu, \gamma}}$ 
denote an optimal PADR to the constrained ERM problem~\eqref{eq:model:constrained ERM-H}. 
The idea of imposing the restricted constraints on finite scenarios follows the scenario approximation method in \citet{calafiore2005uncertain}. This ensures that the obtained PADR $f^*_{\Xi,\widetilde{\cH}^K_{\mu, \gamma}}$ at a new scenario is feasible to the original constraint set with high probability according to \citet[Theorem 10]{luedtke2008sample}.
{The result is provided in the following proposition.
\begin{proposition}\label{prop:feasibility prob}
Suppose that the support set $\cX \subseteq \RR^p$ of $X$ is a compact set with $\max_{x \in \cX} \norm{x}_\infty = \bar{X}$, and that $\psi_j:\RR^d \rightarrow \RR$ is $L_\psi$-Lipschitz continuous for every $j \in [J]$. Then, with probability at least $1-\delta$ over the draw of $\Xi$, we have:
\begin{align*}
\mathbb{P}_X\left[ \max_{j \in [J]}\left\{\psi_j\left(f^*_{\Xi, \widetilde{\cH}_{\mu, \gamma}^K}(X)\right)\right\}> 0 \right] <
&\ \sqrt{\frac{4dK(p+1)}{n}\log\left(\left\lceil \frac{8 L_\psi}{\gamma} \sqrt{2K(p \bar{X}^2 + 1)}\cdot d(p+1)\mu \right\rceil\right)} \\
&+ \sqrt{\frac{\log(2n)}{n}} + \sqrt{\frac{2\log(1/\delta)}{n}}.
\end{align*}
\end{proposition}

\begin{proof}{Proof of Proposition~\ref{prop:feasibility prob}}
Denote $G(\theta, x) \coloneqq \max_{j\in [J]}\psi_j(f(\theta, x))$, where $\theta \in [-\mu, \mu]^{2dK(p+1)}$ is the parameter of PADR $f \in \HK$. Since each $\psi_j$ for $j\in[J]$ is $L_\psi$-Lipschitz continuous, we have that for any $x \in \cX$ and $\theta_1, \theta_2 \in [-\mu, \mu]^{2dK(p+1)}$,
\[
\begin{aligned}
\abs{G(\theta_1, x) - G(\theta_2, x)} 
&\leq L_\psi \norm{f(\theta_1, x) - f(\theta_2, x)} \\
&\leq L_\psi \sqrt{d} \norm{f(\theta_1, x) - f(\theta_2, x)}_\infty \\
&\leq 2 L_\psi \sqrt{d(p \bar{X}^2 + 1)(p+1)} \cdot \norm{\theta_1 - \theta_2}_\infty.
\end{aligned}
\]
Define the sets
\[
\begin{aligned}
&\Theta^\Xi_{\gamma} \coloneqq \left\{\theta \in [-\mu, \mu]^{2dK(p+1)} \mid G(\theta, x^s) \leq 0 ,\ s=1,\dots,n\right\},\\
&\Theta_\epsilon \coloneqq \left\{\theta \in [-\mu, \mu]^{2dK(p+1)} \mid \mathbb{P}_X(G(\theta, X) \leq 0) \geq 1 - \epsilon \right\}.
\end{aligned}
\]
By \citet[Theorem~10]{luedtke2008sample}, it follows that
\[
\mathbb{P}_\Xi \left(\Theta^\Xi_{\gamma} \subseteq \Theta_\epsilon\right) 
\geq 1 - \left\lceil \frac{2}{\epsilon}\right\rceil \left\lceil \frac{8 L_\psi}{\gamma} \sqrt{2K(p \bar{X}^2 + 1)}\cdot d(p+1)\mu \right\rceil^{2dK(p+1)} e^{ - n\epsilon^2/2}.
\]
With $\epsilon \in (0, 1]$, let
\[
M \coloneqq\left\lceil \frac{8 L_\psi}{\gamma} \sqrt{2K(p \bar{X}^2 + 1)}\cdot d(p+1)\mu \right\rceil^{2dK(p+1)},\
\delta \coloneqq \left\lceil \frac{2}{\epsilon}\right\rceil M e^{ - n\epsilon^2/2} \geq \frac{2}{\epsilon} M e^{ - n\epsilon^2/2}.
\]
Without loss of generality, we can assume that $L_\psi \geq \gamma$, which leads to $M > 2$ and $\log(\delta / (2M))< \log(1/4) < -1$.
Therefore, we have that
\[
\frac{\delta}{2M} \geq \frac{1}{\epsilon} e^{ - \frac{n}{2}\epsilon^2} 
> e^{(\log(\frac{\delta}{2M}) + 1)\frac{n}{2}\epsilon^2} \cdot e^{ - \frac{n}{2}\epsilon^2} = e^{\log(\frac{\delta}{2M})\frac{n}{2} \epsilon^2},
\]
implying $\epsilon > \sqrt{2 / n}$.
Furthermore, denote $\delta' \coloneqq M e^{-n\epsilon^2/2}$,
it follows that with probability at least $1-\delta$ over $\Xi$, the solution $f^*_{\Xi, \widetilde{\cH}_{\mu, \gamma}^K}$ satisfies:
\begin{align*}
&\mathbb{P}_{X}\left[
\max_{j\in [J]} \left\{\psi_j \left(f^*_{\Xi, \widetilde{\cH}_{\mu, \gamma}^K}(X)\right)\right\} > 0
\right]\\
&\leq
\sqrt{\frac{2}{n}\log\left(\left\lceil \frac{8 L_\psi}{\gamma} \sqrt{2K(p \bar{X}^2 + 1)}\cdot d(p+1)\mu\right\rceil^{2dK(p+1)}/\delta'\right)}\\
&\leq \sqrt{\frac{4dK(p+1)}{n}\log\left(\left\lceil \frac{8 L_\psi}{\gamma} \sqrt{2K(p \bar{X}^2 + 1)}\cdot d(p+1)\mu\right\rceil\right)} + \sqrt{\frac{2\log(1/\delta')}{n}}\\
&\leq \sqrt{\frac{4dK(p+1)}{n}\log\left(\left\lceil \frac{8 L_\psi}{\gamma} \sqrt{2K(p \bar{X}^2 + 1)}\cdot d(p+1)\mu \right\rceil\right)} + \sqrt{\frac{2\log(1/\delta)}{n}} + \sqrt{\frac{2\log(2/\epsilon)}{n}}\\
&< \sqrt{\frac{4dK(p+1)}{n}\log\left(\left\lceil \frac{8 L_\psi}{\gamma} \sqrt{2K(p \bar{X}^2 + 1)}\cdot d(p+1)\mu \right\rceil\right)} + \sqrt{\frac{2\log(1/\delta)}{n}} + \sqrt{\frac{\log(2n)}{n}},
\end{align*}
which completes the proof. 
\qedsymbol
\end{proof}

It can be observed that the out-of-sample feasibility of the PADR-based ERM solution exhibits the same order of complexity $\widetilde{\cO}\left(\sqrt{\frac{4dK(p+1)}{n}}\right)$ in $n$ and $K$ as the generalization bound in the unconstrained case. 
Similar to the excess risk result for unconstrained problems, with a larger piece number $K$,  model misspecification reduces but the likelihood of out-of-sample infeasibility simultaneously increases. 
}

We next show that the PADR-based ERM model \eqref{eq:model:constrained ERM-H} is asymptotically consistent with~\eqref{eq:model:constrained contextualSP-full E}.

\subsection{Asymptotic consistency of constrained PADR-based ERM}\label{EC:subsec:asymptotic consistency of constrained ERM with PADR}
To establish the consistency result of the PADR-based ERM model for the constrained SP problem~\eqref{eq:model:constrained contextualSP-full E}, we impose the following assumption on functions $\{\psi_j\}_{j \in [J]}$:
\begin{enumerate}[label = (A\arabic*), start=5, labelwidth=!, labelindent=0pt]
    \item\label{assum:A5:Lipschitz of psi} For every $j \in [J]$, $\psi_j:\RR^d \rightarrow \RR$ is $L_{\psi}$-Lipschitz continuous.
\end{enumerate}
By strengthening the constraints by a constant parameter $\gamma_0 > 0$, we obtain a restricted SP problem as follows:
\begin{miniinf}
    {\substack{f\in \mathcal F}}{R(f) \coloneqq \EE_{X, Y}\left[\varphi(f(X), Y)\right]}{}{}
    \label{eq:EC:model:Constrained-P-strengthened}
    \addConstraint{\psi_j(f(x)) + \gamma_0}{\leq 0,\quad}{\text{a.e.}\ x\in \mathcal X,\ \forall\,j\in[J],}
\end{miniinf}
for which we consider the following assumption, as an adaptation of Assumption~\ref{assum:A4:opt lipschitz solution f*} for constrained problems:
\begin{enumerate}[label=(A\arabic*$'$), start=4, ref=(A\arabic*$'$), labelwidth=!, labelindent=0pt]
    \item\label{assum:Ac4:constrained opt lipschitz solution f*} Given $\gamma_0>0$, there exists a Lipschitz continuous function $f^*_{\gamma_0}$ that is optimal to problem~\eqref{eq:EC:model:Constrained-P-strengthened}.
\end{enumerate}
By Assumptions~\ref{assum:A2:compactness of X} and~\ref{assum:Ac4:constrained opt lipschitz solution f*}, 
$f^*_{\gamma_0} \in \cF_{L_z', M_z'}$ for some $L_z'>0$ and $M_z' > 0$.
We first provide the excess risk bound of the restricted ERM model~\eqref{eq:model:constrained ERM-H} to the restricted SP problem \eqref{eq:EC:model:Constrained-P-strengthened}.
\begin{proposition}[Excess risk bound of the restricted constrained ERM]\label{prop:excess-risk-constr}
Suppose that Assumptions~\ref{assum:A1:iid dataset S}-\ref{assum:A3:Lipschitz of varphi} and~\ref{assum:A5:Lipschitz of psi} hold, %
    and that Assumption~\ref{assum:Ac4:constrained opt lipschitz solution f*} holds for some $\gamma_0 > \Gamma \coloneqq 2 d^{1/2} (p^{1/2} + 3) p^{1/2} L_\psi L_z' \bar{X} K^{-1/p}$. With $\gamma \in (0,\gamma_0 - \Gamma )$, 
     $\mu \geq \max\{\frac{1}{2}L_z' K^{1/p}, \frac{1}{4}pL_z' \bar{X}K^{1/p} + \frac{1}{2}M_z'\}$, and $M_{p,d}(\mu) \coloneqq C_\varphi + L_\varphi\left(2\sqrt{d}\mu(p\bar{X} + 1) + \norm{z_0}\right)$, the following holds with probability at least $1-\delta$ over the draw of $\Xi$,
    \begin{align*}
        R(f^*_{\Xi, \widetilde{\cH}^K_{\mu, \gamma}}) 
        \leq& R(f^*_{\gamma_0}) + 2 d^{1/2} (p^{1/2} + 3) p^{1/2} L_\varphi L_z' \bar{X} K^{-1/p}\\
        &+\frac{16\mu L_\varphi\sqrt{2dK(p+1)(p\bar{X}^2 + 1)}}{n}
        + 4M_{p,d}(\mu)\sqrt{\frac{2dK(p+1)\log n}{n}}
        + 2\sqrt{2}M_{p,d}(\mu) \sqrt{\frac{\log(2 / \delta)}{n}}.
    \end{align*}
\end{proposition}

\proof{Proof of Proposition~\ref{prop:excess-risk-constr}.}
    \label{proof:model:theorem 2:constrained excess risk}
    Denote by $f^*_{\widetilde{\cH}^K_{\mu, \gamma}}$ the optimal solution to the following problem
    \begin{mini}
        {\substack{f\in \HK}}{\EE_{X, Y}\left[\varphi(f(X), Y)\right]}{}{}
        \addConstraint{\psi_j(f(x^s)) + \gamma}{\leq 0,\quad}{\forall\, s\in [n],\ \forall\, j\in [J].}
    \end{mini}
    Similar to analysis for unconstrained problems, we deduce the following decomposition:
    \begin{align*}
            R(f^*_{\Xi,\widetilde{\cH}^K_{\mu, \gamma}}) - R(f^*_{\gamma_0})
            &= R(f^*_{\Xi, \widetilde{\cH}^K_{\mu, \gamma}}) - \RS(f^*_{\Xi, \widetilde{\cH}^K_{\mu, \gamma}})
            + \RS(f^*_{\Xi, \widetilde{\cH}^K_{\mu, \gamma}}) - R(f^*_{\widetilde{\cH}^K_{\mu, \gamma}})
            + R(f^*_{\widetilde{\cH}^K_{\mu, \gamma}}) - R(f^*_{\gamma_0})\\
            &\leq R(f^*_{\Xi, \widetilde{\cH}^K_{\mu, \gamma}}) - \RS(f^*_{\Xi, \widetilde{\cH}^K_{\mu, \gamma}})
            + \RS(f^*_{\widetilde{\cH}^K_{\mu, \gamma}}) - R(f^*_{\widetilde{\cH}^K_{\mu, \gamma}})
            + R(f^*_{\widetilde{\cH}^K_{\mu, \gamma}}) - R(f^*_{\gamma_0})\\
            &\leq \underbrace{2\max_{f \in \HK} \abs{R(f) - \RS(f)}}_{\cE_1} + \underbrace{R(f^*_{\widetilde{\cH}^K_{\mu, \gamma}}) - R(f^*_{\gamma_0})}_{\cE_2},
    \end{align*}
    where 
    the last inequality holds since $\widetilde{\cH}^K_{\mu, \gamma} \subseteq \HK$. 
    The first term $\cE_1$ can be directly bounded by Proposition~\ref{prop:model:Uniform generalization error bound}, that is,
    given $\delta > 0$, with probability at least $1-\delta$ over the draw of $\Xi$,
    \[
        \cE_1 
        \leq 
        \frac{16\mu L_\varphi\sqrt{2dK(p+1)(p\bar{X}^2 + 1)}}{n}
        + 4M_{p,d}(\mu)\sqrt{\frac{2dK(p+1)\log n}{n}}
        + 2\sqrt{2}M_{p,d}(\mu) \sqrt{\frac{\log(2 / \delta)}{n}}.
    \]
    For the second term $\cE_2$, we further decompose it by
    \[
        \cE_2 = \underbrace{R(f^*_{\widetilde{\cH}^K_{\mu, \gamma}}) - R(f_{app,\gamma_0})}_{\cE_{21}} + \underbrace{R(f_{app, \gamma_0}) - R(f^*_{\gamma_0})}_{\cE_{22}},
    \]
    where $f_{app, \gamma_0} \in \mathop{\arg\min}_{f \in \HK}\norm{f - f^*_{\gamma_0}}_\infty$.
    Then for almost every $x\in\cX$, 
    \begin{align*}
            {\psi_j}(f_{app, \gamma_0}(x)) + \gamma 
            &= {\psi_j}(f^*_{\gamma_0}(x)) + \left({\psi_j}(f_{app, \gamma_0}(x)) - {\psi_j}(f^*_{\gamma_0}(x))\right) + \gamma\\
            &\leq {\psi_j}(f^*_{\gamma_0}(x)) + L_\psi \sqrt{d} \norm{f_{app, \gamma_0} - f^*_{\gamma_0}}_\infty + \gamma\\
            &\leq {\psi_j}(f^*_{\gamma_0}(x)) + 2 d^{1/2} (p^{1/2} + 3) p^{1/2} L_\psi L_z' \bar{X} K^{-1/p} + \gamma\\
            &< {\psi_j}(f^*_{\gamma_0}(x)) + \gamma_0,
    \end{align*}
    where the second inequality holds by Proposition~\ref{prop:model:universal approximation of PADR}.
    This implies that $\psi_j(f_{app, \gamma_0}(x^s)) + \gamma \leq 0$ for any $s \in [n]$ and $j \in [J]$, i.e., $f_{app, \gamma_0} \in \widetilde{\cH}^K_{\mu, \gamma}$, with probability 1. 
    It follows that $\cE_{21} = \min_{f \in \widetilde{\cH}^K_{\mu, \gamma}} R(f) - R(f_{app, \gamma_0}) \leq 0$ with probability 1.
    Moreover, $\cE_{22}$ can be bounded as
    \[
        \cE_{22} \leq L_\varphi \sqrt{d} \norm{f_{app, \gamma_0} - f^*_{\gamma_0}}_{\infty} 
        \leq 2 d^{1/2} (p^{1/2} + 3) p^{1/2} L_\varphi L_z' \bar{X} K^{-1/p}.
    \]
    By combining bounds of $\cE_1$, $\cE_{21}$, and $\cE_{22}$, we derive the desired bound.\qedsymbol
\endproof
With the excess risk bound in Proposition~\ref{prop:excess-risk-constr}, and the probabilistic feasibility guarantee according to Proposition~\ref{prop:feasibility prob}, 
we obtain the asymptotic consistency of restricted PADR-based ERM model to the original constrained SP \eqref{eq:model:constrained contextualSP-full E}.
\begin{theorem}[Asymptotic consistency of constrained PADR-based ERM]\label{thm:asymptotic-consistency-constr}
    Suppose that Assumptions~\ref{assum:A1:iid dataset S}-\ref{assum:A3:Lipschitz of varphi} and~\ref{assum:A5:Lipschitz of psi} hold, and that $f^*_{\gamma} \in \cF_{L_z', M_z'}$ for any $\gamma \in [0, \bar{\gamma}]$ with some $\bar{\gamma} > 0$. Then for any $n > 0$, there exist $\delta_n >0$ and $\epsilon_n > 0$ such that
    $f^*_{\Xi,\widetilde{\cH}^K_{\mu, \gamma}}$ satisfies $\mathbb{P}_X\big\{\psi_j(f^*_{\Xi,\widetilde{\cH}^K_{\mu, \gamma}}(X))\leq 0,\,j\in [J]\big\} \geq 1 - \epsilon_n$ with probability $1-\delta_n$, where $\epsilon_n$ and $\delta_n$ converge to 0 as $n$ goes to infinity. Moreover, with $K = \cO(n^{p/(p+4)})$, $\mu = \cO(K^{1/p})$, and $\gamma = \cO(K^{-1/p})$, $R(f^*_{\Xi,\widetilde{\cH}^K_{\mu, \gamma}})$ converges to the optimal cost of~\eqref{eq:model:constrained contextualSP-full E} in probability, as $n$ goes to infinity.
\end{theorem}

\subsection{Convergence of ESMM to B-stationarity}\label{EC:subsec:extension:convergence-to-b-stat}
The constrained ERM problem~\eqref{eq:model:constrained ERM-H}can be reformulated as follows with $\widetilde{\xi}_n$ representing a random variable uniformly distributed on $\Xi$:
\begin{mini}
    {\substack{\theta \in \Theta}}
    {F(\theta) \coloneqq \EE\, \left[\varphi(f(\theta, \widetilde\xi_n), \widetilde\xi_n)\right]\phantom{+\frac{1}{n}\lambda\max++++}}{}{}
    \addConstraint{G_\gamma(\theta)\coloneqq \EE\, \left[\sum_{j=1}^J
    \max\{\psi_j(f(\theta, \widetilde\xi_n)) + \gamma, 0\}\right] }{\leq 0.}{}\label{eq:algo:Constrained-ERM-P-H-equiv}
\end{mini}
Suppose that functions $\varphi$ and $\{\psi_j\}_{j \in [J]}$ are Bouligand-differentiable (B-differentiable), i.e., they are locally Lipschitz continuous and directionally differentiable on $\RR^d$,
and that $\Theta$ is a compact convex set.
Since $G_\gamma$ is nonconvex, we are interested in computing Bouligand-stationary (B-stationary) points of~\eqref{eq:algo:Constrained-ERM-P-H-equiv}, the defining condition of which is an extension of d-stationarity to nonconvex-constrained optimization problems. With $\cI > 0$, we construct the following penalized problem
\begin{mini}
    {\substack{\theta \in \Theta}}
    {V(\theta; \lambda) \coloneqq F(\theta) + \lambda G_\gamma(\theta).}{}{}
    \label{eq:algo:penalized problem}
\end{mini} 
According to \citet[Proposition 9.2.2]{cui2021modern}, problem~\eqref{eq:model:constrained ERM-H} with some regularity conditions 
has the exact penalty property such that with a finite penalty parameter $\lambda$, any d-stationary point of~\eqref{eq:algo:penalized problem} is a B-stationary point of~\eqref{eq:algo:Constrained-ERM-P-H-equiv}. 
Since $\max\{\cdot,0\}$ is nondecreasing and convex, it is clear that the upper  surrogate functions of $G_\gamma(\cdot)$ can be constructed by plugging surrogates of $\{\psi_j(f(\cdot, \xi))\}_{j\in [J]}$ into $\max\{\cdot,0\}$, maintaining the corresponding surrogation properties.
Therefore, with the exact penalty property and Assumptions~\ref{assum:B1:Touching}-\ref{assum:B6:generalized quadratic growth} imposed in the context of problem~\eqref{eq:algo:penalized problem}, the ESMM algorithm is ready to solve the penalized problem~\eqref{eq:algo:penalized problem} with convergence to B-stationary points of~\eqref{eq:algo:Constrained-ERM-P-H-equiv}. According to \citep{ECcui2021modern}, for a constrained problem with a B-differentiable objective, a B-stationary point is defined as follows.
\begin{definition}
    We say that $\theta^b$ is a B-stationary point of the objective $F(\cdot)$ on a closed set $\bar{\Theta}\subseteq \Theta$ if $\theta^b \in \bar{\Theta}$ and $F'(\theta^b; v) \geq 0$ for any $v \in \cT(\theta^b; \bar{\Theta})$.
    where $\cT(\theta^b; \bar{\Theta})$ is the tangent cone of $\bar{\Theta}$ at $\theta^b$.
\end{definition}
The B-stationarity coincides with d-stationarity if $\bar{\Theta}$ is convex.
Suppose that functions $\varphi$ and $\{\psi_j\}_{j \in [J]}$ in problem~\eqref{eq:algo:Constrained-ERM-P-H-equiv} are B-differentiable.
With $\bar{\Theta} = \left\{\theta\in\Theta \mid G_\gamma(\theta)\leq 0\right\}$, we denote the set of B-stationary points of~\eqref{eq:algo:Constrained-ERM-P-H-equiv} by $\cS^b$.

With Assumptions~\ref{assum:B1:Touching}-\ref{assum:B6:generalized quadratic growth} imposed in the context of the penalized problem~\eqref{eq:algo:penalized problem},
we can apply the ESMM algorithm to solve~\eqref{eq:algo:penalized problem} and obtain the same convergence results in terms of d-stationarity as in Theorems~\ref{thm:convergence:convergence results}-\ref{thm:EC:nonasymptotic convergence for general cases}.
To further obtain the convergence to B-stationarity, we intend to connect d-stationarity of~\eqref{eq:algo:penalized problem} with B-stationarity of~\eqref{eq:algo:Constrained-ERM-P-H-equiv} by applying the following exact penalty result from \citet{ECcui2021modern}.
Note that by the finiteness of $\Xi$ and compactness of $\Theta$, for any $\xi \in \Xi$, $f(\cdot, \xi)$ is Lipschitz continuous with modulus $L_f\coloneqq 2\sqrt{\max_{s\in[n]}\norm{x^s}^2 + 1}$ and uniformly bounded with $\norm{f(\cdot, \xi)}_\infty \leq M_f\coloneqq L_f\max_{\theta\in\Theta}\norm{\theta}$, and $\varphi(\cdot, \xi)$ is Lipschitz continuous on $[-M_f, M_f]^d$ with uniform modulus denoted by $L'_\varphi$.

\begin{lemma}[\citet{ECcui2021modern}, Proposition 9.2.2]\label{lemma:cui2021modern}
    Let $\Theta$ be a compact convex set. 
    For $\gamma >0$, if there exists $l_c > 0$ such that
    \begin{equation}
        \label{eq:convergence:d-d condition}
        \sup_{\theta \in \Theta\setminus\bar{\Theta}}\left[
            \min_{v \in \cT(\theta; \Theta); \norm{v}\leq 1} G_\gamma'(\theta; v)
        \right] \leq - l_c,
    \end{equation}
    then with $\lambda > L'_\varphi L_f / l_c$, any d-stationary point of~\eqref{eq:algo:penalized problem} is feasible to~\eqref{eq:algo:Constrained-ERM-P-H-equiv} and thus is a B-stationary point of~\eqref{eq:algo:Constrained-ERM-P-H-equiv}.
\end{lemma}

Condition~\eqref{eq:convergence:d-d condition} requires that the maximum descent rates of $G_\gamma$ at points infeasible to~\eqref{eq:algo:Constrained-ERM-P-H-equiv} are bounded away from 0.
Thus by the Lipschitz continuity of $F$, at any infeasible point,
we can always find a descent direction in the penalized problem~\eqref{eq:algo:penalized problem} with a sufficiently large (but still finite) penalty parameter $\lambda$, ensuring that the point can not be a d-stationary point of~\eqref{eq:algo:penalized problem}. 
Nonetheless, it is not easy to verify the technical condition~\eqref{eq:convergence:d-d condition} when $G_\gamma$ is in a composite structure.
The following assumption on $\{\psi_j\}_{j \in [J]}$ provides a sufficient condition for~\eqref{eq:convergence:d-d condition}:
\begin{enumerate}[label=(B\arabic*), resume*=B, labelwidth=!, labelindent=0pt]
    \item\label{assum:B5:constraint condition for dd condition} Given $\gamma > 0$, there exist positive constants $c_{\psi}$ and $l_{\psi}$ such that for $z \in \RR^d$, if $\psi_{j_0}(z) + \gamma \geq 0$ for some $j_0 \in [J]$, then $\norm{z} \geq c_{\psi}$ and $\psi_{j_0}'\left(z; -\frac{z}{\norm{z}}\right) \leq -l_{\psi}$.
\end{enumerate}
Assumption~\ref{assum:B5:constraint condition for dd condition} is more intuitive than~\eqref{eq:convergence:d-d condition}.
For example, consider a decision-making problem with a constraint function $\psi(z) = \norm{z}_1 - C_0$, where $C_0 > 0$ represents the total (absolute) capacity. %
Then the constraint $\psi(z) + \gamma = \sum_{\iota=1}^d \abs{z_\iota} - (C_0 - \gamma) \leq 0$ with $\gamma < C_0$ satisfies~\ref{assum:B5:constraint condition for dd condition} with $c_{\psi} = (C_0 - \gamma) / \sqrt{d}$ and $l_{\psi} = 1$. 
Based on Lemma~\ref{lemma:cui2021modern}, we derive the following exact penalty result under Assumption~\ref{assum:B5:constraint condition for dd condition} by utilizing the piecewise linearity of the inner function. %
\begin{proposition}\label{prop:convergence:sufficient cond for dd condition}
    Suppose that Assumption~\ref{assum:B5:constraint condition for dd condition} holds, and that $\Theta\subseteq \RR^q$ is a compact convex set with $0 \in \Theta$ and $\max_{\theta \in \Theta}\norm{\theta}_\infty = \mu$.
    Then condition~\eqref{eq:convergence:d-d condition} holds with $l_c = \frac{c_{\psi}l_{\psi}}{n \sqrt{q}\mu}$. Moreover, with $\lambda > \bar\lambda\coloneqq \frac{n \sqrt{q}\mu L'_\varphi L_f}{c_\psi l_\psi}$, any d-stationary point of~\eqref{eq:algo:penalized problem} is a B-stationary point of problem~\eqref{eq:algo:Constrained-ERM-P-H-equiv}.
\end{proposition}

\proof{Proof of Proposition~\ref{prop:convergence:sufficient cond for dd condition}.}\label{proof:appendix:prop:convergence:sufficient cond for dd condition}
    By Lemma~\ref{lemma:cui2021modern}, it suffices to show that condition~\eqref{eq:convergence:d-d condition} holds with $l_c = \frac{c_{\psi} l_{\psi}}{n \sqrt{q}\mu}$.
    For any $\xi \in \Xi$ and $j \in [J]$, denote $G_\gamma(\theta, \xi; j) \coloneqq \max\left\{\psi_j(f(\theta, \xi))+\gamma, 0\right\}$.
    For a point $\bar{\theta} \in \Theta$ infeasible to the constraint 
    $G_\gamma(\theta) \coloneqq \EE\left[\sum_{j=1}^J G_\gamma(\theta, \widetilde\xi_n; j)\right] \leq 0$, 
    there exist $\bar{j} \in [J]$ and $\bar{\xi} \in \Xi$ such that $\psi_{\bar j}(f(\bar{\theta}, \bar\xi))+\gamma > 0$. 
    
    For any $\xi \in \Xi$, we denote $\bar{z} = f(\bar{\theta}, \xi)$ and $\bar{v} = - \frac{\bar{\theta}}{\norm{\bar{\theta}}}$ for simplicity.
    Since $f(\cdot, \xi)$ is piecewise linear,  
    there exists $A \in \RR^{d \times q}$ dependent on $\xi$ and $\bar{\theta}$ such that $f(\bar{\theta}, \xi) = A\bar{\theta}$ and $f(\cdot, \xi)'(\bar{\theta}; \bar{v}) = A\bar{v}$. Let
    \[ 
    \mathcal A_< \coloneqq \{j \in [J]: \psi_j(\bar z) + \gamma < 0\}, \ 
    \mathcal A_> \coloneqq \{j \in [J]: \psi_j(\bar z) + \gamma > 0\}, \text{ and } \mathcal A_= \coloneqq \{j \in [J]: \psi_j(\bar z) + \gamma = 0\}.
    \]
    For $j \in \mathcal A_<$, it is clear that $G_\gamma(\cdot, \xi; j)'(\bar{\theta}; \bar{v}) = 0$. 
    By Assumption~\ref{assum:B5:constraint condition for dd condition}, we have that 
    for $j \in \mathcal A_>$,
    \[
        G_\gamma(\cdot, \xi; j)'(\bar{\theta}; \bar{v}) 
        = \psi_j'(f(\bar{\theta}, \xi); f(\cdot, \xi)'(\bar{\theta}; \bar{v})) 
        = \psi_j'\left(\bar{z}; -\frac{\bar{z}}{\norm{\bar{\theta}}}\right)\leq -\frac{c_{\psi}l_{\psi}}{\norm{\bar{\theta}}}\leq - \frac{c_{\psi}l_{\psi}}{\sqrt{q}\mu},
    \]
    and that for $j \in \mathcal A_=$, 
    \[
        G_\gamma(\cdot, \xi; j)'(\bar{\theta}; \bar{v}) = \max \big\{\psi_j'(f(\bar{\theta}, \xi); f(\cdot, \xi)'(\bar{\theta}; \bar{v})) , \, 0 \big\} = \max\left\{\psi_j'\left(\bar{z}; -\frac{\bar{z}}{\norm{\bar{\theta}}}\right), 0\right\} = 0.
    \]
    Since $\bar{j} \in A_>$, we obtain that
    \[
        \min_{v \in \cT(\bar{\theta}; \Theta); \norm{v}\leq 1} G_\gamma'(\bar{\theta}; v) 
        \leq G_\gamma'(\bar{\theta}; \bar{v})
        \leq \frac{1}{n}G_\gamma(\cdot, \bar\xi; \bar{j})'(\bar{\theta}; \bar{v})
        \leq -\frac{c_{\psi} l_{\psi}}{n \sqrt{q}\mu},
    \]
    which completes the proof.\qedsymbol
\endproof

Finally, from Proposition~\ref{prop:convergence:sufficient cond for dd condition} and convergence results in Sections~\ref{subsec:the surrogation and the algorithm} and~\ref{EC:sec:nonasymptotic convergence for general cases}, we summarize the convergence properties of the ESMM with exact penalization in the theorem below, where Assumptions~\ref{assum:B1:Touching} to~\ref{assum:B5:surrogation quadratic gap} are adjusted to surrogate families of $F(\cdot, \xi)$ and $\psi_j(f(\cdot, \xi))$ for any $\xi \in \Xi$ and $j \in [J]$, and Assumption~\ref{assum:B6:generalized quadratic growth} is adjusted to $V(\cdot; \lambda)$. The proof is similar to those of Theorems~\ref{thm:convergence:convergence results}-\ref{thm:EC:nonasymptotic convergence for general cases}, and thus is omitted here.

\begin{theorem}[Convergence of ESMM with exact penalization]\label{thm:EC:convergence of ESMM:constrained}
    Under Assumptions~\ref{assum:B1:Touching} to~\ref{assum:B4:minSurr dd-consistency} and the same settings of Proposition~\ref{prop:convergence:sufficient cond for dd condition},
    for problem~\eqref{eq:algo:penalized problem} with $\lambda \geq \bar \lambda$, 
    let $\{\widetilde{\theta}^{\,t}\}_t$ be a sequence of points produced by Algorithm~\ref{ALGO:unconstrained case}
    with $\varepsilon\geq 0$, $\eta \geq 0$, and $N_{\nu} = \cO(\nu^\alpha)$ for some $\alpha > 0$.
    The following statements hold:
    \begin{enumerate}[label={(\alph*)}, labelwidth=!, labelindent=0pt]
        \item With $\varepsilon > 0$, any limit point of $\{\widetilde{\theta}^t\}_t$ is a composite $(\varepsilon', \rho)$-strongly d-stationary point of~\eqref{eq:algo:penalized problem} with probability 1 for any $\varepsilon' \in [0, \varepsilon)$ and $\rho \in (\eta, +\infty)$, and thus a B-stationary point of~\eqref{eq:algo:Constrained-ERM-P-H-equiv}.
        
        \item With $\varepsilon > 0$, under Assumptions~\ref{assum:B5:surrogation quadratic gap} and~\ref{assum:B6:generalized quadratic growth},
        $\EE\,[\mathrm{dist}(\widetilde{\theta}^T, \cS^b)] = \widetilde{\cO}(T^{-\frac{1}{2}\min\{\frac{\alpha}{2}, 1\}})$. In particular, if $\varphi$ and $\{\psi_j\}_{j \in [J]}$ are piecewise affine and $\Theta$ is a polytope, the same convergence rate holds without the two assumptions.

        \item With $\varepsilon = 0$, functions $\varphi$ and $\{\psi_j\}_{j \in [J]}$ being piecewise affine, and $\Theta$ being a polytope, 
        $\EE\,[\mathrm{dist}(\widetilde{\theta}^T, \cS^b)] = \widetilde{\cO}(T^{-\frac{1}{2}\min\{\frac{\alpha}{2}, 1\}})$;
        moreover, with $\eta > 0$, $\EE\,[\mathrm{dist}(\widetilde{\theta}^T, \cS^b)] = \cO(T^{-\frac{1}{2}\min\{\alpha, 1\}})$, which is up to a logarithmic factor when $\alpha = 1$.
    \end{enumerate}
\end{theorem}
 For the ESMM algorithm with exact penalization, the limit point of generated sequence is a composite strong d-stationarity to~\eqref{eq:algo:penalized problem} with probability 1 according to Theorem~\ref{thm:EC:convergence of ESMM:constrained}, and thus is  a global minimizer of~\eqref{eq:algo:penalized problem} when $\varepsilon$ is sufficiently large and $\eta$ is close to 0,  by the equivalence between composite strong d-stationarity and global optimality in Proposition~\ref{prop:EC:relationship-S* equiv Sder} with appropriate assumptions when needed.
Furthermore, such a global minimizer of the penalized problem~\eqref{eq:algo:penalized problem}  should be also a global minimizer of~\eqref{eq:algo:Constrained-ERM-P-H-equiv} since it is feasible according to Proposition~\ref{prop:convergence:sufficient cond for dd condition}. Under such settings, the asymptotic consistency of the PADR framework for the constrained case in Theorem~\ref{thm:asymptotic-consistency-constr} turns into asymptotic optimality.

\subsection{Experiments on constrained newsvendor problem}\label{subsec:exp:cvx-constrained-nv}
We compare PADR and other benchmarks on two-product constrained Newsvendor problems, assuming that the uncensored historical demand data is available. 
We denote the unit back-order cost, unit holding cost, and order decision by $c_{bi}$, $c_{hi}$, and $z_i$, respectively, for product $i=1,2$. The demand models for the two products are set as follows:
\[
    \overline{Y}_1(X) = 15 X_1 - 5 X_2 + 30,\quad \overline{Y}_2(X) = 15 X_1 + 5 X_2 + 30.
\]
{In this subsection, we implement the prescriptive methods, which ensure feasibility by solving the individual constrained problem for each sample.
Additionally, we implement NN-DR equipped with $\ell_1$-penalization and RKHS-DR methods following the penalized formulation proposed by \citet[Section 6.2]{bertsimas2022data}.}

\subsubsection{Newsvendor problem with a linear constraint.} 
\label{subsec:Newsvendor_linear_constraint}
We first consider an inventory capacity constraint for the two products, which is formulated as $z_1 + z_2 \leq C_0$.
In this subsection, we set $(c_{b1}, c_{h1}, c_{b2}, c_{h2}) = (8, 2, 2, 8)$, $C_0 = 60$, such that the probability that the constraint is active at the optimal decision is around 0.5.

Since DR-based methods (i.e., PADR, NN-DR, and RKHS-DR) may yield infeasible decisions, 
in Figure~\ref{fig:exp:C-CVX-dim} we report the test cost of the order decisions, projected onto the convex feasible region when necessary, along with the feasible frequency of their original decisions.
\begin{figure}[htbp]
    \FIGURE
    {{
      \subfloat[Basic results]{%
          \includegraphics[width=0.33\textwidth]{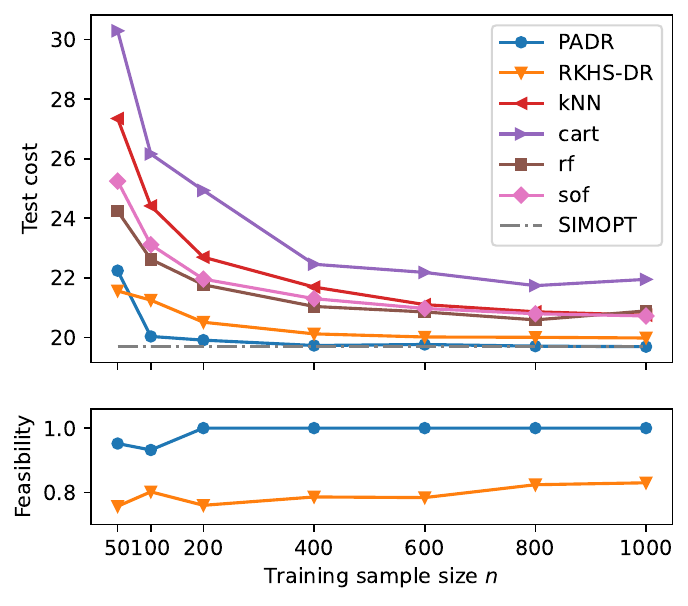}\label{subfig:constr:basic}
      }
      \subfloat[Sparse demand model]{%
          \includegraphics[width=0.33\textwidth]{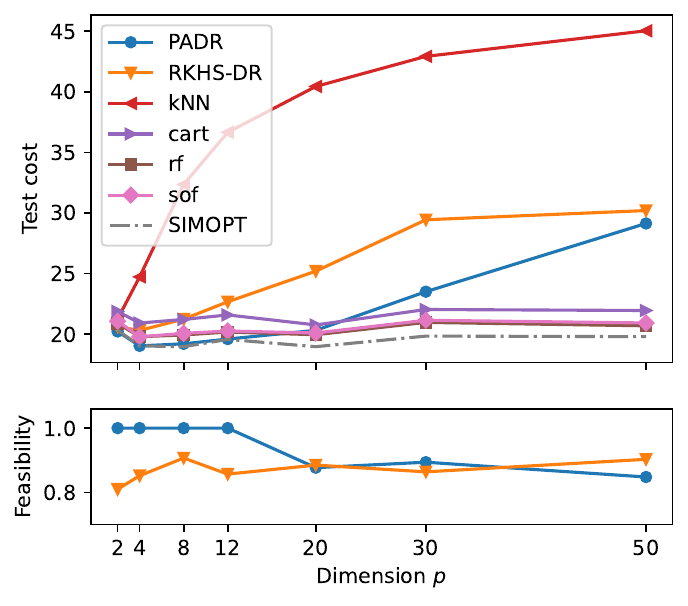}\label{subfig:constr:sparse}
      }
      \subfloat[Dense demand model]{%
          \includegraphics[width=0.33\textwidth]{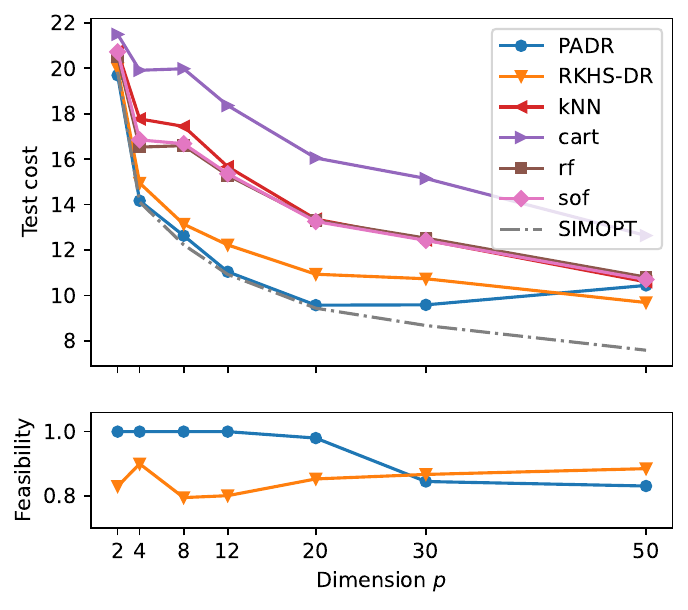}\label{subfig:constr:dense}
      }
    }}
    {Newsvendor problem with a linear constraint\label{fig:exp:C-CVX-dim}}
    {}
\end{figure}

{Figure~\ref{subfig:constr:basic} shows that with the basic demand model, PADR approaches SIMOPT and achieves 100\% feasibility when $n \geq 200$. In contrast, other benchmarks exhibit noticeable performance gaps relative to SIMOPT. Particularly, RKHS-DR achieves roughly 80\% feasibility and NN-DR exhibits unstable feasibility performance.
Figures~\ref{subfig:constr:sparse} and~\ref{subfig:constr:dense} present the results under varying feature dimensions. 
Consistent with the unconstrained setting, PADR outperforms other methods in the sparse model when $p < 20$ and in the dense model when $p < 50$.
Notably, by the utilization of $\ell_1$ penalization, PADR achieves 100\% feasibility when $p < 20$.
Moreover, in the dense demand model, PADR exhibits the clear edge over prescriptive methods by achieving better empirical performance on test costs, while it fails to be feasible at 10\% - 20\% probability.}

\subsubsection{Newsvendor problem with a nonconvex constraint.}
We next consider a capacity constraint $C(z_1) + C(z_2) \leq C_0$ with the concave capacity cost function~\eqref{eq:exp:nonconvex cost} for the two products. 
We set $(c_{b1}, c_{h1}, c_{b2}, c_{h2}) = (7, 7, 3, 3)$ and $C_0 = 50$ with the same probability for the constraint being active as the experiment in Section \ref{subsec:Newsvendor_linear_constraint}.
Since the projection of infeasible decisions to the nonconvex feasible region may be far away to the optimal solutions, we only report the average test cost of feasible decisions in Figure~\ref{fig:exp:UNC-NCVX-constr-dim}, with a careful hyperparameter selection to ensure high feasible frequencies.

We use the standard MM algorithm to solve the penalized ERM problem of RKHS-DR, referred to here as RKHS-DR(MM). 
For prescriptive methods and SIMOPT, we solve the two-dimensional nonconvex-constrained problems using grid search. 
The SOF method in \citet{kallus2022stochastic}is not suitable for this nonconvex-constrained problem, and it may require the study of the second-order perturbation analysis for stochastic optimization, which is beyond the scope of this paper.
For the other prescriptive methods, we solve two-dimensional nonconvex constrained problems by grid search. 

For the problem with the basic demand model, both PADR and RKHS-DR attain 100\% feasibility on the test set. PADR converges to SIMOPT with significantly fewer training samples than other benchmarks. 
Figures~\ref{subfig:ncvx-constr-sparse} and~\ref{subfig:ncvx-constr-dense} show the results across varying feature dimensions. PADR outperforms the other methods in both sparse and dense cases when $p \leq 30$, achieving  90 \% $\sim$ 100\% feasibility in the sparse case when $p \leq  20$ and in the dense case when $p \leq  30$.
Under the dense demand model, RKHS-DR performs well in low-dimensional cases, but its performance deteriorates as $p$ increases. 
Prescriptive methods---kNN (in the dense case), CART, and RF---maintain a stable yet noticeable gap to SIMOPT, which is possibly due to their limitations to approximate conditional distributions with the sample set that is not sufficiently large.
The relatively small gap of PADR to SIMOPT shows PADR's capability to learn the underlying decision rule structure in problems with nonconvex constraints while maintaining a high probability of feasibility when  $p \leq  30$. 
{NN-DR obtains relatively lower test costs than PADR and RKHS-DR; however, its feasibility rate is substantially lower than those of the other two methods. This suggests that the decision rule yielded by NN-DR may face difficulties in effectively addressing nonconvex-constrained optimization problems.}

\begin{figure}[htbp]
    \FIGURE
    {{
      \subfloat[Basic results]{%
            \includegraphics[width=0.33\textwidth]{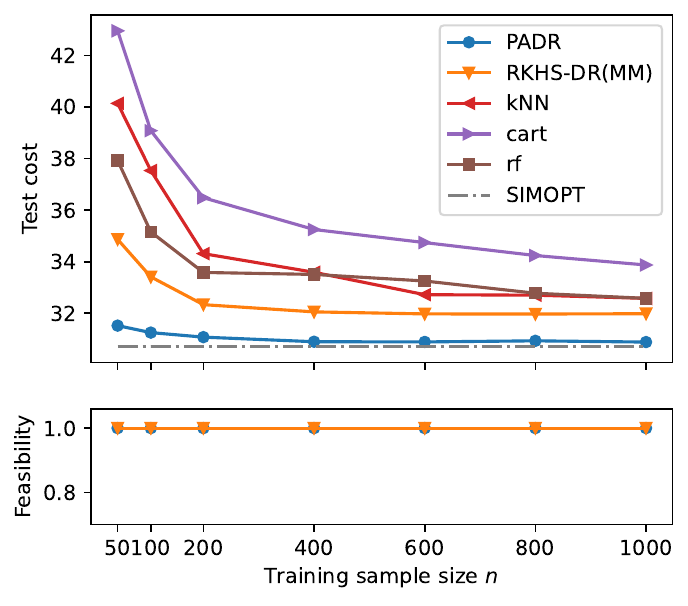}\label{subfig:ncvx-constr-basic}
      }
      \subfloat[Sparse demand model]{%
          \includegraphics[width=0.33\textwidth]{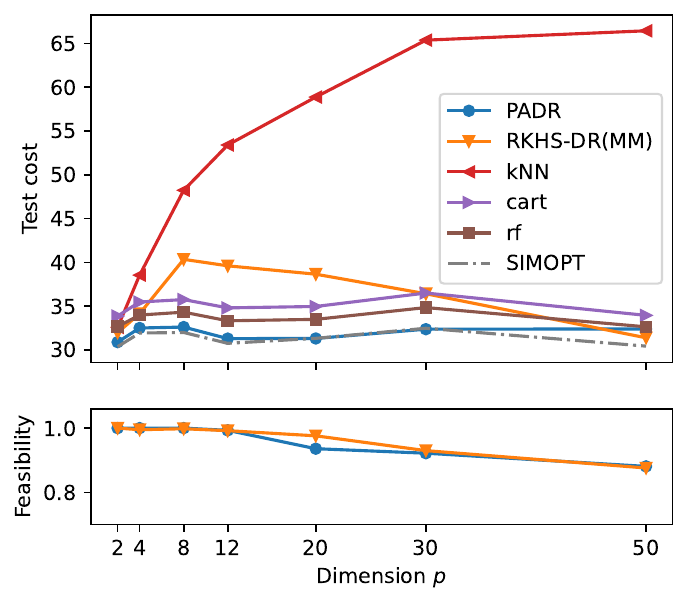}\label{subfig:ncvx-constr-sparse}
      }
      \subfloat[Dense demand model]{%
          \includegraphics[width=0.336\textwidth]{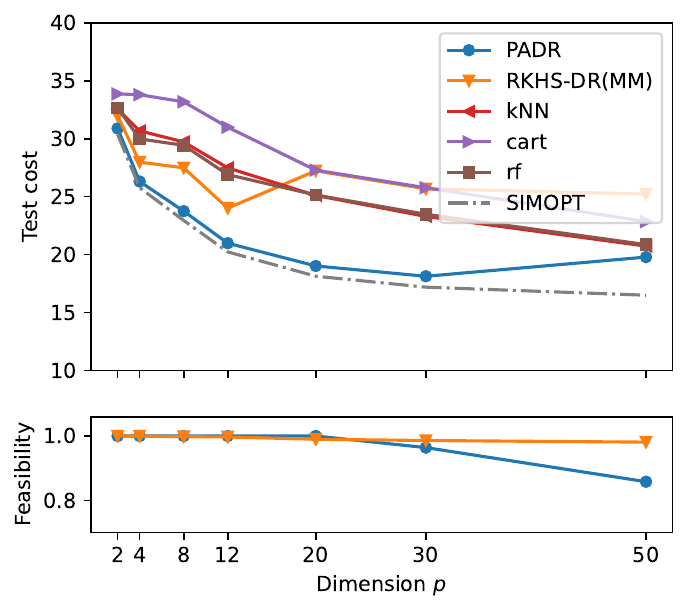}\label{subfig:ncvx-constr-dense}
      }
    }}
    {Newsvendor problem with a nonconvex constraint\label{fig:exp:UNC-NCVX-constr-dim}}
    {%
    }
\end{figure}

{
Thanks to the Associate Editor who suggested a hybrid idea during the review process to address the feasibility issue of DR-based methods, we implement a hybrid strategy by combining the decision rule with an auxiliary method that ensures feasibility (such as RF). Specifically, for the Newsvendor problem with a nonconvex constraint, we combine PADR and NN-DR with the prescriptive method RF, referred to as PADR(hybrid) and NN-DR(hybrid). 
Each method applies the corresponding decision rule first; if the output is infeasible, the RF solution is used as a fallback. We compare the performance of these hybrid methods with pure RF and CART methods, as shown in Figure~\ref{fig:exp:hybrid}. 
Both PADR(hybrid) and NN-DR(hybrid) consistently outperform RF, with the exception when $p = 50$. 
Under the sparse demand model, NN-DR(hybrid) outperforms PADR(hybrid) and shows greater improvements over RF. Under the dense demand model, the performance of both hybrid methods is close to SIMOPT, with PADR(hybrid) slightly better. 
Overall, these preliminary results demonstrate the promise of hybrid approaches in improving both feasibility and performance. This is especially beneficial for contextual SP problems with complex, nonconvex constraints, and we believe this avenue deserves further investigation. 

\begin{figure}[hbtp]
    \FIGURE
    {
      {
      \subfloat[Sparse demand model$\phantom{xxxxxxxx}$]{
          \includegraphics[width=0.46\textwidth]{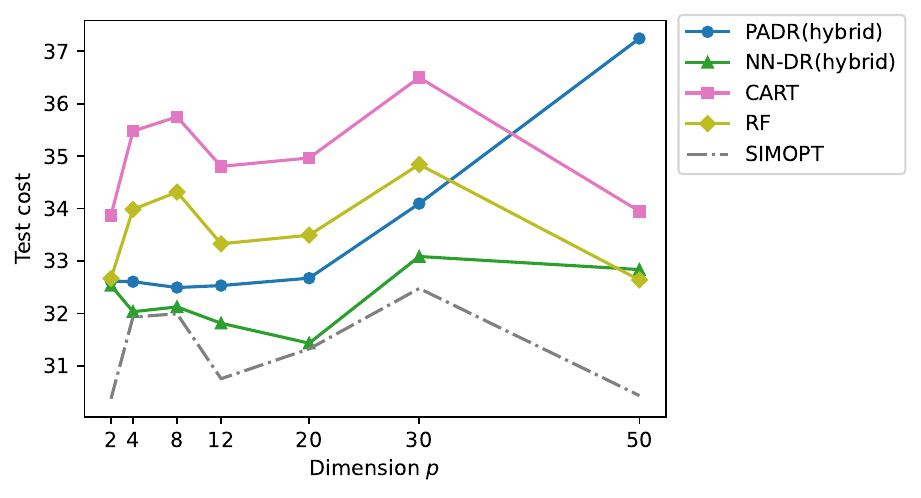}
      }
      \subfloat[Dense demand model$\phantom{xxxxxxx}$]{
          \includegraphics[width=0.46\textwidth]{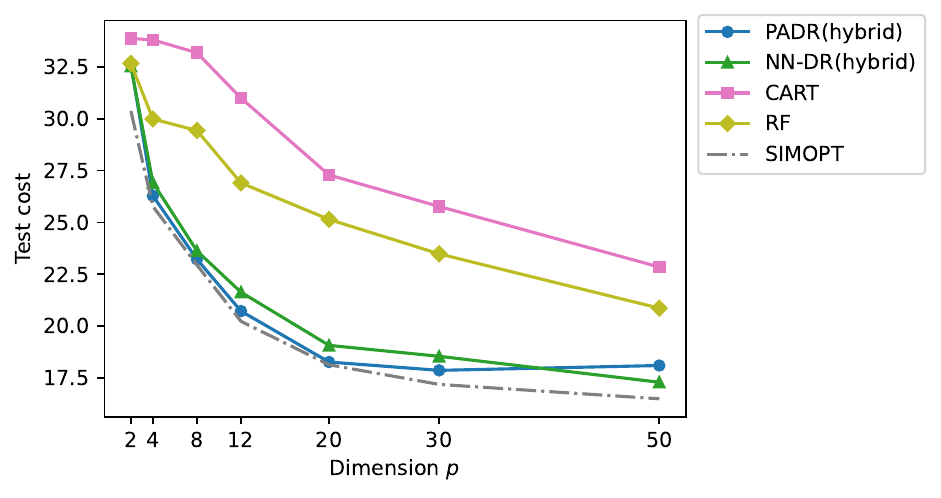}
      }
      }
    }
    {Newsvendor problem with a nonconvex constraint: hybrid approaches\label{fig:exp:hybrid}}
    {}
  \end{figure}

}

\section[Examples Satisfying or Violating Assumption (A4)]%
{Examples Satisfying or Violating Assumption~\ref{assum:A4:opt lipschitz solution f*}}
\label{EC:sec:examples-ensuring-Lipschitz-assumption}

In this section, we present several examples in which the SP problems either satisfy or violate Assumption~\ref{assum:A4:opt lipschitz solution f*}. 
We conclude with a brief discussion on a potential analytical framework to establish sufficient conditions for the Lipschitz continuity of the optimal decision rule.

\subsection[Problems satisfying Assumption (A4)]%
{Problems satisfying Assumption~\ref{assum:A4:opt lipschitz solution f*}}
\label{ECsubsec:problems where assum4 holds}

A canonical class of problems ensuring Assumption~\ref{assum:A4:opt lipschitz solution f*} is quantile regression, such as the Newsvendor problem, which in fact served as the motivation for the Lipschitz continuity assumption in our paper. Additionally, as extensively discussed in \citet[Chapter 6]{shapiro2021lectures}, risk-averse optimization also yields quantile-based solutions. For example, Average Value-at-Risk (AV@R) attains its minimum at the Value-at-Risk (V@R), which by definition is a quantile of the underlying distribution (see Section~6.2.4). 
Section~6.4.2 in the above reference further shows that in the risk-averse Newsvendor problem, where the objective is a composition of a risk measure and the inventory cost, if the risk measure is real-valued and coherent (e.g., AV@R or the mean-upper-semideviation of order $p$ in Section~6.3.2), the optimal order quantity remains a quantile of a distribution that depends on the risk measure and the inventory parameters. In these situations, determining $z^*$ reduces to estimating (or “regressing”) the conditional quantile function based on $X$. Thus, the Lipschitz continuity of $z^*$ hinges on that of the conditional quantile function. 
Under the additive model $Y(x) = \bar{Y}(x) + \epsilon$, where $\epsilon$ is an independent and zero-mean random variable, and $X$ is compact, the continuity of the mean function $\bar{Y}(x)$ implies the desired Lipschitz continuity.

Another class of problems satisfying Assumption \ref{assum:A4:opt lipschitz solution f*} is a two-stage stochastic optimization problem:
\[
\min_{z \in \cZ}\quad g(z) + \EE\left[Q(Az, Y) \mid X = x\right],
\]
where
\[
Q(u, Y) = \min\{q^\top v : Wv = Y - u, v \geq 0\}.
\]
Here, $g: \RR^d \rightarrow \RR$ is convex, $\cZ$ is a closed convex set, and $\mathbb{P}_{Y|X=x}$ is a Borel probability measure with finite first moment. 
Under assumptions of complete recourse and dual feasibility for the second-stage problem, \citet[Theorem~2.2]{romisch1993stability} show that if the function $\EE\left[Q(\cdot, Y) \mid X = x\right]$ is strongly convex for all $x \in \cX$, then the optimal solution set is a singleton, and there exist constants $L > 0$ and $\delta > 0$ such that
\[
\norm{Af^*(x_1) - Af^*(x_2)} \leq L \mathcal{W}_1(\mathbb{P}_{Y\mid x_1}, \mathbb{P}_{Y\mid x_2})^{1/2},
\]
whenever $\mathcal{W}_1(\mathbb{P}_{Y\mid x_1}, \mathbb{P}_{Y\mid x_2}) \leq \delta$, where $\mathcal{W}_1$ denotes the $L_1$-Wasserstein distance. Consequently, if $A$ is invertible, $\cX$ is compact, and $Y(x) = \bar{Y}(x) + \epsilon$ with continuity of $\bar{Y}(\cdot)$, the Lipschitz continuity of the optimal policy follows.
A subsequent work by \citet{romisch1996lipschitz} further investigates the Lipschitz stability of the optimal solution sets with respect to a pseudo-metric over the space of probability measures, which may facilitate the establishment of Lipschitz continuous solution mappings in problems with non-unique solutions.

\subsection[Problems not satisfying Assumption (A4)]%
{Problems not satisfying Assumption~\ref{assum:A4:opt lipschitz solution f*}}

In this part, we present examples where Assumption~\ref{assum:A4:opt lipschitz solution f*} does not hold. Consider the additive random variable model $Y(x) = \bar{Y}(x) + \epsilon$ with a discontinuous  mean function $\bar{Y}(x)$ on $\cX$. Then the resulting optimal decision rule will be discontinuous. A notable instance arises in regression problems with the mean squared error loss as the cost function, where the optimal decision rule coincides with $\bar{Y}(x)$ under normally distributed noise $\epsilon$. 
If $\bar{Y}(x)$ is discontinuous or discrete, e.g., due to constraints that depend on $x$, the corresponding optimal decision rule will also be discontinuous.

For many problems, the set of optimal decisions is inherently disconnected across scenarios. One case is when the feasible region $\cZ$ itself is disconnected, so that the optimal solutions conditional on different features may lie in different connected components. Another case is  when the problem structure yields discrete optimal solutions, such as when $\cZ$ is a polytope and the optimal decision always lies at an extreme point. 
In these situations, small perturbations of the input can cause abrupt jumps in the optimal decision rule, thereby not satisfying Assumption~\ref{assum:A4:opt lipschitz solution f*}. Addressing such cases calls for the development of Heaviside function classes \citep{YueLiuPang25treatment} including decision trees to capture the discontinuous components of the decision rule. Advancing such methods constitutes a promising and important direction for future research.

To conclude this section, we briefly outline a potential direction for analyzing sufficient conditions for the Lipschitz continuity of the optimal decision rule.
A principled approach is to first characterize conditions on the function $\Phi(z, x) \coloneqq \EE_Y [\varphi(z, Y) \mid X = x]$ (defined over $\cZ \times \cX$) that ensure Lipschitz continuity of the solution mapping $z^*(x) \in \arg\min_{z \in \cZ} \Phi(z, x)$, using implicit function theorems for variational inequalities. 
In particular, with the function $\Phi(z, x)$ being continuously differentiable in $z$ and the generalized equation $0 \in \nabla_z \Phi(z, x) + \cN_\cZ(z)$ being sufficient for optimality, \citet{dontchev2009implicit} provide conditions on the mapping $x \mapsto \nabla_z\Phi(\bar{z}, x)$ under which the mapping $z^*(x)$ is locally Lipschitz continuous (e.g., \citet[Theorem~2B.5]{dontchev2009implicit}). 
It thus reduces to identify structural properties of $\varphi(z, Y)$ and conditional distributional properties of $Y$  such that those regularity conditions hold, such  as growth properties of $\Phi$.
While we recognize that verifying these conditions is nontrivial in general, this framework offers a promising direction for theoretically analyzing the regularity of optimal policies in contextual SP, which deserves further exploration.

{
\section{Surrogation Construction for Problems in Section~\ref{section:experiments}}\label{EC:sec:surrogation construction for exps}
In this section we provide the specific surrogation construction of $F(\cdot, \xi^s)$ with $\xi^s \in \Xi$ for problems in Section~\ref{section:experiments}, including the Newsvendor problem, the Newsvendor problem with additional concave cost, and the product placement problem.

For the Newsvendor problem, the cost function $\varphi(z, \xi^s) = c_b\max\{y^s - z, 0\} + c_h\max\{z - y^s, 0\}$ is the sum of two monotonic convex functions. By utilizing surrogation construction for \ref{enum:algo:phi2}, the surrogate function of $F(\theta, \xi^s) = \varphi(f(\theta, \xi^s), \xi^s)$ at $\theta^\nu$ is constructed as follows,
\[
\widehat{F}(\theta, \xi^s; \theta^\nu, I) = c_b\max\left\{\widehat{f}(\theta,\xi^s; I_2)-y^s, 0\right\} + c_h \max\left\{y^s - \widecheck{f}(\theta, \xi^s; I_1), 0\right\},
\]
where $I = (I_1, I_2) \in \cI^{\varepsilon}(\theta^\nu, \xi^s)$ and $\widehat{f}(\theta, \xi^s; I_2)$ and $\widecheck{f}(\theta, \xi^s; I_1)$ are the upper and lower surrogate functions of $f(\theta, \xi^s)$ respectively.

For the Newsvendor problem with the additional concave piecewise-affine cost function,
we denote the concave cost function by the general form $C(z) \coloneqq \min_{i \in [K_C]} \{a_i z + b_i\}$. In this case, the cost function is
\begin{align*}
\varphi(z, \xi^s) 
&= c_b\max\{y^s - z, 0\} + c_h\max\{z - y^s, 0\} + \min_{i \in [K_C]} \{a_i z + b_i\}\\
&= c_b\max\{y^s - z, 0\} + c_h\max\{z - y^s, 0\} - \max_{i \in [K_C]} \{-a_i z - b_i\},
\end{align*}
which satisfies~\ref{enum:algo:phi1}. Following the construction for~\ref{enum:algo:phi1}, the surrogate function of $F(\theta, \xi^s) = \varphi(f(\theta, \xi^s), \xi^s)$ at $\theta^\nu$ is constructed as below,
\begin{align*}
\widehat{F}(\theta, \xi^s; \theta^\nu, I) = & 
c_b\max\left\{\widehat{f}(\theta,\xi^s; I_2)-y^s, 0\right\} + c_h \max\left\{y^s - \widecheck{f}(\theta, \xi^s; I_1), 0\right\}\\
&- \left(\max\{(-a_{I_3}), 0\}^\top \widecheck{f}(\theta, \xi^s; I_1) + \min\{(-a_{I_3}), 0\}^\top \widehat{f}(\theta, \xi^s; I_2) - b_{I_3}\right),
\end{align*}
where $I_1 \in \cI^\varepsilon_h(\theta^\nu, \xi^s)$, $I_2 \in \cI^\varepsilon_g(\theta^\nu, \xi^s)$, and $I_3 \in \cI^\varepsilon_{C}(f(\theta^\nu, \xi^s), \xi^s)$. 
Consistent with the surrogation for~\ref{enum:algo:phi1}, the surrogation index $I = (I_1, I_2, I_3)$ in this case includes the surrogation index for the outer function $C(\cdot)$ at the point $f(\theta^\nu, \xi^s)$.

Finally, for the product placement problem, we first convert the two-stage problem into a single-stage convex problem by extending the auxiliary variable $f$:
\[
\min_{\substack{z \geq 0\\v^s \geq0,\  s\in[n]}}\quad\frac{1}{n}\sum_{s=1}^n 
\left\{c^\top z + g^\top v^s + b^\top \max\{y^s + Wv^s - z, 0\}\right\}.
\]
Since $c \geq 0$ and $b \geq 0$, for the component $\varphi(z, v, \xi^s) = c^\top z + g^\top v + b^\top \max\{y^s + Wv - z, 0\}$, we construct the surrogate function of $F(\theta, v, \xi^s)\coloneqq \varphi(f(\theta,\xi^s),v,\xi^s)$ at $(\theta^\nu)$ as
\[
\widehat{F}(\theta, v, \xi^s; \theta^\nu, I) = c^\top \widehat{f}(\theta, \xi^s, I_2) + g^\top v + b^\top \max\left\{y^s + Wv - \widecheck{f}(\theta, \xi^s; I_1), 0\right\},
\]
where $I = (I_1, I_2) \in \cI^{\varepsilon}(\theta^\nu, \xi^s)$.
}

{\section{Additional Experimental Results}\label{EC:sec:experiments}
In this section, we further conduct experiments on the Newsvendor problem using alternative demand models that share a fixed linear component but with different nonlinear structures, including quadratic, cubic, piecewise-affine and sine functions. 
These experiments are designed to assess the robustness of PADR and benchmark methods across varying nonlinear models.
By incorporating an explicit linear term in all models, we aim to examine whether the nonlinear effects still manifest when combined with a linear component.

Specifically, we consider the following demand models:
\begin{align}
&\bar{Y}_{\text{ec1}}(X) = 3(X_{1}+X_{2}) + \nonlinearity \cdot \max\{5X_{1}-10X_{2}, -10X_{1}+5X_{2}, 15X_{1}\} + 20,\label{eq:demand-model:pa}
\\
&\bar{Y}_{\text{ec2}}(X) = 3(X_{1}+X_{2}) + \nonlinearity \cdot 5(X_{1}-X_{2})^2 + 10,\label{eq:demand-model:quad}
\\
&\bar{Y}_{\text{ec3}}(X) = 3(X_{1}+X_{2}) + \nonlinearity \cdot [5(X_{1}-X_{2})^{3} - 15(X_{1}-X_{2})] + 50,\label{eq:demand-model:cube}
\\
&\bar{Y}_{\text{ec4}}(X) = 3(X_{1}+X_{2}) + \nonlinearity \cdot [4\sin(\pi X_{1}) + \max\{16X_{2}, -20X_{2}\}] + 30.\label{eq:demand-model:psin}
\end{align}

We report the results under the above models in Figures~\ref{fig:exp:nv-lnrMA3dense-k}-\ref{fig:exp:nv-lnrSinPAdense-k}. 
\begin{figure}[htbp]
  \FIGURE
  {
    {
      \subfloat[$p=2 \phantom{xxxxxxx}$]{%
        \includegraphics[width=0.426\textwidth]{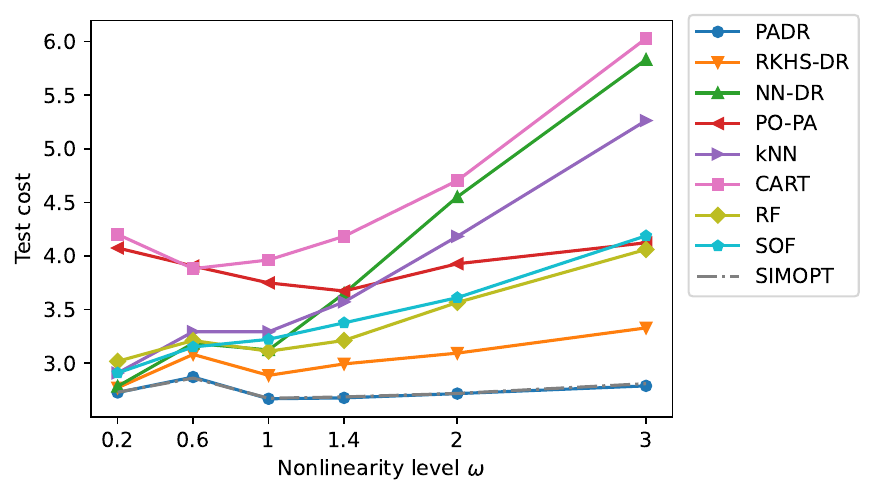}\label{fig:exp:nv-lnrMA3dense-k-dim2}
      }
      \subfloat[$p=20$  $\phantom{xxxxxx}$]{%
        \includegraphics[width=0.426\textwidth]{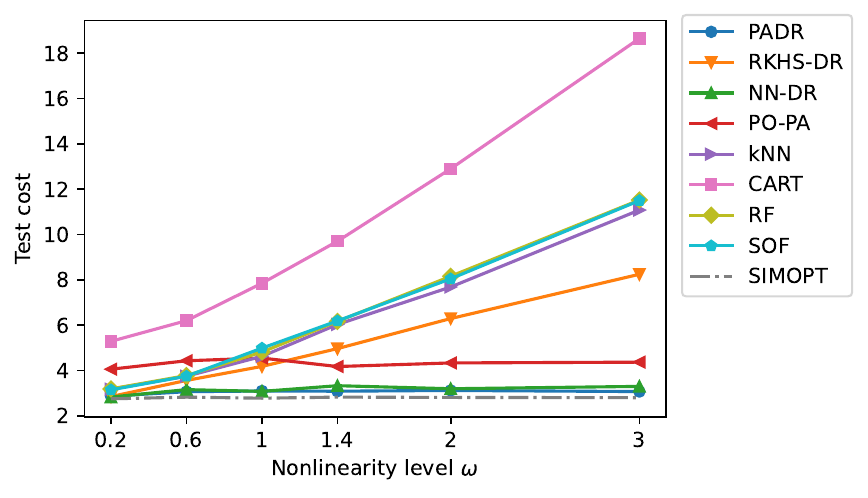}\label{fig:exp:nv-lnrMA3dense-k-dim20}
      }
    }
  }
  {Newsvendor problem with varying nonlinearity in demand model~\eqref{eq:demand-model:pa}\label{fig:exp:nv-lnrMA3dense-k}}
  {}
\end{figure}

\begin{figure}[htbp]
  \FIGURE
  {
    {
      \subfloat[$p=2 \phantom{xxxxxxx}$]{%
        \includegraphics[width=0.426\textwidth]{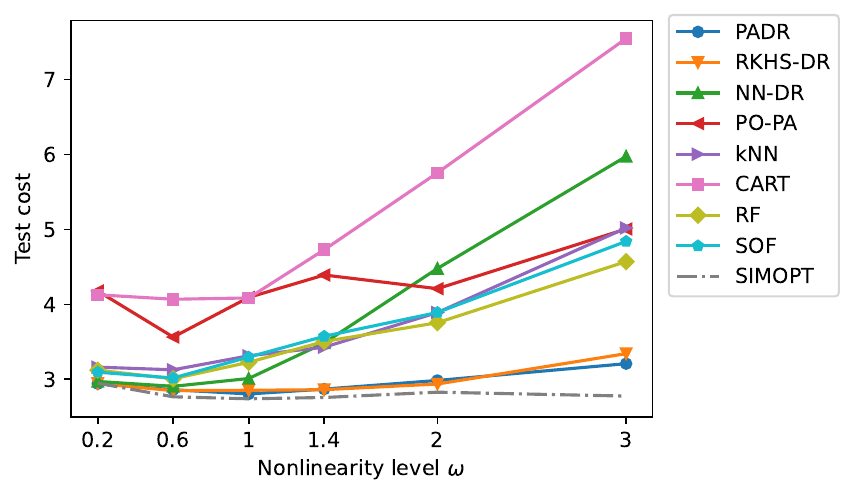}\label{fig:exp:nv-lnrX2dense-k-dim2}
      }
      \subfloat[$p=20$  $\phantom{xxxxxx}$]{%
        \includegraphics[width=0.426\textwidth]{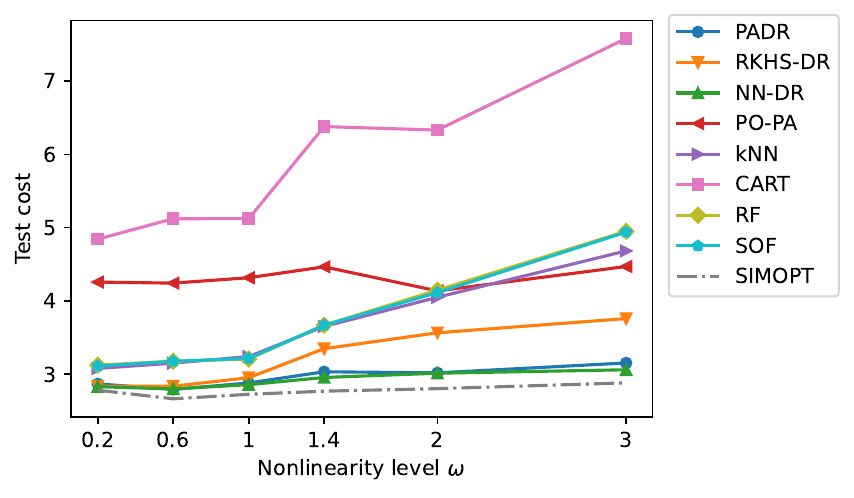}\label{fig:exp:nv-lnrX2dense-k-dim20}
      }
    }
  }
  {Newsvendor problem with varying nonlinearity in demand model~\eqref{eq:demand-model:quad}\label{fig:exp:nv-lnrX2dense-k}}
  {}
\end{figure}

\begin{figure}[htbp]
  \FIGURE
  {
    {
      \subfloat[$p=2 \phantom{xxxxxxx}$]{%
        \includegraphics[width=0.426\textwidth]{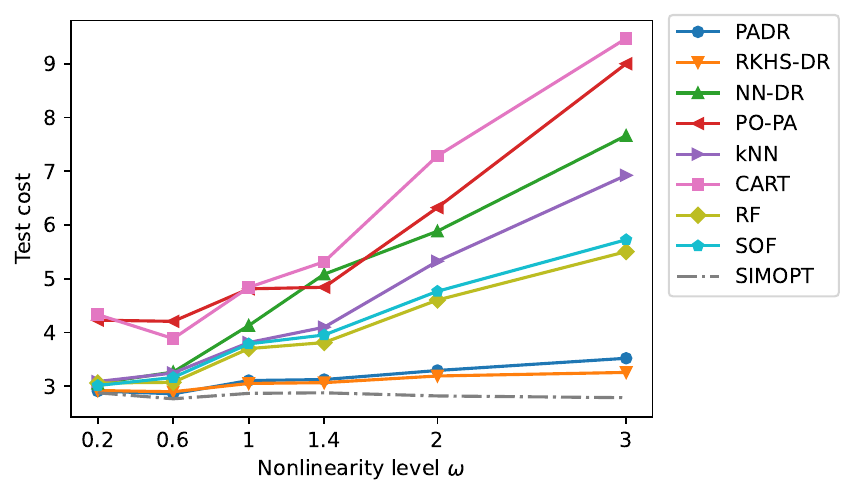}\label{fig:exp:nv-lnrX3dense-k-dim2}
      }
      \subfloat[$p=20$  $\phantom{xxxxxx}$]{%
        \includegraphics[width=0.426\textwidth]{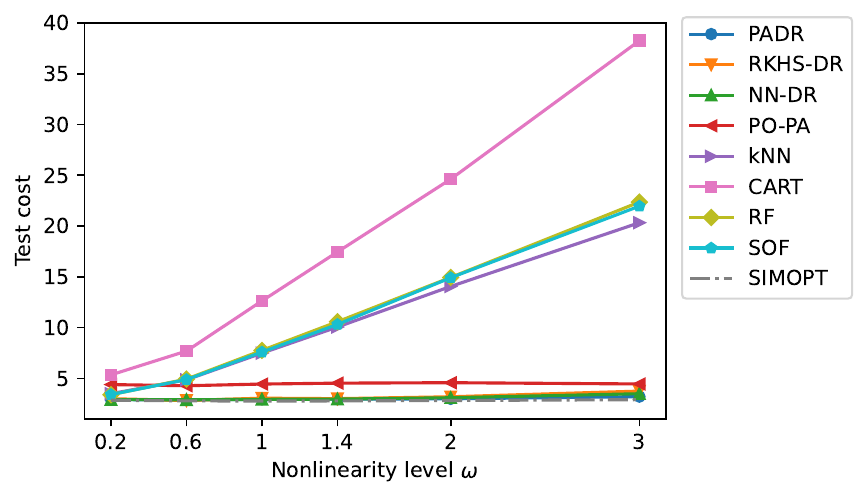}\label{fig:exp:nv-lnrX3dense-k-dim20}
      }
    }
  }
  {Newsvendor problem with varying nonlinearity in demand model~\eqref{eq:demand-model:cube}\label{fig:exp:nv-lnrX3dense-k}}
  {}
\end{figure}

\begin{figure}[htbp]
  \FIGURE
  {
    {
      \subfloat[$p=2 \phantom{xxxxxxx}$]{%
        \includegraphics[width=0.426\textwidth]{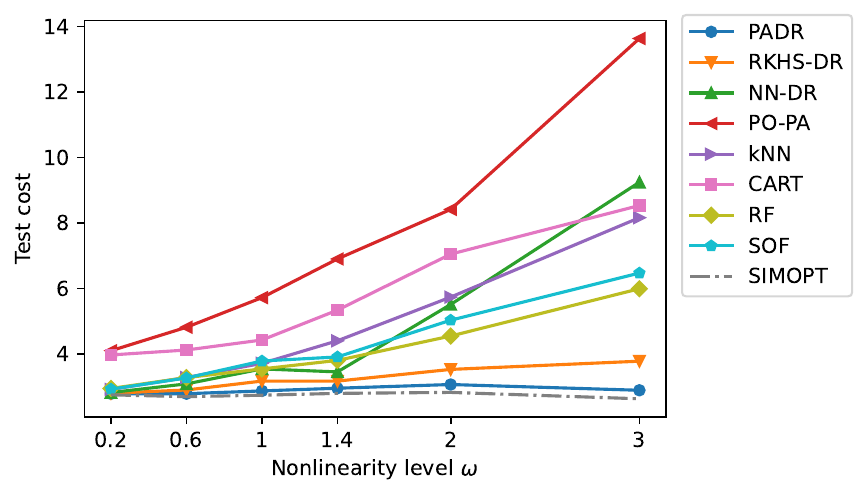}\label{fig:exp:nv-lnrSinPAdense-k-dim2}
      }
      \subfloat[$p=20$  $\phantom{xxxxxx}$]{%
        \includegraphics[width=0.426\textwidth]{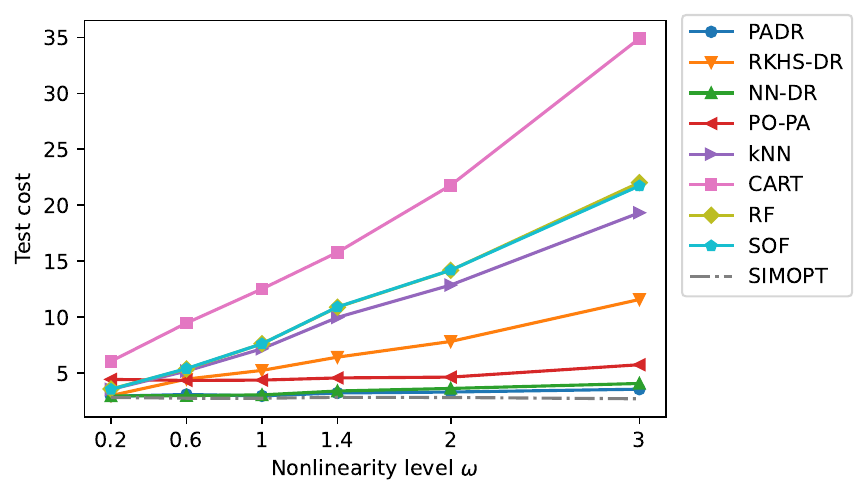}\label{fig:exp:nv-lnrSinPAdense-k-dim20}
      }
    }
  }
  {Newsvendor problem with varying nonlinearity in demand model~\eqref{eq:demand-model:psin}\label{fig:exp:nv-lnrSinPAdense-k}}
  {}
\end{figure}
In general, the comparison results remain similar to those observed in Section~\ref{subsubsec:nonlinearity}.
Despite the inclusion of an explicit linear component in each demand model, the nonlinear terms continue to significantly affect the performance of benchmark methods. Across all settings, PADR consistently ranks among the top-performing methods. The observed degradation in NN-DR's performance in low-dimensional cases persists. In contrast, RKHS-DR is competitive at $p=2$, but its performance deteriorates as the feature dimension increases, likely due to the higher sample complexity required for kernel-based methods to generalize in high-dimensional spaces. 
Furthermore, RKHS-DR performs noticeably worse under demand models with nonsmooth components (i.e., models~\eqref{eq:demand-model:pa} and~\eqref{eq:demand-model:psin}) than under smooth ones (i.e., models~\eqref{eq:demand-model:quad} and~\eqref{eq:demand-model:cube}), suggesting that kernel-based smooth approximations struggle to capture nonsmooth structures without sufficient data.
These results confirm the robustness of PADR to various nonlinear demand models with superior performance, compared to all other benchmarks, including NN-DR and RKHS-DR.

}

\end{APPENDICES}

\ACKNOWLEDGMENT{The second author gratefully acknowledge the support by National Natural Science Foundation of China under grant number 72201151, and 72188101. }

\end{document}